\documentclass[11pt]{amsart}
\usepackage{mathrsfs, amsmath, amssymb, array,graphics, booktabs,diagbox,listings,color,textcomp,picinpar,wrapfig,cutwin,longtable,slashbox,amsthm,tikz,threeparttable,makecell,cases}

\usepackage[backref]{hyperref}
\usepackage{multicol}
\usepackage{multirow}
\usepackage[page,title,titletoc,header]{appendix}
\usepackage[font=small,labelfont=bf]{caption}
\DeclareCaptionFont{tiny}{\tiny}
\DeclareCaptionFont{normalsize}{\normalsize}

\usepackage{geometry}
\definecolor{listinggray}{gray}{0.9}

\definecolor{lbcolor}{rgb}{0.9,0.9,0.9}

\usepackage{algorithm}
\usepackage[noend]{algpseudocode}

\newtheorem{theorem}{Theorem}[section]

\newtheorem{corollary}[theorem]{Corollary}
\newtheorem{proposition}[theorem]{Proposition}

\theoremstyle{definition}
\newtheorem{definition}[theorem]{Definition}

\theoremstyle{remark}
\newtheorem{remark}[theorem]{Remark}

\numberwithin{equation}{section}

\newcommand{\HH}{\mathbb{H}}

\begin{document}

\title{Compact hyperbolic Coxeter four-dimensional polytopes with eight facets}

\author{Jiming Ma}
\address{School of Mathematical Sciences \\Fudan University\\Shanghai 200433, China} \email{majiming@fudan.edu.cn}

\author{Fangting Zheng}
\address{Department of Mathematical Sciences \\ Xi'an Jiaotong Liverpool University\\ Suzhou 200433,	China  }
\email{Fangting.Zheng@xjtlu.edu.cn}

\date{Nov. 22, 2022}
\thanks{Jiming Ma was partially supported by  NSFC  11771088  and 12171092. Fangting Zheng was supported by NSFC 12101504 and XJTLU Research Development Fund RDF-19-01-29}

\begin{abstract}
	In this paper, we obtain the complete classification for compact hyperbolic Coxeter four-dimensional polytopes with eight facets.
\end{abstract}

\maketitle

\section{Introduction}

A Coxeter polytope in the spherical, hyperbolic or Euclidean space is a polytope whose dihedral angles are all integer sub-multiples of $\pi$. Let $\mathbb{X}^d$ be $\mathbb{E}^d$, $\mathbb{S}^d$, or $\mathbb{H}^d$. If $\Gamma\subset Isom(\mathbb{X}^d)$ is a finitely generated discrete reflection group, then its fundamental domain is a Coxeter polytope in $\mathbb{X}^d$. On the other hand, if $\Gamma=\Gamma(P)$ is generated by reflections in the bounding hyperplanes of a Coxeter polytope $P\subset\mathbb{X}^d$, then $\Gamma$ is a discrete group of isometries of $\mathbb{X}^d$ and $P$ is its fundamental domain.

There is an extensive body of literature in this field. In early work, \cite{Coxeter:1934} has proved that any spherical Coxeter polytope 
is a simplex and any Euclidean Coxeter polytope is either a simplex or a direct product of simplices. See, for example, \cite{Coxeter:1934,Bourbaki: 1968}  for full lists of spherical and Euclidean Coxeter polytopes.

However, for hyperbolic Coxeter polytopes, the classification remains an active research topic. It was proved by Vinberg \cite{Vinberg:1985} that no compact hyperbolic Coxeter polytope exists in dimensions $d\geq 30$, and non-compact hyperbolic Coxeter polytope of finite volume does not exist in dimensions $d \geq 996$ \cite{Prokhorov:1987}. These bounds, however, may not be sharp. Examples of compact polytopes are known up to dimension $8$ \cite{Bugaenko:1984,Bugaenko:1992}; non-compact polytopes of finite volume are known up to dimension $21$ \cite{Vinberg:1972,VK:1978,Borcherds: 1998}. As for the classification, complete results are only available in dimensions less than or equal to three. Poincare finished the classification of 2-dimensional hyperbolic polytopes in \cite{Poincare: 1882}. That result was important to the work of Klein and Poincare on discrete groups of isometries of the hyperbolic plane. In 1970, Andreev proved an analogue for the $3$-dimensional hyperbolic convex finite volume polytopes \cite{Andreev1: 1970,Andreev2: 1970}. This theorem played a fundamental role in Thurston's work on the  geometrization of 3-dimensional Haken manifolds.

In higher dimensions, although a complete classification is not available, interesting examples have been displayed in \cite{Makarov: 1965,Makarov: 1966,Vinberg: 1967,Makarov: 1968,Vinberg: 1969,Rusmanov: 1989,ImH: 1990,Allcock: 2006}. In addition, enumerations are reported for the cases in which the differences between the numbers of facets $m$ and the  dimensions $d$ of polytopes are fixed to some small number. When $m-d=1$, Lann\'{e}r classified all compact hyperbolic Coxeter simplices \cite{Lanner: 1950}. The enumeration of non-compact hyperbolic simplices with finite volume has been reported by several authors, see e.g. \cite{Bourbaki: 1968,Vinberg: 1967,Koszul:1967}. For $m-d=2$, Kaplinskaja described all compact or non-compact but of finite volume hyperbolic Coxeter simplicial prisms \cite{Kaplinskaya: 1974}. Esselmann \cite{Esselmann: 1996} later enumerated other compact possibilities in this family, which are named \emph{Esselmann polytopes}. 
Tumarkin \cite{Tumarkin: n2} classified all other non-compact but of finite volume hyperbolic Coxeter $d$-dimensional polytopes with
$n + 2$ facets. In the case of $m-d=3$, Esselman proved in 1994 that compact hyperbolic Coxeter $d$-polytopes with $d+3$ facets  only exist when $d\leq 8$ \cite{Esselmann}. By expanding the techniques derived by Esselmann in \cite{Esselmann} and \cite{Esselmann: 1996}, Tumarkin completed the classification of compact hyperbolic Coxeter $d$-polytopes with $d+3$ facets \cite{Tumarkin: n3}. In the non-compact case, Tumarkin proved in \cite{Tumarkin: n3fv,Tumarkin: n3fvs} that such polytopes do not exist in dimensions greater than or equal to $17$; there is a unique polytope in dimension of $16$. Moreover, the author provided in the same papers the complete classification of  a special family of  pyramids over a product of three simplices, which exist only in dimension of $4,5,\cdots,9$ and $13$. The classification for the case of finite volume has not completed yet. Regarding this sub-family, Roberts provided a list with exactly one non-simple vertex \cite{Roberts:15}. In the case of $m-d=4$,
Flikson-Tumarkin showed in \cite{FT:08} that compact hyperbolic Coxeter $d$-polytope with $d + 4$ facets does not exist when $d$ is greater than or equal to $8$. This bound is sharp because of the example constructed by Bugeanko \cite{Bugaenko:1984}. In addition, Flikson-Tumarkin showed that Bugeanko's example is the only $7$-dimensional polytope with $11$ facets. However, complete classifications for $d=4,5,6$ are not presented there. 

Besides, some scholars have also considered polytopes with small numbers of disjoint pairs
\cite{FT:08s,FT:09,FT:14} or of  certain combinatoric types, such as $d$-pyramid \cite{Tumarkin: n2,Tumarkin: n3fv} and $d$-cube \cite{Jacquemet2017,JT:2018}. An updating overview of the current knowledge for hyperbolic Coxeter polytopes is available on Anna Felikson's webpage \cite{Annahomepage}.

In this paper, we classify all the compact hyperbolic Coxeter $4$-polytope with $8$ facet. The main theorem is as follows:

\begin{theorem}\label{thm:main}
	There are exactly $348$ compact hyperbolic Coxeter $4$-polytopes with $8$ facets. In particular, $P_{21}$ has two dihedral angles of $\displaystyle\frac{\pi}{12}$, and $P_{17,8}$ has an dihedral angle of $\displaystyle\frac{\pi}{7}$. Among hyperbolic Coxeter polytopes of dimensions larger than or equal to $4$, these two values of dihedral angles appear for the first time and $\displaystyle\frac{\pi}{12}$ is the smallest dihedral angle known so far.
\end{theorem}

We remark that Burcroff has also carried out the same list independently almost in the same time \cite{Amanda:2022}. Mutually comparing the results benefit both authors. Burcroff communicated with us when our preprint appear. She kindly pointed out several typos about conveying the information into Coxeter diagrams. Ours also help she find out two lost or double-counted combinatoric types that admit hyperbolic structure. We now all agree that $348$ is the correct number. The correspondence between the notions of our and Burcroff's polytopes is presented in Section \ref{section:vadilation}.

The paper \cite{JT:2018} is the main inspiration for our recent work on enumerating hyperbolic Coxeter polytopes. In comparison with \cite{JT:2018}, we use a more universal ``block-pasting" algorithm, which is first introduced in \cite{MZ:2018}, rather than the ``tracing back" attempt. More geometric obstructions are adopted and programmed to considerably reduce the computational load. Our algorithm efficiently enumerates hyperbolic Coxeter polytopes over arbitrary combinatoric type rather than merely the $n$-cube. 

Last but not the least, our main motivation in studying the hyperbolic Coxeter polytopes is for  the construction of high-dimensional hyperbolic manifolds. However, this is not the theme here. Readers can refer to, for example,  \cite{KM: 2013}, for interesting hyperbolic manifolds built via special hyperbolic Coxeter polytopes.

 The paper is organized as follows. We provide in Section $2$ some preliminaries about (compact) hyperbolic (Coxeter) polytopes. In Section $3$, we recall the 2-phases procedure and related terminologies introduced by Jacquemet and Tschantz \cite{Jacquemet2017,JT:2018} for numerating all hyperbolic Coxeter $n$-cubes. The $37$ combinatorial types of simple $4$-polytopes with $8$ facets are reported in Section $4$. Enumeration of all the ``SEILper"-potential matrices are explained in Section $5$. The ``6-rounds" procedure are applied to the ``SEILper"-potential matrices for the Gram matrices of actual hyperbolic Coxeter polytopes in Section $6$. Validations and the complete lists of the resulting Coxeter diagrams of the Theorem \ref{thm:main} are shown in Section $7$. 

\textbf{Acknowledgment}

We would like to thank Amanda Burcroff a lot for communicating with us about her result and pointing out several confusing drawing typos in the first arXiv version. The computations is pretty delicate and complex, and the list now is much more  convincing due to the  mutual check. We are also grateful to Nikolay Bogachev for his interest and discussion about the results, and noting the missing of a hyperparallel distance data and some textual mistakes in previous version. The computations throughout this paper are performed on a cluster of server of PARATERA, engrid12, line priv$\_$para (CPU:Intel(R) Xeon(R) Gold 5218 16 Core v5@2.3GHz).

\section{preliminary} \label{section:cchp}
In this section, we recall some essential facts about compact Coxeter hyperbolic polytopes, including Gram matrices, Coxeter diagrams, characterization theorems, etc. Readers can refer to, for example, \cite{Vinberg:1993} for more details.  

\subsection{Hyperbolic space, hyperplane and convex polytope}
We first describe a hyperboloid model of the $d$-dimensional hyperbolic space $\mathbb{H}^d$. Let $\mathbb{E}^{d,1}$ be a $d+1$-dimensional Euclidean vector space equipped with a Lorentzian scalar product $\langle\cdot,\cdot \rangle$ of signature $(d,1)$. We denote by $C_+$ and $C_-$ the connected components of the open cone $$C=\{x=(x_1, ...,x_d,x_{d+1})\in \mathbb{E}^{d,1}:\langle x,x\rangle<0\}$$ 
with $x_{d+1}>0$ and $x_{d+1}<0$, respectively. Let $R_{+}$ be the group of positive numbers acting on $\mathbb{E}^{d,1}$ by homothety. The hyperbolic space $\mathbb{H}^d$ can be identified with the quotient set $C_+/R_+$, which is a subset of $P\mathbb{S}^d=(\mathbb{E}^{d,1}\backslash \{0\})/R_+.$ There is a natural projection  $$\pi:(\mathbb{E}^{d,1}\backslash \{0\})\rightarrow P\mathbb{S}^d.$$ 

\noindent We denote $\overline{\mathbb{H}^d}$ as the completion of $\mathbb{H}^d$ in $P\mathbb{S}^d$. The points of the boundary $\partial \mathbb{H}^d= \overline{\mathbb{H}^d}\backslash \mathbb{H}^d$ are called the \emph{ideal points}. 
The affine subspaces of $\mathbb{H}^d$ of dimension $d-1$ are \emph{hyperplanes}. In particular, every hyperplane of $\mathbb{H}^d$ can be represented as $$H_e=\{\pi(x):x\in C_+,\langle x,e\rangle=0\},$$ where $e$ is a vector with $\langle e,e \rangle=1$. The half-spaces separated by $H_e$ are denoted by $H_e^+$ and $H_e^-$, where
\begin{equation}
H_e^-=\{\pi(x):x\in C_+,\langle x,e\rangle\leq 0\}. \label{1}
\end{equation} 

The \emph{mutual disposition of hyperplanes} $H_{e}$ and $H_{f}$ can be described in terms of the corresponding two vectors $e$ and $f$ as follows: 
\begin{itemize}
	\item The hyperplanes $H_{e}$ and $H_{f}$ intersect if $\vert\langle e,f \rangle\vert<1$. The value of the dihedral angle of $H_{e}^-\cap H_{f}^-$, denoted by $\angle H_e H_f$, can be obtained via the formula $$\cos \angle H_e H_f=-\langle e,f\rangle;$$
	\item The hyperplanes $H_{e}$ and $H_{f}$ are ultra-parallel if  $\vert\langle e,f\rangle\vert=1$;
	\item The hyperplanes $H_{e}$ and $H_{f}$ diverge if $\vert\langle e,f\rangle\vert>1$. The distance $\rho(H_e,H_f)$ between $H_e$ and $H_f$, when $H_e^+\subset H_f^-$ and $H_f^+\subset H_e^-$, is determined by  $$\cosh \rho(H_e,H_f)=-\langle e,f \rangle.$$
\end{itemize}


We say a hyperplane $H_e$ \emph{supports} a closed bounded convex set $S$ if $H_e\cap S \ne \empty 0$ and $S$
lies in one of the two closed half-spaces bounded by $H_e$. If a hyperplane $H_e$ supports $S$, then $H_e \cap S$ is called a \emph{face} of $S$. 

\begin{definition}
	A $d$-dimensional convex hyperbolic polytope is a subset of the form 
	
	\begin{equation}
	P=\overline{\mathop{\cap}\limits_{i\in \mathcal{I}} H_{i}^-}\in\overline{\mathbb{H}^d}, \label{2}
	\end{equation}
	where $H_i^-$ is the negative half-space bounded by the hyperplane $H_i$ in $\mathbb{H}^d$ and the line ``--- "  above the intersection  means taking the completion in $\overline{\mathbb{H}^d}$, 
	under the following assumptions:
	\begin{itemize}
		\item  $P$ contains a non-empty open subset of $\mathbb{H}^d$ and is of finite volume;
		\item  Every bounded subset $S$ of $P$ intersects only finitely many $H_{i}$.
	\end{itemize}
\end{definition}

A convex polytope of the form (\ref{2}) is called \emph{acute-angled} if for distinct $i,j$, either $\angle H_iH_j\leq \frac{\pi}{2}$ or $H_i^+\cap H_j^+=\emptyset$. It is obvious that Coxeter polytopes are acute-angled. We denote  $e_i$ as the corresponding unit vector of $H_i$, namely $e_i$ is orthogonal to $H_i$ and point away from $P$. The polytope $P$ has the following form in the hyperboloid model.
$$P=\pi(K)\cap \overline{\mathbb{H}^d},$$ where $K=K(P)$ is the convex polyhedral cone in $\mathbb{E}^{d,1}$ given by $$K=\{x\in \mathbb{E}^{d,1}: \langle x,e_i\rangle \leq 0~ \text{for all}~ i \}.$$ 

In the sequel, a $d$-dimensional convex polytope $P$ is called a \emph{$d$-polytope}. A $j$-dimensional face is named a $j$-face of $P$. In particular, a $(d-1)$-face is called a \emph{facet} of $P$.  We assume that each of the hyperplane $H_i$ intersects with $P$ on its facet. In other words, the hyperplane $H_i$ is uniquely determined by $P$ and is called a \emph{bounding hyperplane} of the polytope $P$. A hyperbolic polytope $P$ is called \emph{compact} if all of its $0$-faces, i.e., vertices, are in $\mathbb{H}^d$. It is called of  \emph{finite volume} if some vertices of $P$ lie in $\partial\mathbb{H}^d$.

\subsection{Gram matrices, Perron-Frobenius Theorem, and Coxeter diagrams}\label{hcp}

Most of the content in this subsection is well-known by  peers in this field. We present them for the convenience of the readers. 
In particular, Theorem \ref{thm:signature} and \ref{Vinberg:thm3.1} are extremely important throughout this paper. 

For a hyperbolic Coxeter $d$-polytope $P=\overline{\mathop{\cap}\limits_{i\in \mathcal{I}} H_{i}^-}$, we define the Gram matrix of polytope $P$ to be the Gram matrix $(\langle e_i,e_j\rangle)$ of the system of vectors $\{e_i\in \mathbb{E}^{d,1}\vert i\in\mathcal{I}\}$ that determine $H_i^-$s. 
The Gram matrix of $P$ is the $m\times m$ symmetric matrix $G(P)=(g_{ij})_{1\leq i,j\leq m}$  defined as follows:

\begin{center}
	$g_{ij}=
	\left\{
	\begin{array}{ccl}1  \hspace{0.7cm} &\mbox~~~~~{{\rm if}}& j=i, \\
	-\cos\frac{\pi}{k_{ij}}&\mbox ~~~~~{{\rm if}}&  H_i~ \text{and}~ H_j~ \text{intersect} \text{ at~a~dihedral~ angle~}~ \frac{\pi}{k_{ij}},\\
	-\cosh \rho_{ij}&\mbox ~~~~~{{\rm if}}&  H_i~ \text{and}~ H_j~ \text{divergeat~at~a~distance}~\rho_{ij},\\
	-1  \hspace{0.7cm} &\mbox~~~~~{{\rm if}}& H_i~\text{and}~H_j~\text{are~ultra-parallel}.
	\end{array}\right.
	$
\end{center}

Other than the Gram matrix,  a Coxeter polytope $P$ can also be described by its \emph{Coxeter graph} $\Gamma=\Gamma(P)$. Every node $i$ in $\Gamma$ represents the bounding hyperplane $H_i$ of $P$. Two nodes $i_1$ and $i_2$ are joined by an edge with weight $2\leq k_{ij} \leq \infty$ if $H_i$ and $H_j$ intersect in $\mathbb{H}^n$ with angle $\frac{\pi}{k_{ij}}$. If the hyperplanes $H_i$ and $H_j$ have a common perpendicular of length $\rho_{ij}>0$ in $\mathbb{H}^n$, the nodes $i_1$ and $i_2$ are joined by a dotted edge, sometimes labelled $\cosh \rho_{ij}$. In the following, an edge of weight $2$ is omitted, and an edge of weight $3$ is written without its weight. The rank of $\Gamma$ is defined as the number of its nodes. In the compact case, $k_{ij}$ is not $\infty$, and we have $2\leq k_{ij}< \infty$. 

A square matrix $M$ is said to be the direct sum of the matrices $M_1, M_2,\cdots,M_n$ if by some permutation of the rows and of columns, it can be brought to the form
\begin{center}
	$ \begin{pmatrix}
	M_1&&&&0 \\
	&M_2&&&\\
	&& \cdot&&\\
	&&&\cdot& \\
	0&&&&M_n
	\end{pmatrix}_. $
	
\end{center}

\noindent A matrix $M$ that cannot be represented as a direct sum of two matrices is said to be \emph{indecomposible}\footnote{It is also referred to as ``irreducible" in some references.}. Every matrix can be represented uniquely as a direct sum of indecomposible matrices, which are called (indecomposible) components. 
We say a polytope is \emph{indecomposible} if its Gram matrix $G(P)$ is indecomposible.

\begin{figure}[h]
	\scalebox{0.35}[0.35]{\includegraphics {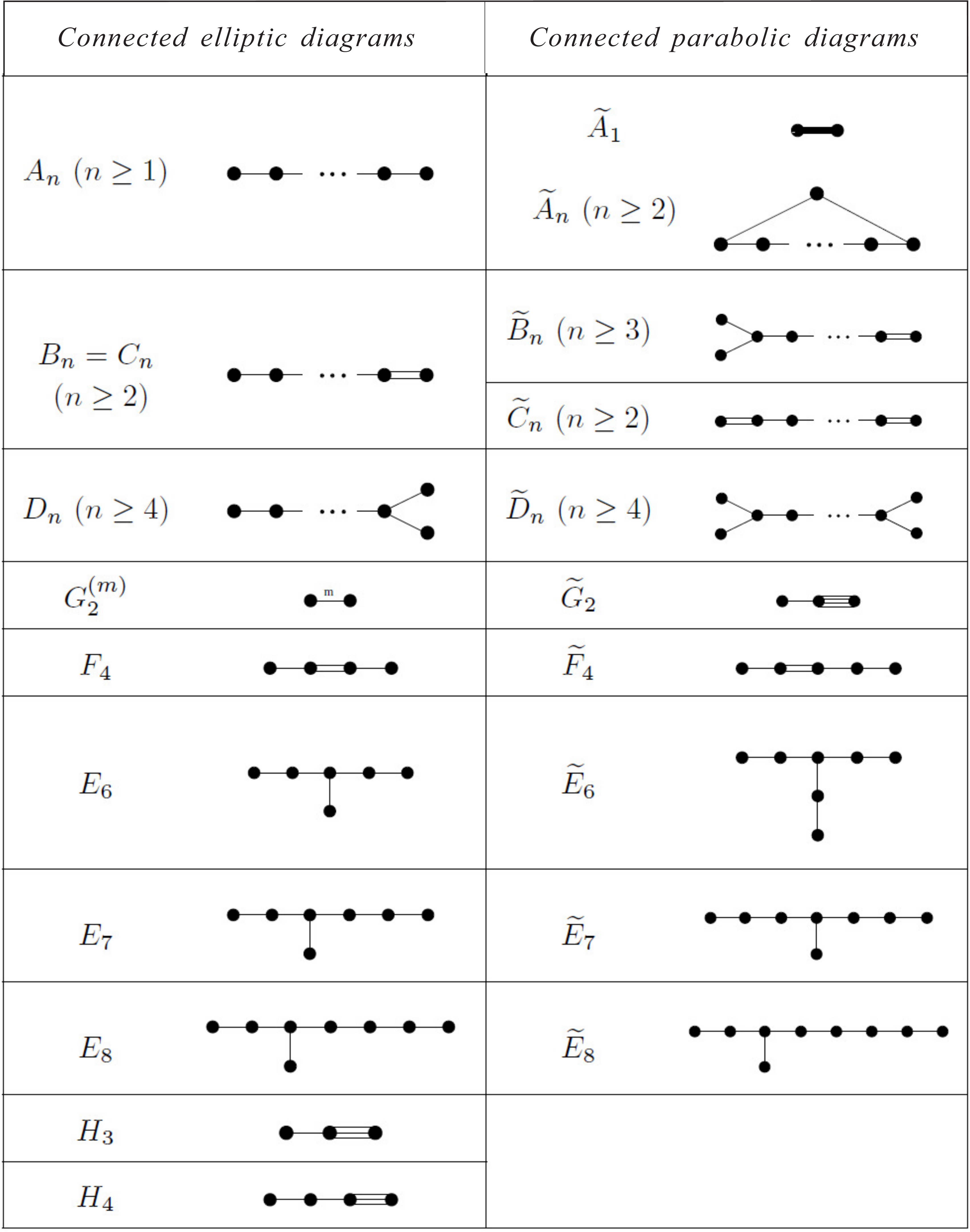}}
	\caption{Connected elliptic (left) and connected parabolic (right) Coxeter diagrams.}
	\label{figure:coxeter}
\end{figure}

In 1907, Perron found a remarkable property of the the eigenvalues and eigenvectors of matrices with positive entries. Frobenius later generalized it by investigating the spectral properties of indecomposible non-negative matrices.

\begin{theorem}[Perron-Frobenius, \cite{G:1959}]\label{thm:PF}
	An indecomposible matrix $A=(a_{ij})$ with non-positive entries always has a single positive eigenvalue $r$ of $A$. The corresponding eigenvector has positive coordinates. The moduli of all of the other eigenvalues do not exceed $r$.
\end{theorem}

It is obvious that Gram matrices $G(P)$ of an indecomposible Coxeter polytope is an indecomposible symmetric matrix with non-positive entries off the diagonal. Since the diagonal elements of $G(P)$ are all $1$s, $G(P)$ is either positive definite, semi-positive definite or indefinite. According to the Perron-Frobenius theorem, the defect of a connected semi-positive definite matrix $G(P)$ does not exceed $1$, and any proper submatrix of it is positive definite. For a Coxeter $n$-polytope $P$, its Coxeter diagram $\Gamma(P)$ is said to be \emph{elliptic} if $G(P)$ is positive definite; $\Gamma (P)$ is called \emph{parabolic} if any indecomposable component of $G(P)$ is degenerate and every subdiagram is elliptic. The elliptic and connected parabolic diagrams are exactly the Coxeter diagrams of spherical and Euclidean Coxeter simplices, respectively. They are classified  by Coxeter \cite{Coxeter:1934} as shown in Figure \ref{figure:coxeter}.

\begin{figure}[h]
	\scalebox{0.35}[0.35]{\includegraphics {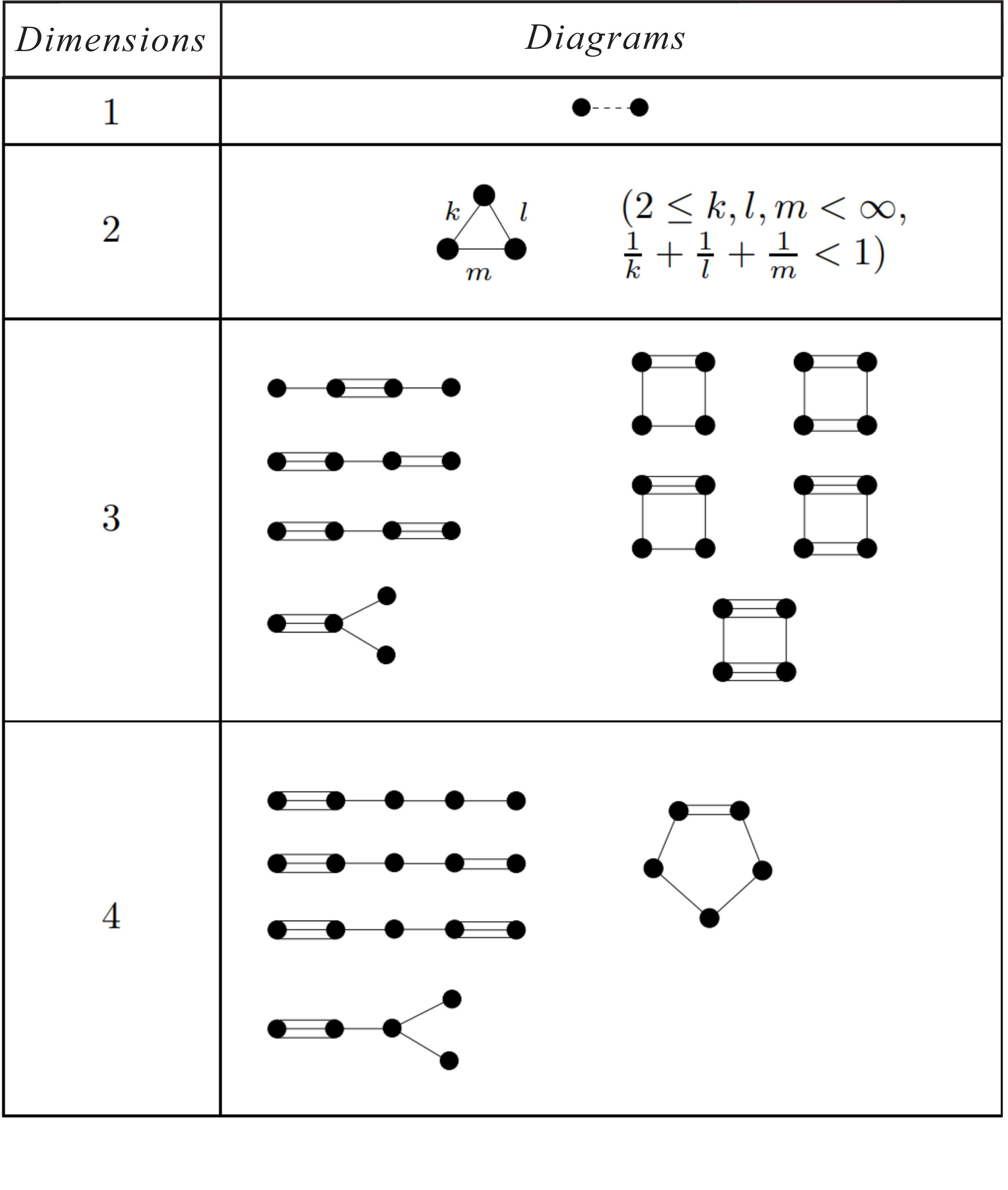}}
	\caption{The Lann\'{e}r diagrams.}
	\label{figure:lanner}
\end{figure}

A connected diagram $\Gamma$ is a \emph{Lann\'{e}r diagram} if $\Gamma$ is neither elliptic nor parabolic; any proper subdiagram of $\Gamma$ is elliptic. Those diagrams are irreducible Coxeter diagrams of compact hyperbolic Coxeter simplices. All such diagrams, reported by Lann\'{e}r \cite{Lanner: 1950}, are listed in Figure \ref{figure:lanner}.

Although the full list of hyperbolic Coxeter polytopes remains incomplete, some powerful algebraic restrictions are known \cite{Vinberg:1985}:

\begin{theorem} \label{thm:signature}
	(\cite{Vinberg:1985s}, Th. 2.1). Let $G=(g_{ij})$ be an indecomposable symmetric matrix
	of signature $(d,1)$, where $g_{ii}=1$ and $g_{ij}\leq 0$ if $i\ne j$. Then there exists a unique (up to isometry of $\mathbb{H}^d$) convex hyperbolic polytope $P\subset\mathbb{H}^d$, whose Gram matrix coincides with $G$.
\end{theorem}

\begin{theorem} \label{Vinberg:thm3.1} 
	(\cite{Vinberg:1985s}, Th. 3.1, Th. 3.2)
	Let $P=\mathop{\cap}\limits_{i\in I} H_i^- \in \mathbb{H}^d$ be a compact acute-angled polytope and $G=G(P)$ be the Gram matrix. Denote $G_J$ the principal submatrix of G formed from the rows and columns whose indices belong to $J\subset I$. Then, 
	\begin{enumerate}
		\item The intersection $\mathop{\cap}\limits_{j\in J}H_j^-,J\subset I$, is a face $F$ of $P$ if and only if the matrix $G_J$ is positive definite
		\item   For any $J\subset I$ the matrix,  $G_J$ is not parabolic.
	\end{enumerate}
\end{theorem}

A convex polytope is said to be \emph{simple} if each of its faces of codimension $k$ is contained in exactly $k$ facets. 
By Theorem \ref{Vinberg:thm3.1}, we have the following corollary:

\begin{corollary}
	Every compact acute-angled polytope is simple.
\end{corollary}

\section{Potential hyperbolic Coxeter matrices}\label{section:potential}

In order to classify all of the compact hyperbolic Coxeter $4$-polytopes with $8$ facets, we firstly enumerate all Coxeter matrices of simple $4$-polytope with $8$ facets that satisfy spherical conditions around all of the vertices. These are named \emph{potential hyperbolic Coxeter matrices} in \cite{JT:2018}. Almost all of the terminology and theorems in this section are proposed by Jacquemet and Tschantz. We recall them here for reference, and readers can refer to \cite{JT:2018} for more details.

\subsection{Coxeter matrices}\label{subsection:coxeterMatrix}

The \emph{Coxeter matrix} of a hyperbolic Coxeter polytope $P$ is a symmetric matrix $M=(m_{ij})_{1\leq i,j\leq N}$ with entries in $\mathbb{N}\cup\{\infty\}$ such that
\[m_{ij}=\left\{\begin{array}{cl} 
1,&\text{if } j=i,\\ 
k_{ij}, &\text{if }H_i\text{ and }H_j\text{ intersect in }\HH^n\text{ with angle }\frac{\pi}{k_{ij}},\\
\infty, & \text{otherwise}.
\end{array}\right.\]
Note that, compared with Gram matrix, the Coxeter matrix does not involve the specific information of the distances of the disjoint pairs.

\begin{remark}
	In the subsequent discussions, we refer to \textit{the Coxeter matrix $M$ of a graph $\Gamma$} as the Coxeter matrix $M$ of the Coxeter polyhedron $P$ such that $\Gamma=\Gamma(P)$.\\
\end{remark}

\subsection{Partial matrices}\label{subsection:partial}

\begin{definition}
	Let $\Omega=\{n\in \mathbb{Z}\,|\,n\geq 2\}\cup\{\infty\}$ and let $\bigstar$ be a symbol representing an undetermined real value. A \textit{partial matrix of size $m\geq 1$} is a symmetric $m\times m$ matrix $M$ whose diagonal entries are $1$, and whose non-diagonal entries belong to $\Omega\cup\{\bigstar\}$.
\end{definition}

\begin{definition}
	Let $M$ be an arbitrary $m\times m$ matrix, and $s=(s_1,s_2,\cdots,s_k)$, $1\leq s_1<s_2<\cdots<s_k\leq m$. Let $M^{s}$ be the $k\times k$ submatrix of $M$ with $(i,j)$-entry $m_{s_i,s_j}$. 
\end{definition}

\begin{definition}
	We say that a partial matrix $M=(p_{ij})_{1\leq i,j,\leq m}$ is a \emph{potential matrix} for a given polytope $P$ if 
	
	$\bullet$ There are no entries with the value $\bigstar$;
	
	$\bullet$ There are entries $\infty$ in positions of $M$ that correspond to disjoint pair;
	
	$\bullet$ For every sequence $s$ of indices of facets meeting at a vertex $v$ of $P$, the matrix, obtained by replacing value $n$ with $\cos\frac{\pi}{n}$ of submatrix $M^s$, is elliptic.
	
\end{definition}

For brevity, we use a \emph{potential vector} $$C=(p_{12},p_{13},\cdots p_{1m},p_{23},p_{24},\cdots,p_{2m},\cdots p_{ij},\cdots p_{m-1,m}),~ p_{ij}\ne \infty$$ to denote the potential matrix, where $1\leq i<j\leq m$ and non-infinity entries are placed by the subscripts lexicographically. The potential matrix  and potential vector $C$  can be constructed one from each other easily. In general, an arbitrary Coxeter matrix corresponds to a Coxeter vector following the same manner. We mainly use the language of \emph{vectors} to explain the methodology and  report the enumeration results. It is worthy to remark that, for a given Coxeter diagram, the corresponding (potential / Gram) matrix and vector are not unique in the sense that they are determined under a given labeling system of the facets and may vary when the system changed. In Section \ref{chapter:algorithm}, we apply a permutation group to the nodes of the diagram and remove the duplicates to obtain all of the distinct desired vectors

For each rank $r\geq 2$, there are infinitely many finite Coxeter groups, because of the infinite 1-parameter family of all dihedral groups, whose graphs consist of two nodes joined by an edge of weight $k\geq 2$. However, a simple but useful truncation can be utilized:

\begin{proposition}
	There are finitely many finite Coxeter groups of rank $r$ with Coxeter matrix entries at most seven.
\end{proposition}

It thus suffices to enumerate potential matrices with entries at most seven, and the other candidates can be obtained from substituting integers greater than seven with the value seven. In other words, we now have more variables, that are restricted to be integers larger than or equal to seven, besides length unknowns.  In the following, we always use the terms ``Coxeter matrix" or ``potential matrix" to mean the one with integer entries less than or equal to seven unless otherwise mentioned.  

In \cite{JT:2018}, the problem of finding certain hyperbolic Coxeter polytopes is solved in two phases. In the first step, potential matrices for a particular hyperbolic Coxeter polytope are found; the ``Euclidean-square obstruction" is used to reduce the number. Secondly, relevant algebraic conditions are solved for the admissible distances between non-adjacent facets.  
In our setting, additional universal necessary conditions, except for the vertex spherical restriction and Euclidean square obstruction, are adopted and programmed to reduce the number of the potential matrices.

\section{ combinatorial type of simple 4-polytopes with 8 facets}\label{section:4d8f}
In 1909, Br\"{u}ckner reported the enumeration of all the different combinatorial types of simple 4-polytopes with 8 facets. Br\"{u}ckner used the Schlegel diagrams to represent all of the combinatorial types of simple 4-polytope. However, not every Br\"{u}ckner's diagram is a
Schlegel diagram. Gr\"{u}nbaum and Sreedharan then used the  "beyond-beneath" technique, for example see [ \cite{G:1967}, Section 5.2], and the Gale diagram developed by M. A. Perles to enumerate once more and corrected some results of Br\"{u}ckner's. Here is the main theorem:

\begin{theorem}[Br\"{u}ckner,Gr\"{u}nbaum-Sreedharan] There are $37$ different combinatorial types of simple $4$-polytopes with $8$ facets.
\end{theorem}

\begin{figure}[H]
	\scalebox{0.15}[0.15]{\includegraphics {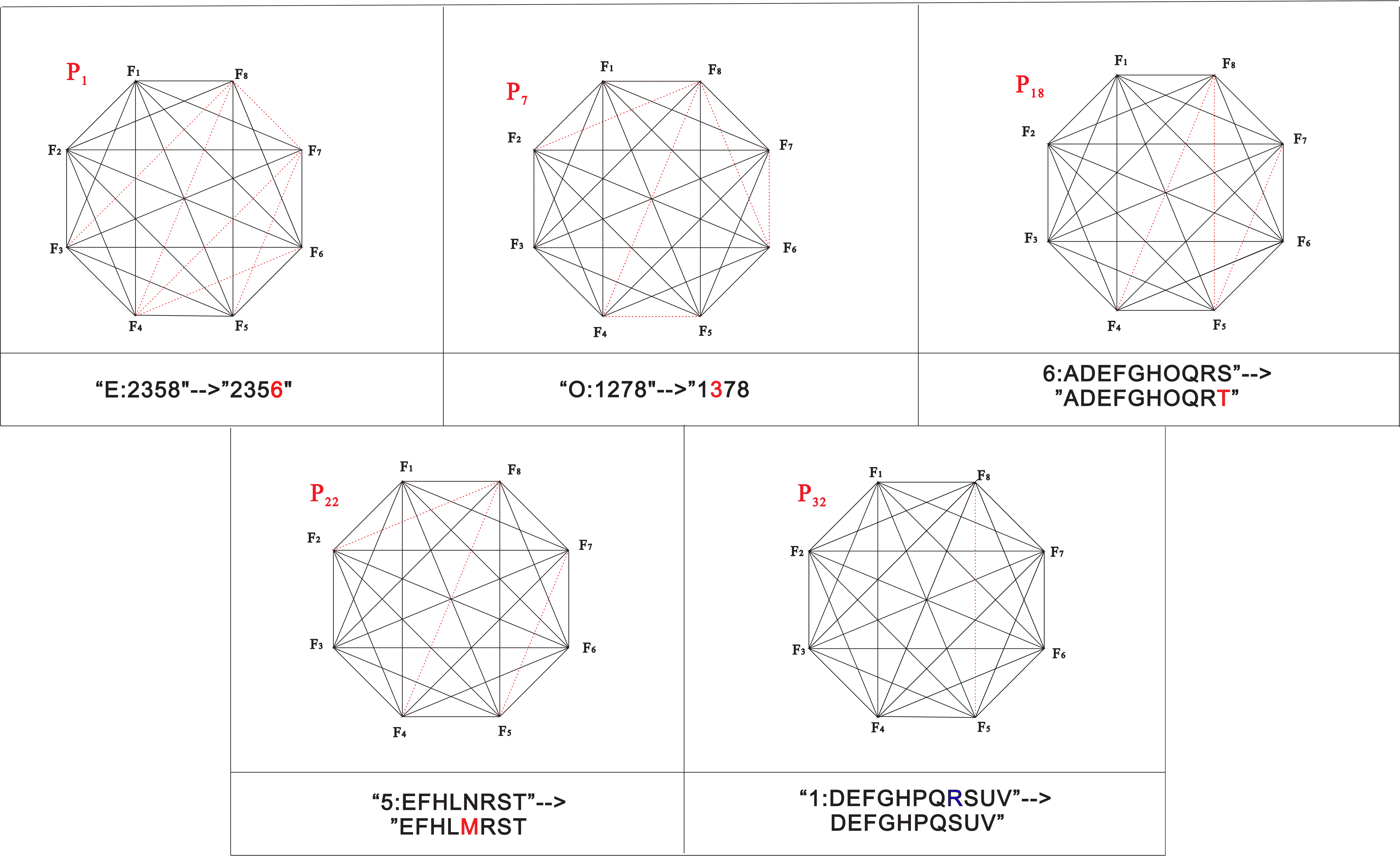}}
	\caption{Corrections to Table 4 in \cite{G:1967}} \label{figure:errs}
\end{figure}

We correct some minor errors in their list as in Figure \ref{figure:errs}, where the polytopes $P^8_i$ in \cite{GS:1967} is now represented by $P_i$ instead.

 Each line in the third column of Table $4$ in \cite{GS:1967} is for one simple polytope with $8$ facets. The data on the first line is as follows:
 
 {\color{blue}
 	\noindent [1,2,4,5] [1,2,3,4] [1,3,4,5] [1,3,5,6] [2,3,5,6] [1,2,3,7] [1,2,6,7] [1,2,5,8] [1,2,6,8] [2,3,4,5] [1,3,6,7] [2,3,6,7] [1,5,6,8] [2,5,6,8]
 	
 }
 
 \noindent where the number $1$, $2$, $\cdots$, $8$ denote the eight facets and each square bracket corresponds to one vertex that is incident to the enclosed four facets. For example, there are $14$ vertices of the above polytope $P_{1}$.
 
 From the original information, we can search out the following \textbf{\emph{data}}  for each polytope:
 
 \begin{enumerate}
 	\item The permutation subgroup $g_i$  of $S_8$ that is isomorphic to the symmetry group of $P_k$;
 	\item The set $d_k$ of pairs of the disjoint facets;
 	\item The set $l_4$ of sets of four facets of which the intersection is of the combinatorial type of a tetrahedron. 
 	\item The set $l_4\_basis$  of sets of four facets of which bound a 3-simplex facet. The label of the bounded 3-simplex facet is recorded as well. Note that $l_4\_basis$  is a subset of $l_4$ and its can be non-empty only for a $4$-dimensional polytope.
 	\item The set $i_2$ of sets of facets, where the intersection is of the combinatorial type of a $2$-cube.
 	\item The set $s_3~ /~ s_4$ of sets of three / four facets of which the intersection is not an edge / a vertex of $P_k$, and no disjoint pairs are included.
 	\item The set $e_3~ /~ e_4$ of sets of three / four facets of which no disjoint pairs are included.
 	\item The set $se_5~ /~ se_6$ of sets of five / six facets of which no disjoint pairs are included.
 	
 \end{enumerate}  
 
 For example, for the polytope $P_{1}$, the above sets are as shown in Table \ref{table:p1data}.

\begin{table}[h]
	{\footnotesize
		\begin{tabular}{|c|c|l|}
			\Xcline{1-3}{1.2pt}
			\multicolumn{3}{|c|}{\textbf{$P_{1}$ }}\\
			\hline
			\multirow{2}{*}{Vert}& \multirow{2}{*}{14}&  $\{\{1, 2, 4, 5\}, \{1, 2, 3, 4\}, \{1, 3, 4, 5\}, \{1, 3, 5, 6\}, \{2, 3, 5, 6\}, \{1, 2, 3, 7\}, \{1, 2, 6, 7\}, \{1, 2, 5, 8\}, $\\
			
			&&$	\{1, 2, 6, 8\}, \{2, 3, 4, 5\}, \{1, 3, 6, 7\}, \{2, 3, 6, 7\}, \{1, 5, 6, 8\}, \{2, 5, 6, 8\}\} $\\
		
			\hline
			$d_{1}$ & 6&$\{ \{3, 8\}, \{4, 6\}, \{4, 7\}, \{4, 8\}, \{5, 7\}, \{7, 8\} \}$ \\
			\hline
			$l_4$ & 3 & $\{\{1, 2, 3, 5\}, \{1, 2, 3, 6\}, \{1, 2, 5, 6\}\}$ \\
			\hline
			$l_4\_basis$& 3 & $\{\{4, \{1, 2, 3, 5\}\}, \{7, \{1, 2, 3, 6\}\}, \{8, \{1, 2, 5, 6\}\}\}$ \\
			\hline
			$s_3$ &0& $\emptyset$ \\
			\hline
			$s_4$& 3 & $\{\{1, 2, 3, 5\}, \{1, 2, 3, 6\}, \{1, 2, 5, 6\}\}$\\		 
			\hline
			\multirow{3}{*}{$e_3$}& \multirow{5}{*}{$28$}&   $\{\{1, 2, 3\}, \{1, 2, 4\}, \{1, 2, 5\}, \{1, 2, 6\}, \{1, 2, 7\}, \{1, 2, 8\}, \{1, 3, 4\}, \{1, 3, 5\}, \{1, 3, 6\}, \{1, 3, 7\},$\\
			&&$\{1, 4, 5\}, \{1, 5, 6\}, \{1, 5, 8\}, \{1, 6, 7\}, \{1, 6, 8\}, \{2, 3, 4\}, \{2, 3, 5\}, \{2, 3, 6\}, \{2, 3, 7\}, \{2, 4, 5\},$\\
			&&$ \{2, 5, 6\}, \{2, 5, 8\}, \{2, 6, 7\}, \{2, 6, 8\}, \{3, 4, 5\}, \{3, 5, 6\}, \{3, 6, 7\}, \{5, 6, 8\},$\\
	
			\hline
			\multirow{2}{*}{$e_4$}& \multirow{2}{*}{$17$}&  
			$\{\{1, 2, 3, 4\}, \{1, 2, 3, 5\}, \{1, 2, 3, 6\}, \{1, 2, 3, 7\}, \{1, 2, 4, 5\}, \{1, 2, 5, 6\}, \{8, 1, 2, 5\}, \{1, 2, 6, 7\}, \{8, 1, 2, 6\}, $\\
			&&$ \{1, 3, 4, 5\}, \{1, 3, 5, 6\}, \{1, 3, 6, 7\}, \{8, 1, 5, 6\}, \{2, 3, 4, 5\}, \{2, 3, 5, 6\}, \{2, 3, 6, 7\}, \{8, 2, 5, 6\}\}  $\\
			
			\hline
			$i_2$ &0& $\emptyset$\\
			\hline
			$se_5$ &$3$& $\{\{1, 2, 3, 4, 5, 6\}, \{2, 3, 4, 5, 6, 7\}, \{3, 4, 5, 6, 7, 8\}, \{4, 5, 6, 7, 8, 9\}\}$\\
			\hline
			$se_6$ &0& $\emptyset$\\
			\hline
			\multirow{12}{*}{$g_{1}$}& \multirow{12}{*}{12}& \multicolumn{1}{c|}{$(1 2 3 4 5 6 7 8)$}\\
			&& \multicolumn{1}{c|}{$(	1 2 3 7 6 5 4 8)$}\\
			&& \multicolumn{1}{c|}{$(1 2 5 4 3 6 8 7)$}\\
			&& \multicolumn{1}{c|}{$(1 2 5 8 6 3 4 7)$}\\
			&& \multicolumn{1}{c|}{$(1 2 6 7 3 5 8 4)$}\\
			&& \multicolumn{1}{c|}{$(1 2 6 8 5 3 7 4)$}\\
			&& \multicolumn{1}{c|}{$(2 1 3 4 5 6 7 8)$}\\
			&& \multicolumn{1}{c|}{$(2 1 3 7 6 5 4 8)$}\\
			&& \multicolumn{1}{c|}{$(2 1 5 4 3 6 8 7)$}\\
			&& \multicolumn{1}{c|}{$(2 1 5 8 6 3 4 7)$}\\
			&& \multicolumn{1}{c|}{$(2 1 6 7 3 5 8 4)$}\\
			&& \multicolumn{1}{c|}{$(2 1 6 8 5 3 7 4)$}\\
			\Xcline{1-3}{1.2pt}
			
		\end{tabular}
	}
	
	\caption{Combinatorics of $P_{1}$.}
	\label{table:p1data}
\end{table}

\begin{figure}[h]
	\scalebox{0.27}[0.25]{\includegraphics {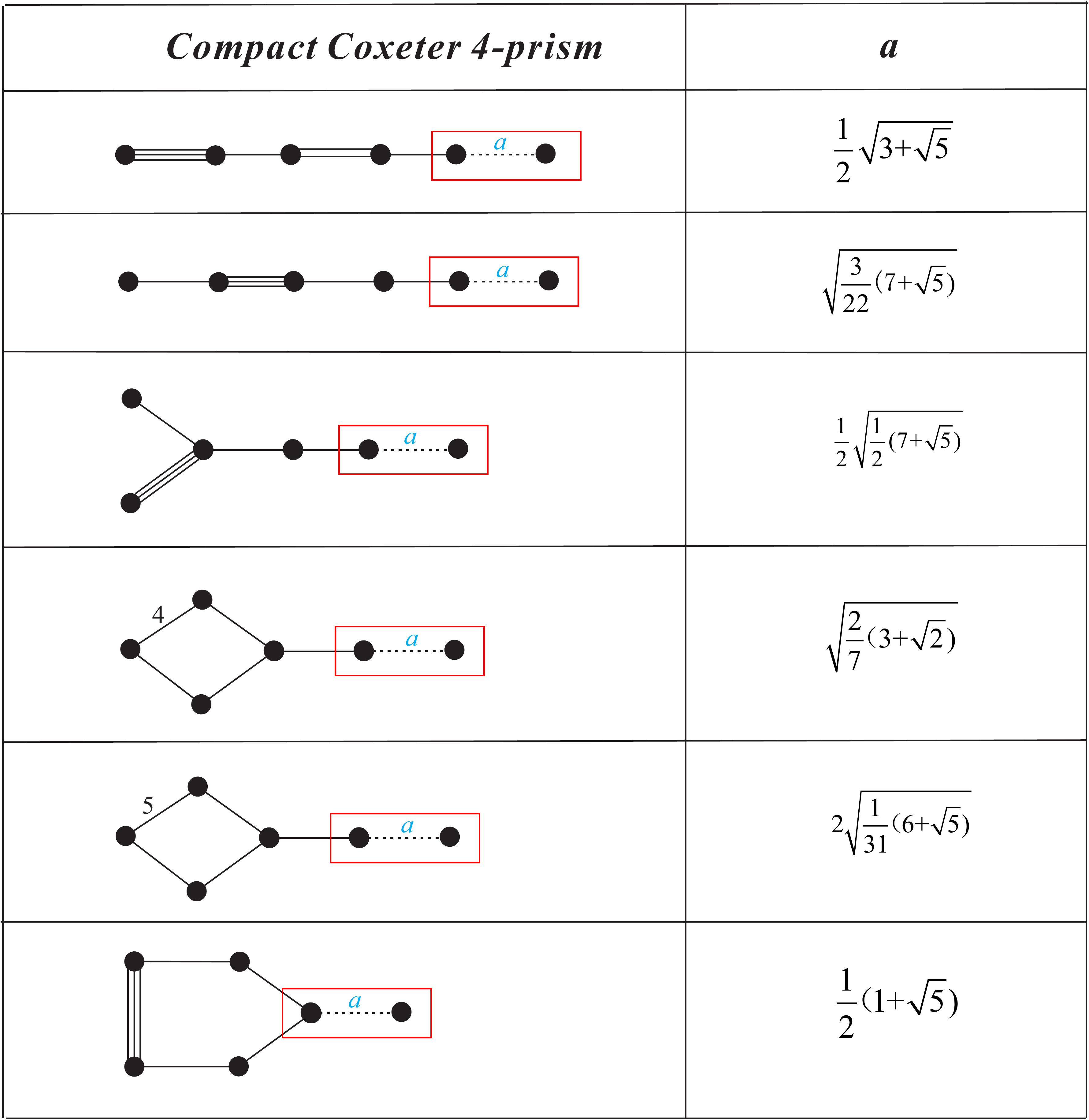}}
	\caption{Compact prisms in $\mathbb{H}^4$.}
	\label{figure:prism4}
\end{figure}

It is worthy to mention that the set $l_4\_basis$, if not empty, can help to reduce the computation since the list of simplicial $4$-prisms is available. For example, suppose $\{1,2,3,5\}$ (referring to facets $F_1,F_2,F_3,F_5$) bound a facet $4$ (means $F_4$) of 3-simplex. Then we can assume that $F_4$ is orthogonal to $F_1$, $F_2$, $F_3$, and $F_5$ (i.e., $m_{14}=m_{24}=m_{34}=m_{54}=2$). The vectors obtained this way can be treated as \emph{bases}, named basis vectors, and all of the other potential vectors that may lead to a Gram vector can be realized by gluing the simplicial $4$-prisms, as shown in Figure \ref{figure:prism4}, at their orthogonal ends.

Moreover, among all of the $37$ polytopes, we only need to study those with number of hyperparall pairs larger than or equal to three due to the following theorems:

\begin{theorem}[\cite{FT:08}, part of Theorem A]
	If $d\leq 4$ and the $d$-polytope $P$ has no pair of disjoint facets then $P$ is either a simplex
	or one of the seven Esselmann polytopes.
\end{theorem}

\begin{theorem}[\cite{FT:09}, Main Theorem A]
	A compact hyperbolic Coxeter $d$-polytope with exactly one pair of non-intersecting facets has at most $d + 3$ facets.
\end{theorem}

\begin{theorem}[\cite{FT:14}, Theorem 7.1]
		Compact hyperbolic Coxeter $4$-polytopes with two pairs of disjoint facets has at most $7$ facets.
\end{theorem}

There are $24$ polytopes with at least three pairs of  hyperparallel facets. We group them by the number $d_k$ of disjoint pairs as illustrated in Table \ref{table: group}.


\begin{table}[h]
	{\footnotesize
		\begin{tabular}{c|c}
			\Xcline{1-2}{1.2pt}
			\textbf{$d_k$} & \textbf{~~ \quad ~~ \quad~~ \quad ~~ \quad ~~ \quad ~~ \quad~~ \quad ~~ \quad labels of polytopes~~ \quad ~~ \quad~~ \quad ~~ \quad ~~ \quad ~~ \quad~~ \quad ~~ \quad}\\
			\hline
			$6$ & 1 2 3 \\
			\hline
			$5$ &4 5 6 7 13 \\
			\hline
			$4$ & 8	9 10 14	15	16	17 34		
			\\
			\hline
			$3$ &11	12	18	19	20	21	22	26
			\\
		
			\Xcline{1-2}{1.2pt}
			
		\end{tabular}
	}
	
	\hspace*{0.5cm}
	\caption{Four groups with respect to different numbers of disjoint pairs.}
	\label{table: group}
\end{table}

\section{Block-pasting algorithms for enumerating all the candidate matrices over certain combinatorial type. } \label{chapter:algorithm}

We now use the \emph{block-pasting algorithm} to determine all of the potential matrices for  the $24$ compact  combinatorial types reported in Section \ref{section:4d8f}. Recall that the entries have only finite options, i.e., $k_{ij}\in\{1,2,3,\cdots, 7\}\cup \{\infty\}$. Compared to the backtracking search algorithm raised in \cite{JT:2018}, ``block-pasting" algorithm is more efficient and universal. Generally speaking, the backtracking search algorithm uses the method of ``a series circuit", where the potential matrices are produced one by one. Whereas, the block-pasting algorithm adopts the idea of ``a parallel circuit", where different parts of a potential matrix are generated simultaneously and then pasted together. 

For each  vertex $v_i$ of a $4$-dimensional hyperbolic Coxeter polytope $P_k$, we define the \emph{i-chunk}, denoted as $k_i$, to be the ordered set of the four facets intersecting at the vertex $v_i$ with increasing subscripts.  For example, for the polytope $P_{1}$ discussed above, there are $14$ chunks as it has $14$ vertices. We may also use $k_i$ to denote the ordered set of subscripts, i.e., $k_i$ is referred to as a set of integers of length four.

Since the compact hyperbolic $4$-dimensional polytopes are simple, each chunk possesses $\tbinom{4}{2}=6$ dihedral angles, namely the angles between every two adjacent facets. For every chunk $k_i$, we define an \emph{i-label} set $e_i$ to be the ordered set $\{10a+b|\{a,b\}\in E_i\}$, where $E_i$, named the \emph{i-index} set, is the ordered set of pairs of facet labels. These are formed by  choosing every two members from the chunk $k_i$ where the labels increase lexicographically. For example, suppose the four facets intersecting at the first vertex are $F_1$, $F_2$, $F_4$, and $F_5$. Then, we have 
$$k_1=\{F_1,F_2,F_4,F_5\} (\text{or} \{1,2,4,5\}),$$ $$E_1=
\{\{1,2\},\{1,4\},\{1,5\},\{2,4\},\{2,5\},\{4,5\}\},$$ $$e_1=
\{12,14,15,24,25,45\}.$$

Next, we list all of the Coxeter vectors of the elliptic Coxeter diagrams of rank $4$. Note that we have made the convention of considering only the diagrams with integer entries less than or equal to seven. The qualified Coxeter diagrams are shown in Figure \ref{figure:elliptic4}:

\begin{figure}[H]
	\scalebox{0.3}[0.3]{\includegraphics {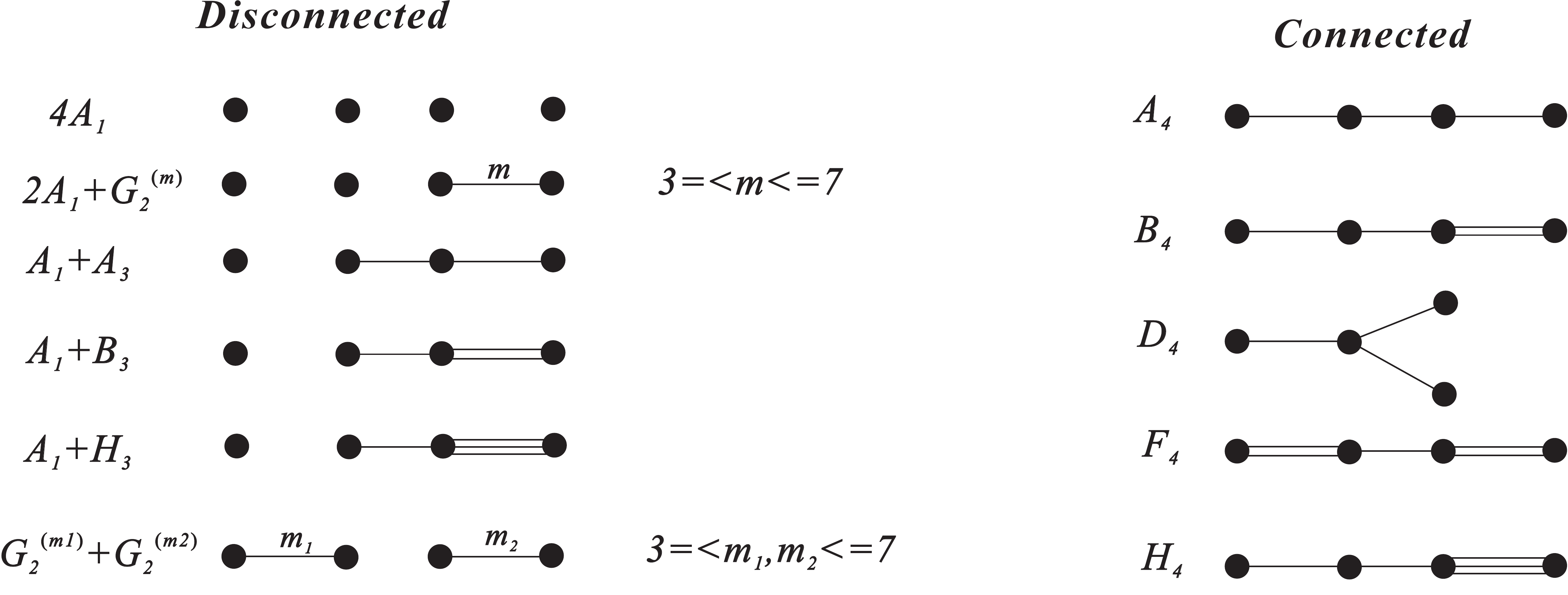}}
	\caption{Spherical Coxeter diagram of rank $4$ with labels less than or equal to seven.}\label{figure:elliptic4}
\end{figure}

\begin{figure}[H]
	\scalebox{0.3}[0.3]{\includegraphics {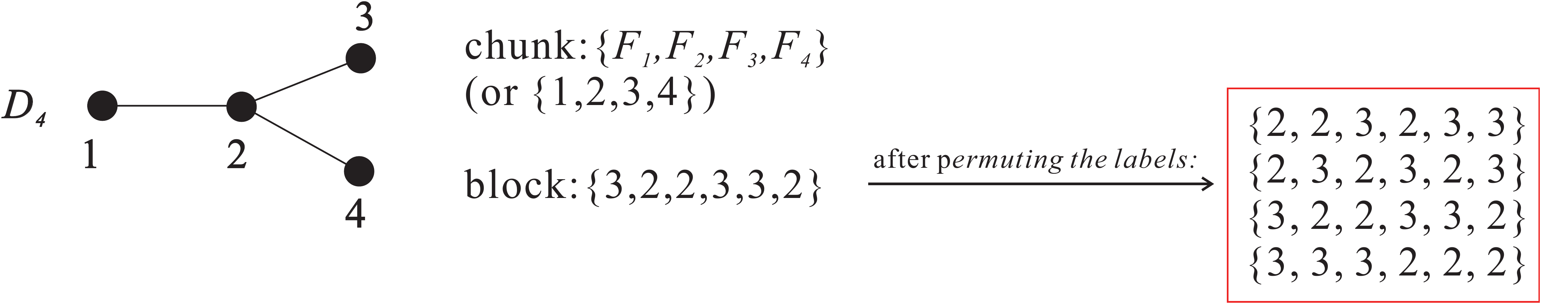}}
	\caption{prepare the pre-block}\label{figure:d4}
\end{figure}

We apply the permutation group on five letters $S_4$ to the labels of the nodes of the Coxeter diagrams in Figure \ref{figure:elliptic4}. This produces all of the possible vectors when varying the order of the four facets. For example, there are $4$ vectors for the single diagram $D_4$ as shown in Figure \ref{figure:d4}.  There are $242$ distinct such vectors of rank $4$ elliptic Coxeter diagrams in total. The set of all of these vectors is called the \emph{pre-block}; it is denoted by $\mathcal{S}$. The set $\mathcal{S}$ can be regarded as a $242\times 6$ matrix in the obvious way. In the following, we do not distinguish these two viewpoints and may refer to $\mathcal{S}$ as either a set or a matrix. 

 We then generate every dataframe $B_i$, named the $i$-block, of size $242\times 6$, corresponding to every chunk $k_i$, for a given polytope $P_k$, where $1\leq i \leq |V_k|$ and $|V_k|$ is the number of vertices of $P_k$. Firstly we evaluate $B_i$ by $\mathcal{S}$ and take the ordered set $e_i$ defined above as the column names of $B_i$. For example, for $e_1=\{12,14,15,24,25,45\}$, the columns of $B_i$ are referred to as $(12)$-, $(14)$-, $(15)$-, $(24)$-, $(25)$-, $(45)$-columns. 
 
 Denote $L$ to be a vector of length $28$ as follows:
 $$L=\{12,13,...,18{\color{red},}~23,24,...,28{\color{red},}~34,35,...,38{\color{red},}~45,46,...,48{\color{red},}~56,57,58{\color{red},}~67,68{\color{red},}~ 78\}.$$
 Then all of the numbers in the label set of $d_k$ (the set of disjoint pair of facets of $P_k$) are excluded from $L$ to obtain a new vector For brevity, the new vector is also denoted by $L$. For example, the numbers excluded for the polytope $P_{1}$ are $38$, $46$, $47$, $48$, $57$, and $78$ as illustrated in Table \ref{table:p1data}. The length of $L$ is denoted by $l$.  
 
 Next, we extend every $242\times 6$ dataframe to a $242\times l$ one, with column names $L$, by simply putting each $(ij)$-column to the position of corresponding labeled column, and filling in the value zero in the other positions. We continue to use the same notation $B_i$ for the extended dataframe. In the rest of the paper, we always mean the extended dataframe when using the notation $B_i$. 

 After preparing all of the the blocks $B_i$ for a given polytope $P_k$, we proceed to paste them up. More precisely, when pasting $B_1$ and $B_2$, a row from $B_1$ is matched up with a row of $B_2$ where every two entries specified by the same index $i$, where $i\in e_1\cap e_2$, have the same values. The index set $e_1\cap e_2$ is called a \emph{linking key} for the pasting. The resulting new row is actually the sum of these two rows in the non-key positions; the values are retained in the key positions. The dataframe of the new data is denoted by $B_1\cup^*B_2$.
 
 We use the following example to explain this process. Suppose
 
 $B_1=\{x_1,x_2\}=\{(1,2,4,4,2,6,0,0,\cdots ,0), (1,2,4,5,2,6,0,0,\cdots,0)\}$,
 
 $B_2=\{y_1,y_2,y_3\}=\{(1,2,4,4,0,0,1,7,0,0,\cdots,0),(1,2,4,4,0,0,6,5,0,0,\cdots,0)$,
 
 \hspace{3.5cm}$(1,2,3,4,0,0,1,7,0,0,\cdots,0)\} .$

 In this example, $x_1$ has the same values with $y_1$ and $y_2$ on the $(12)$-, $(13)$-, $(14)$- and $(15)$- positions. In other words, the linking key here is $\{12,13,14,15\}$. Thus, $y_1$ and $y_2$ can paste to $x_1$, forming the Coxeter vectors $$(1,2,4,4,2,6,1,7,0,...,0)~\text{and}~(1,2,4,4,2,6,6,5,0,...,0), \text{respectively}.$$
 \noindent In contrast, $x_2$ cannot be pasted to any element of $B_2$ as there are no vectors with entry $5$ on the $(15)$-position. Therefore, $$B_1\cup^* B_2=\{(1,2,4,4,2,6,1,7,0,0,\cdots,0),(1,2,4,4,2,6,6,5,0,0,\cdots,0)\}.$$

We then move on to paste the sets $B_1 \cup^*B_2$ and $B_3$. We follow the same procedure with an updated index set. Namely the linking key, is now $e_1\cup e_2\cap e_3$. We conduct this procedure until we finish pasting the final set $B_{\vert V_k\vert}$. The set of linking keys used in this procedure is $$\{e_1\cap e_2, e_1\cup e_2\cap e_3,\cdots,e_1\cup e_2\cup\cdots\cup e_{i-1}\cap e_i,\cdots ,e_1\cup e_2\cup \cdots \cup e_{\vert V\vert -1}\cap e_{\vert V\vert}\}.$$

After pasting the final block $B_{\vert V_k\vert}$, we obtain all of the potential vectors of the given polytope. This approach has been Python-programmed on a PARATERA server cluster.

\begin{figure}[H]
	\scalebox{0.4}[0.4]{\includegraphics {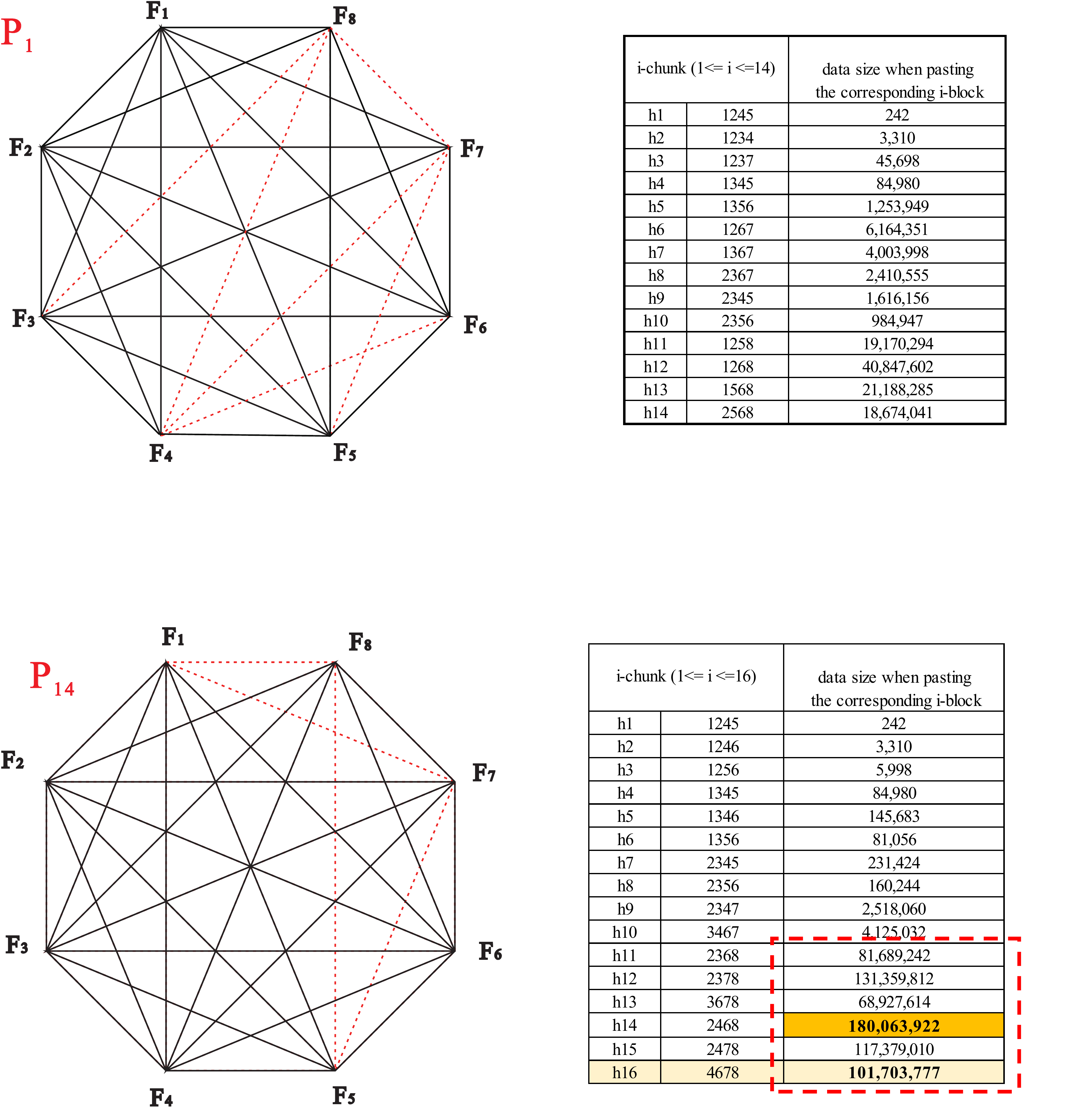}}
	\caption{Use ``block-pasting" algorithm over polytope $P_1$ and $P_{14}$.}\label{figure:p1p14}
\end{figure}

 When apply this approach , it successfully enumerated all truncated candidate for $P_1$. However it usually encounter memory error in solving other case. For example, for polytope $P_{14}$, the computer get stuck at pasting $B_{11}$. We turn to operate it in a serve and it indeed work out finally, see Figure \ref{figure:p1p14} for more details. The peak of the amount of resulting vector has reached $180,063,922$, which far exceed the abilities of storage and computation of an ordinary laptop. Moreover, even using the server, we can not further solve more cases. A refined algorithm is needed to continue this research.

 The philosophy of the refined one is to introduce more necessary conditions other than the vertex spherical restriction, to reduce the amount of vectors in the process of block-pasting. The refined algorithm relies on symmetries of the polytopes and some remarks.

 Firstly, we collect data sets $\mathcal{L}_4$, $\mathcal{L}_4\_basis$, $\mathcal{S}_3$, $\mathcal{S}_4$, $\mathcal{S}_5$, $\mathcal{S}_6$, $\mathcal{E}_3$, $\mathcal{E}_4$,
 $\mathcal{E}_5$, $\mathcal{E}_6$, $\mathcal{I}_2$, as claimed in Table \ref{table:library}, by the following two steps: 
 \begin{enumerate}
 	\item Prepare Coxeter diagrams of rank $r$, as assigned in Table \ref{table:library}, and write down the Coxeter vectors under an arbitrary system of node labeling.
 	\item Apply the permutation group $S_r$ to the labels of the nodes and produce the desired data set consisting of all of the distinct Coxeter vectors under all of the different labelling systems.
 \end{enumerate}

Note that the set $\mathcal{S}_4$ is exactly the pre-block $\mathcal{S}$ we construct before. Readers can refer to the process of producing $\mathcal{S}$ for the details of building these data sets.

\begin{table}[H]
	{\footnotesize
		\begin{tabular}{c|c|c|c}
			\Xcline{1-4}{1.2pt}
			\multirow{2}{*}{\textbf{Types of Coxeter diagrams}}
			& 	\textbf{\# Coxeter}  & \textbf{\# distinct Coxeter Vectors } 
			& \multirow{2}{*}{\textbf{data sets}}\\
			&\textbf{diagrams}& \textbf{after permutation on nodes} & \\   
			\hline	
			Coxeter diagrams of
			compact hyperbolic 3-simplex& 9& 108 &$\mathcal{L}_4$\\
			\hline
			rank $3$ elliptic Coxeter diagrams & 9& 31&$\mathcal{S}_3$ \\
			\hline
			rank $4$ elliptic Coxeter diagrams & 29& 242 & $\mathcal{S}_4$ \\
			\hline
			rank $5$ elliptic Coxeter diagrams &47 & 1946 &$\mathcal{S}_5$\\
			\hline
			rank $6$ elliptic Coxeter diagrams & 117& 20206 &$\mathcal{S}_6$\\
			\hline
			rank $3$ connected parabolic Coxeter diagrams &3 & 10 &$\mathcal{E}_3$ \\
			\hline
			rank $4$ connected parabolic Coxeter diagrams & 3 & 27 &$\mathcal{E}_4$\\	
			\hline
			rank $5$ connected parabolic Coxeter diagrams&5&257&
			$\mathcal{E}_5$\\
			\hline
			rank $6$ connected parabolic Coxeter diagrams&4&870& $\mathcal{E}_6$\\
			\hline
			Coxeter diagrams of Euclidean $2$-cube &4&3&$\mathcal{I}_2$\\
			\Xcline{1-4}{1.2pt}
		\end{tabular}
	}
	
	\caption{Data sets used to reduce the computational load.}
	\label{table:library}
\end{table}


Next, we modify the block-pasting algorithm by using additional metric restrictions. More precisely, remarks \ref{remark:1}--\ref{remark:4}, which are practically reformulated from Theorem \ref{Vinberg:thm3.1}, must be satisfied. 

\begin{remark}(``$l4$-condition") \label{remark:1}
	The Coxeter vector of the six dihedral angles formed by the four facets with the labels indicated by the data in $l_4$ is {\color{red} IN} $\mathcal{L}_4$
\end{remark}

\begin{remark}(``$s3$/$s4$/$s5$/$s6$-condition") \label{remark:2}
	The Coxeter vector of the three/six/ten/fifteen dihedral angles formed by the three/four/five/six facets with the labels indicated by the data in $s_3$/$s_4$/$se_5$/$se_6$ is {\color{red} NOT IN} $\mathcal{S}_3$/$\mathcal{S}_4$/$\mathcal{S}_5$/$\mathcal{S}_6$.
\end{remark}

\begin{remark}(``$e3$/$e4$/$e5$/$e6$-condition") \label{remark:3}
	The Coxeter vector of the three/six/ten/fifteen dihedral angles formed by the three/four/five/six facets with the labels indicated by the data in $e_3$/$e_4$/$se_5$/$se_6$ is {\color{red} NOT IN} $\mathcal{E}_3$/$\mathcal{E}_4$/$\mathcal{E}_5$/$\mathcal{E}_6$.
\end{remark}

\begin{remark}(``$i2$-condition") \label{remark:4}
	The Coxeter vector formed by the four facets with the labels indicated by the data in $i_2$ is {\color{red} NOT IN} $\mathcal{I}_2$.
\end{remark}

The {\color{red}``IN"} and {\color{red}``NOT IN"} tests are called the ``saving" and the ``killing" conditions, respectively. The ``saving" conditions are much more efficient than the ``killing" ones, because the ``what kinds of vectors are qualified" is much more restrictive than the  ``what kinds of vectors are not qualified". Moreover, we remark that the $l3$-condition and the sets $l_3$ and $\mathcal{L}_3$, which can be defined analogously as the $l4$ setting, are not introduced. This is because the effect of using both $s3$- and $e3$- conditions is equivalent to adopting the $l3$-condition.   

We now program these conditions and insert them into appropriate layers during the pasting to reduce the computational load. Here the ``appropriate layer" means the layer where the dihedral angles are non-zero for the first time. For example, for $\{1,2,3\}\in e_3$, we find that after the $j$-th block pasting, the data in columns ($12$-,$13$-,$23$-) of the dataframe $B_1\cup^*B_2\cdots \cup^*B_j$ become non-zero. Therefore, the $e3$-condition for $\{1,2,3\}$ is inserted immediately after the $j$-th block pasting. The symmetries of the polytopes are factored out when the pastes are finished. The matrices (or vectors) after all these conditions (metric restrictions and symmetry equivalence) are called  \emph{``SEILper"-potential matrices (or vectors)} of certain combinatorial types. All of the numbers of the results are reported in Tables \ref{table:gall}. The numbers in red indicate that the corresponding polytopes have a non-empty set $l_4\_basis$; therefore, the results obtained are basis SEILper potential vectors. This calculation is called the \emph{basis approach}. We confirm these cases without using the $l4$-condition, called the \emph{direct approach}, in the validation part presented in Section \ref{section:vadilation}. 

\begin{table}[h]
	{\footnotesize
		\begin{tabular}{c|c|c|c|c|c|c|c|c|c|c}
			\Xcline{1-11}{1.2pt}
			~~ ~~\#~$d_k$~~~~ & ~~~~\textbf{label}~~~~ &~~~~ \# \textbf{SEILper} ~~~~& &~~~~$d_k$~~~~&~~~~\textbf{label}~~~~&~~~~\# \textbf{SEILper} ~~~~&&~~~~$d_k$~~~~&~~~~\textbf{label}~~~~& ~~~~\# \textbf{SEILper} ~~~~\\
			\Xcline{1-3}{0.8pt}\Xcline{5-7}{0.8pt}\Xcline{9-11}{0.8pt}
			6&	{\color{red}1}&	{\color{red} 8}	&&	5&	13&	88,738 && 3&	{\color{red}11}&	{\color{red}0}\\
			\Xcline{1-3}{0.8pt}\Xcline{5-7}{0.8pt}\Xcline{9-11}{0.8pt}
			6&	{\color{red}2}	&{\color{red}12}	&&	4&	{\color{red}9}&	{\color{red}142}	&&	3&	{\color{red}12}&	{\color{red}1,071}\\
			\Xcline{1-3}{0.8pt}\Xcline{5-7}{0.8pt}\Xcline{9-11}{0.8pt}
			6&	{\color{red}3}&	{\color{red}18}	&&	4&	{\color{red}10}&	{\color{red}2}	&&	3&	18&	92,886\\
			\Xcline{1-3}{0.8pt}\Xcline{5-7}{0.8pt}\Xcline{9-11}{0.8pt}
			5&	{\color{red}4}&	{\color{red}231}	&&	4&	14&	0	&&	3&	19&	532\\
		    \Xcline{1-3}{0.8pt}\Xcline{5-7}{0.8pt}\Xcline{9-11}{0.8pt}
			5&	{\color{red}5}&	{\color{red} 398}	&&	4&	15&	4,723&&		3&	20&	138\\
			\Xcline{1-3}{0.8pt}\Xcline{5-7}{0.8pt}\Xcline{9-11}{0.8pt}
			5&	{\color{red}6}&	{\color{red}10}&&		4&	16&	73,006	&&	3&	21&	193,77\\
			\Xcline{1-3}{0.8pt}\Xcline{5-7}{0.8pt}\Xcline{9-11}{0.8pt}
			5&	{\color{red}7}&	{\color{red}4,247}&&		4&	17&	325,957	&&	3&	22&	150,444\\
			\Xcline{1-3}{0.8pt}\Xcline{5-7}{0.8pt}\Xcline{9-11}{0.8pt}
			5&	{\color{red}8}&	{\color{red}2,176}&&		4&	34&	7,608&&		3&	26&	49,599\\
			
			\Xcline{1-11}{1.2pt}
		\end{tabular}
	}
	
	\hspace*{0.5cm}
	\caption{Results of ``SEILper"-potential matrices. Recall that $d_k$ in the table means the number of disjoint pairs as defined before}
	\label{table:gall}
\end{table}

\section{Signature Constraints of hyperbolic Coxeter \texorpdfstring{$n$}-polytopes}\label{section:signature}

After preparing all of the SEILper matrices, we proceed to calculate the signatures of the potential Coxeter  matrices to determine if they lead to the Gram matrix $G$ of an actual hyperbolic Coxeter polytope. 

Firstly, we modify every SEILper matrix $M$ as follows: 

\begin{enumerate}
	\item Replace $\infty$s by length unknowns $x_i$;
	\item Replace $2$, $3$, $4$, $5$, and $6$ by $0$, $-\frac{1}{2}$, $-\frac{l}{2}$, $-\frac{m}{2}$, $-\frac{n}{2}$, where $$l^2-1=2,~l>0,~m^2-m-1=0,~m>0,~n^2-3=0,~n>0;$$
	\item Replace $7$s by angle unknowns of $-\frac{y_i}{2}$.
\end{enumerate} 

By Theorem \ref{thm:signature}, the resulting Gram matrix must have signature $(4, 1)$. This implies that the determinant of every $6\times 6$ minor of each modified $8\times 8$ SEILper matrix is zero. Therefore, we have the following system of $28$ equations and inequality on $x_i$, $l$, $m$, $n$, and $y_i$ to further restrict and lead to the Gram matrices of the desired polytopes:
$$(6.1) ~~~
\begin{cases}
	2\det (M_i)=0,\text{for~ any of the~}\tbinom{8}{2}=28 ~6\times 6~\text{minor}~M_i~\text{of}~ M. \\
	1.8<y_i<2 ~\text{for~ all}~y_i\\
	x_i<-1~\text{for~ all}~x_i\\
	l^2-1=2,~l>0,~m^2-m-1=0,~m>0,~n^2-3=0,~n>0
\end{cases}$$

The above conditions are initially stated by Jacquemet and Tschantz in \cite{JT:2018}. Due to practical constraints in \emph{Mathematica}, we denote $2cos(\frac{\pi}{4}),~2cos(\frac{\pi}{5}),~2cos(\frac{\pi}{6})$ by $\frac{l}{2},~\frac{m}{2},~\frac{n}{2}$, rather than $l,~m,~n$ and set $2\det(M_i)=0$ rather than $\det(M_i)=0$. Delicate reasons for doing so can be found in \cite{JT:2018}. Moreover, we first find the \emph{Gr\"{o}bner bases} of the polynomials involved, i.e., $2\det(M_i),~l^2-1,~m^2-m-1,~ n^2-3$, before solving the system. This might help to quickly pass over the cases that have no solution. However, when dealing with some combinatorial types, like $P_{17},~ P_{22},~P_{13},~P_{16}$, etc., the computation cannot be accomplished in reasonable amount of time. In some cases, a single matrix can require more than two hours to compute, which is costly and impedes the validation process. Hence, we introduce the following $6$-round strategy to make the computation much more feasible and efficient. 

\begin{enumerate}
	\item ``One equation killing"
	
	\noindent Select $2$--$4$ equations, where each of them corresponds to a $6\times 6$ minor and the deleted two rows and columns containing $d_i$  (or $d_i-1$) \footnote{Only if there is no hope to have at least 2 cases where deleted two rows and columns containing all $d_i$ length unknowns, we compromise to find those including $d_i-1$ length unknowns.} length unknowns $x_i$. We use each equation together with the inequalities corresponding to the unknowns in the minor as a condition set and solve them sequentially with time constraint of 1s. 
	\\
	
	\noindent The result consists of ``out set" (the solution is non-empty after the killing), ``left set" (the solution is aborted due to the time constraint), and ``break set" (the solution is empty). We pass the SEILper matrices whose results are either in the ``out set" or the ``left set" to the second round.\\
	\item ``Twenty-eight equations killing"
	
	\noindent We now apply the condition system (6.1) to the SEILper matrices that pass the first round with time constraints 10s, where all of the 28 equations are involved. The result also consists of ``out / left / break" sets. We save the ``out set" to the ``pre-result set" and  pass those in ``left set" to the third round.\\
		
	\item ``Seven equations killing"
		
	\noindent Select $2$--$4$ groups of $7$ equations, where each of the groups corresponds to a $7\times 7$ minor.  There are $\tbinom{7}{6}=7$ $6\times 6$ minors and the deleted one row and column containing as many length unknowns $x_i$ as possible. We use each group of equations together with the inequalities corresponding to the unknowns left in the minor as condition set and solve with time constraint of 1s.\\
	
	\noindent Note that, what is left from the second round are those that can not be solved in 10s when allowing all the $28$ equations in. We therefore reduce the number of constraint equations from $28$ to $7$ to move forward. We pass the SEILper matrices whose results are either in the  ``out set" or the ``left set" to the fourth round.\\
	\item ``Non-Gr\"{o}bner killing"
	
	 \noindent There are not many candidates left up to now. We find in the practice that the function  \emph{GroebnerBasis} in \emph{Mathematica} either works quite quickly or consumes an unaffordable amount of time. We therefore drop this step and proceed to solve the system (6.1) directly with time constraint of 300s.  We save the ``out set" to the ``pre-result set" and pass the ``left set" into next round. \\
	 
	\item ``Range Analysis" 
	
	\noindent What we expect about the angle unknown $y_i$ is not an arbitrary integer in the interval $(1.8,2)$, it should be twice an cosine value of an angle of the form $\displaystyle\frac{\pi}{n},\text{where}~n\geq 7$. Namely, $ \displaystyle \frac{\pi}{\arccos(y_i/2)}\in \mathbb{Z}_{\geq 7}.$ We now choose a small number of equations involving angle and (as less as possible) length unknowns, where the upper bound for $y_i$ is strictly less than $2$, e.g. $1.99$. Then we can run out all the possibilities for the angle unknowns due to the integrality restriction. Alternatively, we might solve for the unknowns explicitly. We use the two methods to analyze the left cases. And for all the polytopes, no new  candidates left from the previous step survived.\\

    \item ``Pre-result checking"
    
    \noindent The last step is to check whether the signature is indeed $(4,1)$ and whether $\displaystyle \frac{\pi}{\arccos~ (y_i/2)}$ is an integer among the pre-result set.

\end{enumerate}

\begin{table}[h]
	{\footnotesize
		\begin{tabular}{c|c|c|c|c|c|c|c}
			\Xcline{1-8}{1.2pt}
			
			\# \textbf{SEILper} & \textbf{R6}&	\#\textbf{Pre-result} &	\textbf{R1} &\textbf{R2} & \textbf{R3} & \textbf{R4} & \textbf{R5}\\
			&{\color{red}\#\textbf{Result}}&&$\#$~left+out &$\#$out / $\#$left &$\#$~left+out &$\#$out / $\#$left&\\
			&&&(excluded rows) &&	(excluded rows)&&\\
			&&&[\#$x_i$ / \#$x_i$ excluded]&&&&\\
			
		    \Xcline{1-8}{0.8pt}
		    325,957	&{\color{red}8}	&{\color{red}47+1=48}& 40934~(7,8)~ [4/3]&	{\color{red}47}/ 89	&15 (8)&
		    {\color{red}1}/ {\color{blue}1}&	0\\
		    &&&
		   22963~(4,6)~[4/3]&&13 (6)& &\\
		   &&&9371~(6,8)~[4/3]&&&&\\
		   &&&8899~(4,8)~[4/3] &&&&\\

			\Xcline{1-8}{1.2pt}
		\end{tabular}
	}
	
	\hspace*{0.5cm}
	\caption{The 6-round procedure about the polytope $P_{17}$.}
	\label{table:p17M}
\end{table}

This approach has been Mathematica-programmed on a PARATERA server cluster. And we illustrate the $6$-round procedure on the polytope $P_{17}$ as an example in Table \ref{table:p17M}. In the fifth round, the twice Gram matrix $2M$ of the only unsolved case (marked in blue in Table \ref{table:p17M}) is as below:

{\tiny
	\begin{center}
		$ 2\begin{pmatrix}
			1&	-(1/2)&	0&	-(1/2)&	-(1/2)&	0&	0&	{\color{red}x_1}\\
			-(1/2)&	1&	-(h_1/2)&	0&	0&	0&	0&	0\\
			0&	-(h_1/2)	&1&	0&	0&	0&	-(1/\sqrt{2})&	0\\
			-(1/2)&	0&	0&	1&	{\color{red}x_2}&	-(1/\sqrt{2})&	-(1/2)&	0\\
			-(1/2)&	0&	0&	{\color{red}x_2}&	1&	-(1/\sqrt{2})&	-(1/2)&	0\\
			0&	0&	0&	-(1/\sqrt{2})&	-(1/\sqrt{2})&	1&	{\color{red}x_3}&	{\color{red}x_4}\\
			0&	0&	-(1/\sqrt{2})&	-(1/2)&	-(1/2)	&{\color{red}x_3}	&1	&0\\
			{\color{red}x_1}&	0&	0&	0&	0&	{\color{red}x_4}&	0&	1\\

		\end{pmatrix}_. $
		
	\end{center}
}

\noindent We use $M_{i,j}$ to denote the minor of the Gram matrix $2M$ after excluding the $i$- and $j$- rows and columns. The minors $M_{7,8}$ and $M_{1,5}$ contain the only angle unknown $h_1$ and one length variable $x_2$. The determinants of them are:

\begin{center}
	$-2 + 5b - 3 x_2^2 -  \sqrt{2} h_1 + \sqrt{2} x_2 h_1 - 2 x_2 h_1^2 + 2 x_2^2 h_1^2 = 0$,\\

	$-4 + 10 x_2 - 6 x_2^2 + h_1^2 - 3 x_2 h_1^2 + 2 x_2^2 h_1^2 = 0$, respectively.\\
\end{center}

\begin{figure}[H]
	\scalebox{0.55}[0.55]{\includegraphics {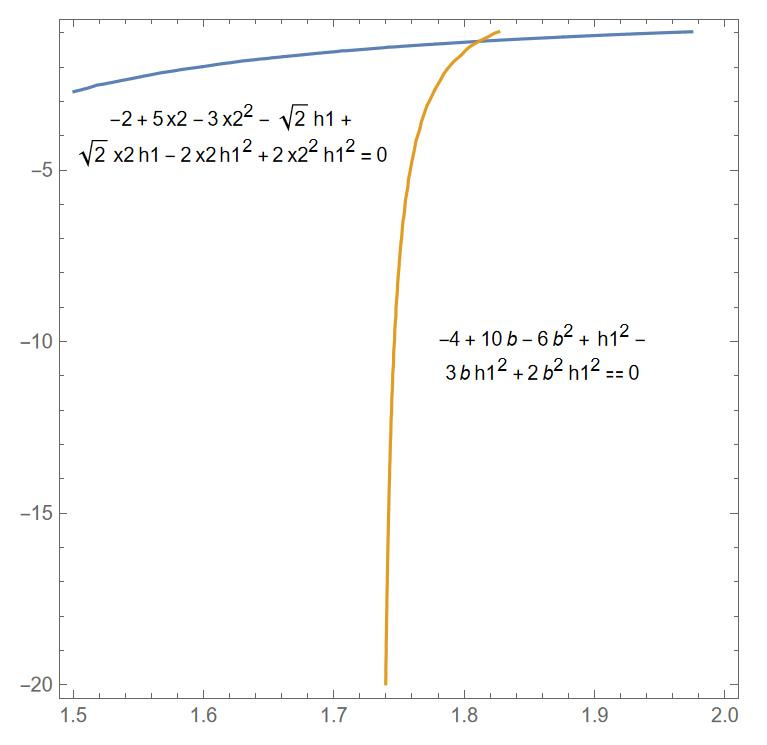}}
	\caption{Tackling with the solved case of $P_{17}$ in the fifth round.} \label{figure:analysis}
\end{figure}

\noindent The graph for these equations are as shown in Figure \ref{figure:analysis}, where  $x_2<-1$ and $1.8<h_1<2$. 

We have two methods to make the claim that no real solutions can be obtained. On one hand, the only solution for the constraints mentioned above is $h_1\approx 1.81129,~ x_2\approx -1.28078$, which means that is no qualified solution for $h_1$, i.e., the solution is not of the form $2\cos \frac{\pi}{n}$. On the other hand, we find out that region bounds for $h_1$ and $x_2$ under the restriction of

$\left\{\begin{array}{l}
	 -2 + 5x_2 - 3 x_2^2 -  \sqrt{2} h_1 + \sqrt{2} x_2 h_1 - 2 x_2 h_1^2 + 2 x_2^2 h_1^2 = 0\\
	 x_2<-1\\
	 1.8<h_1<2\\
\end{array}
\right.$

\noindent are $1.8<h_1< 1.97374$, $-1.3062\leq x_2 <-1$. That means $h_1\in \{2\cos\frac{\pi}{n},n=7,8,9,...,20\}$. We can plug in the value of $h_1$ into some other determinants containing $h_1$ and finally find that the solution set is empty.

After conducting all of these procedures in \emph{Mathematica}, we find that only fourteen of all simple $4$-polytopes with $8$ facets admit compact hyperbolic structure. Besides, we can glue the $3$-prism to seven of them. The results are reported in Table \ref{table:glu}. The polytopes with labels in red are the ``basis" polytopes, and the ones on last line can be obtained by gluing prism ends to those of the second line from the bottom. The Coxeter diagrams and length information are shown in the end of this paper (pages 22--55). 

\begin{table}[h]
	{\footnotesize
			\begin{tabular}{c|c|c|c|c|c|c|c|c|c|c|c|c|c|c}
				\Xcline{1-15}{1.2pt}
				$d_k$ &	\multicolumn{3}{c|}{6} & \multicolumn{5}{c|}{5} & \multicolumn{3}{c|}{4}& \multicolumn{3}{c}{3}\\
				\hline
				polytope&{\color{red}1}\quad& {\color{red}2}\quad	& {\color{red}3}\quad&{\color{red}4}\quad&{\color{red}6}\quad&{\color{red}7}\quad&{\color{red}8}\quad&$13$\quad&$16$\quad&$17$\quad&$34$\quad&$18$\quad&$21$\quad&$26$\quad\\
				\hline
				\# (selected) SEILper potential	&8	&12&	18&	1	&1&	11&4&4&12&48&12&75&1&4\\
				\hline
				\# gram (of basis vectors) &	8&	12&	18&	1	&1&	5&1&3&4&8&12&4&1&2\\
				\hline
				\# gram after suitably gluing & \multirow{2}{*}{130} &\multirow{2}{*}{49}& \multirow{2}{*}{115}&\multirow{2}{*}{2}&\multirow{2}{*}{1} &\multirow{2}{*}{15}&\multirow{2}{*}{2}&\multirow{2}{*}{N}&\multirow{2}{*}{N}&\multirow{2}{*}{N}&\multirow{2}{*}{N}&\multirow{2}{*}{N}&\multirow{2}{*}{N}&\multirow{2}{*}{N}\\
				
				$3$-prisms to the orthoganal ends&&&&&&&&&&&&&&\\	
				
				\Xcline{1-15}{1.2pt}
				
			\end{tabular}
		}
		
		\vspace{0.5cm}
		\caption{Results of the compact hyperbolic Coxeter $4$-polytopes with 8 facets. The value of polytope labeled by $13,16,17,34,18,21,$ or $26$ on the last line is ``N", which means the $l_4\_basis$ set of $P_{13}$/$P_{16}$/$P_{17}$/ $P_{34}$/$P_{18}$/$P_{21}$/$P_{26}$ is empty and we are not be able to glue them with prism ends.}
		\label{table:glu}
	\end{table}

\section{Validation and Results}\label{section:vadilation}

1. ``Basis Approach" vs. ``Direct approach".

We calculate SEILper potential matrices without using  $l_4\_basis$-conditions for those polytopes having $3$-simplex facets. The numbers of Gram matrices corresponding to all of the possible compact hyperbolic polytopes and the results after the \emph{Mathematica} round are reported in Table \ref{table:direct}. They are the same as the previous work accomplished via ``basis approach". 

\begin{table}[H]
	{\footnotesize
			\begin{tabular}{c|ccccccccc}
				\Xcline{1-10}{1.2pt}
				
				grp	&\multicolumn{3}{c|}{6}		&\multicolumn{6}{c}{5}				\\
				\Xcline{1-10}{1.2pt}
				polytopes&	1&	2	&\multicolumn{1}{c|}{3}&	&4&	6&	7&	8 &\\
				\#SEILper&	130&	49&	\multicolumn{1}{c|}{115}&&	571&	26&	8,579	&5,258&\\
				
				\#Gram&	130&	49&	\multicolumn{1}{c|}{115}&&	2&	1&	15&	2&\\
				\Xcline{1-10}{1.2pt}

			\end{tabular}
		}
		
		\hspace*{0.5cm}
		\caption{Results obtained by direct approach. }
		\label{table:direct}
	\end{table} 

\begin{figure}[H]
	\scalebox{0.3}[0.3]{\includegraphics {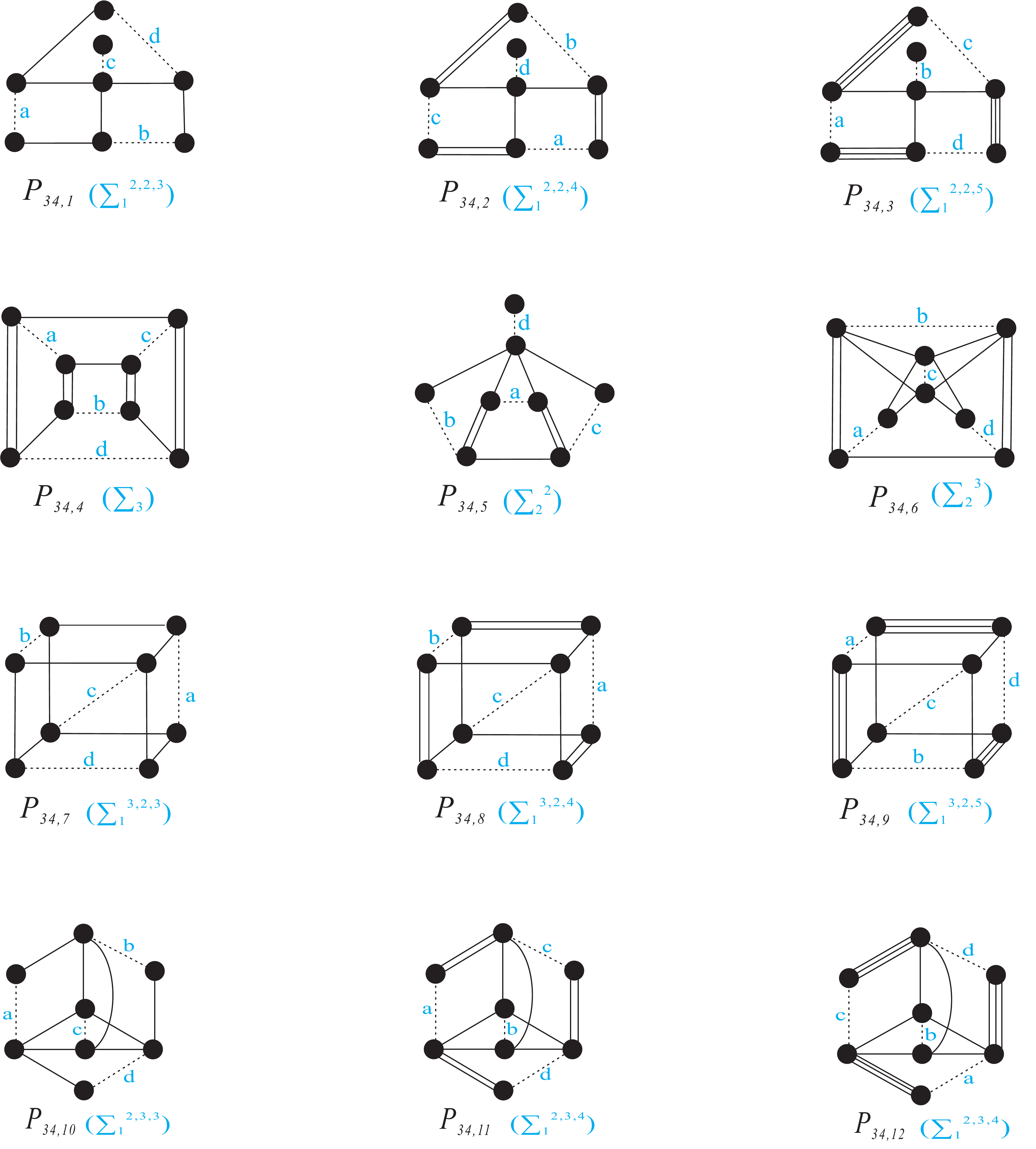}}
	\caption{$P_{34}$} \label{figure:p34}
\end{figure}

2. The compact hyperbolic $4$-cubes are in the family of compact hyperbolic Coxeter $4$-polytope with 8 facets. We obtain the same results of exactly $12$ compact hyperbolic $4$-cubes obtained by Jacquemet and Tschantz in \cite{JT:2018} \footnote{It seems that there is a small typo in Table $5$ of \cite{JT:2018}, where the lengths for $\Sigma_2^2$ and $\Sigma_2^3$ should be swapped. } as shown in Figure \ref{figure:p34}.

\begin{figure}[H]
	\scalebox{0.3}[0.3]{\includegraphics {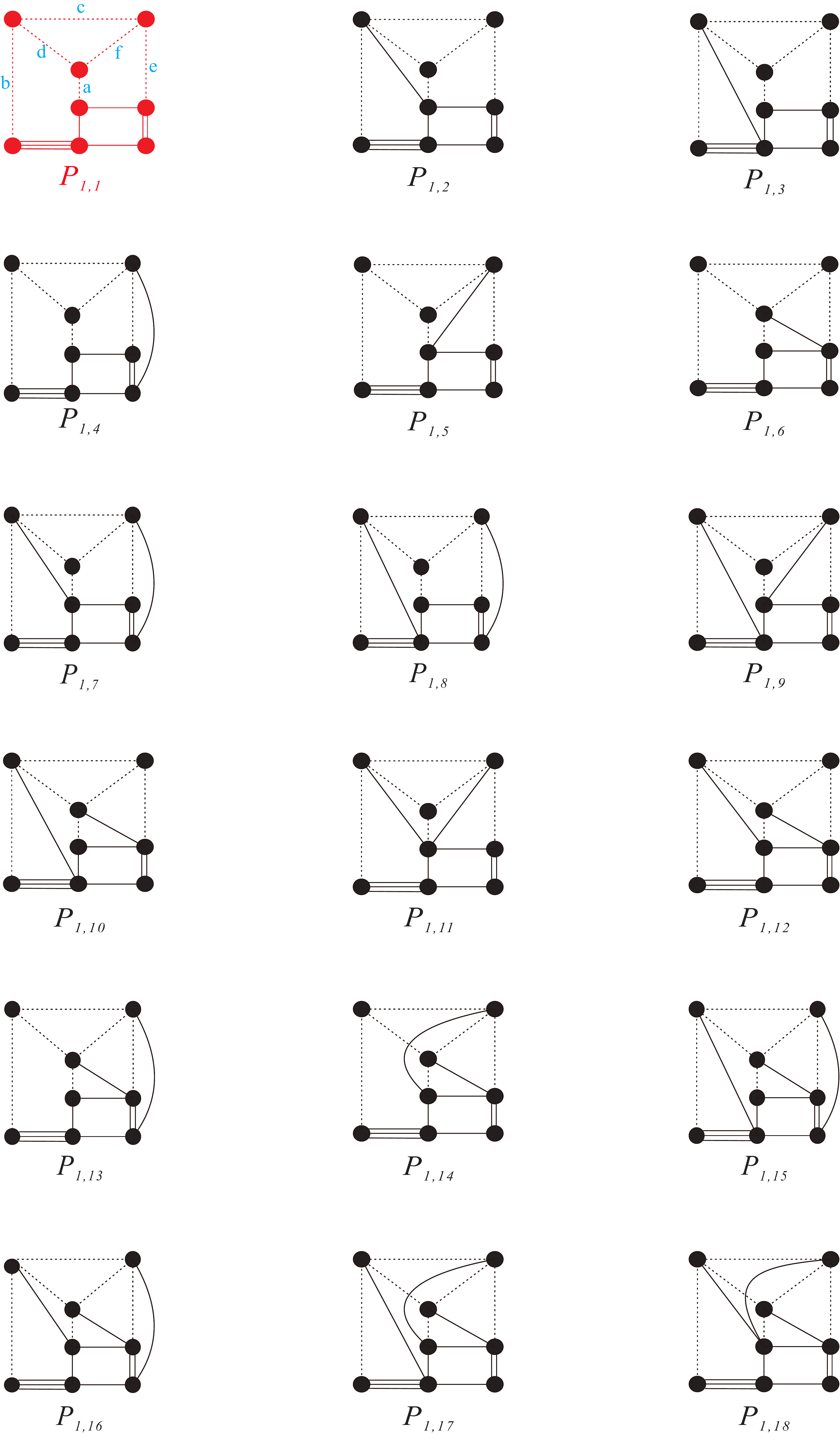}}
	\caption{$P_1$(1/8)} \label{figure:p810}
\end{figure}

\begin{figure}[H]
	\scalebox{0.3}[0.3]{\includegraphics {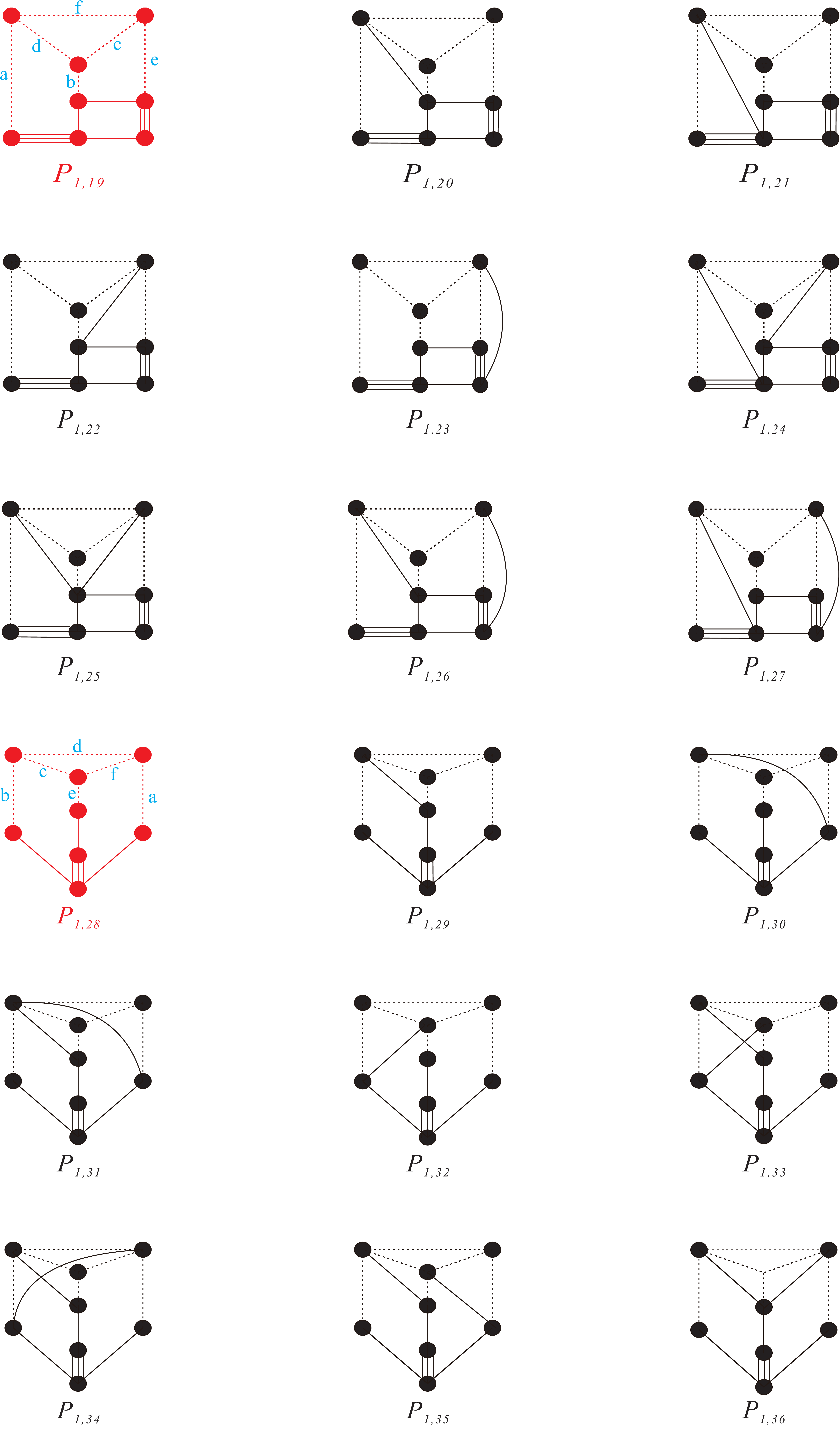}}
	\caption{$P_1$(2/8)} \label{figure:p811}
\end{figure}

\begin{figure}[H]
	\scalebox{0.3}[0.3]{\includegraphics {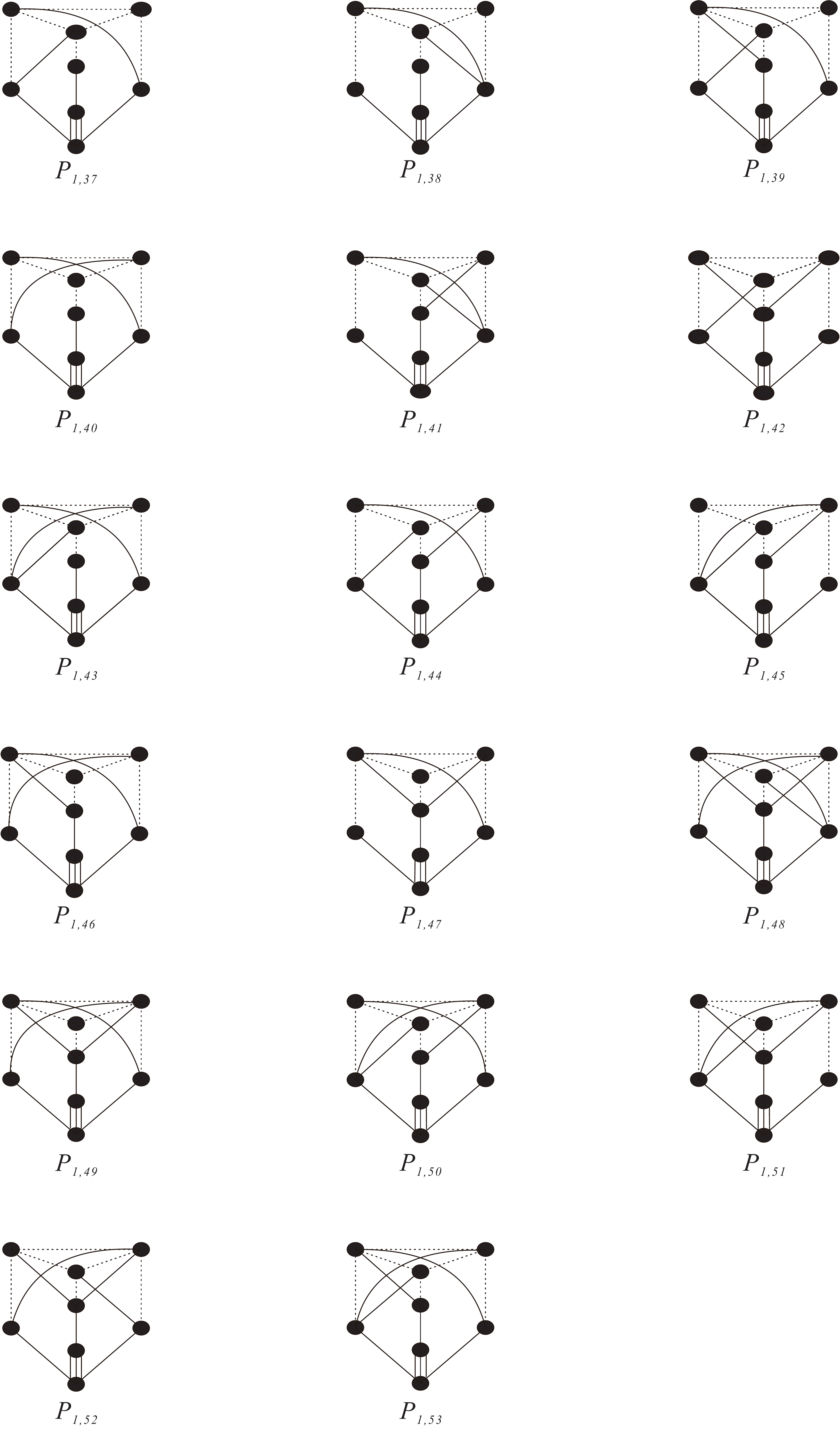}}
	\caption{$P_1$(3/8)} \label{figure:p812}
\end{figure}

\begin{figure}[H]
	\scalebox{0.3}[0.3]{\includegraphics {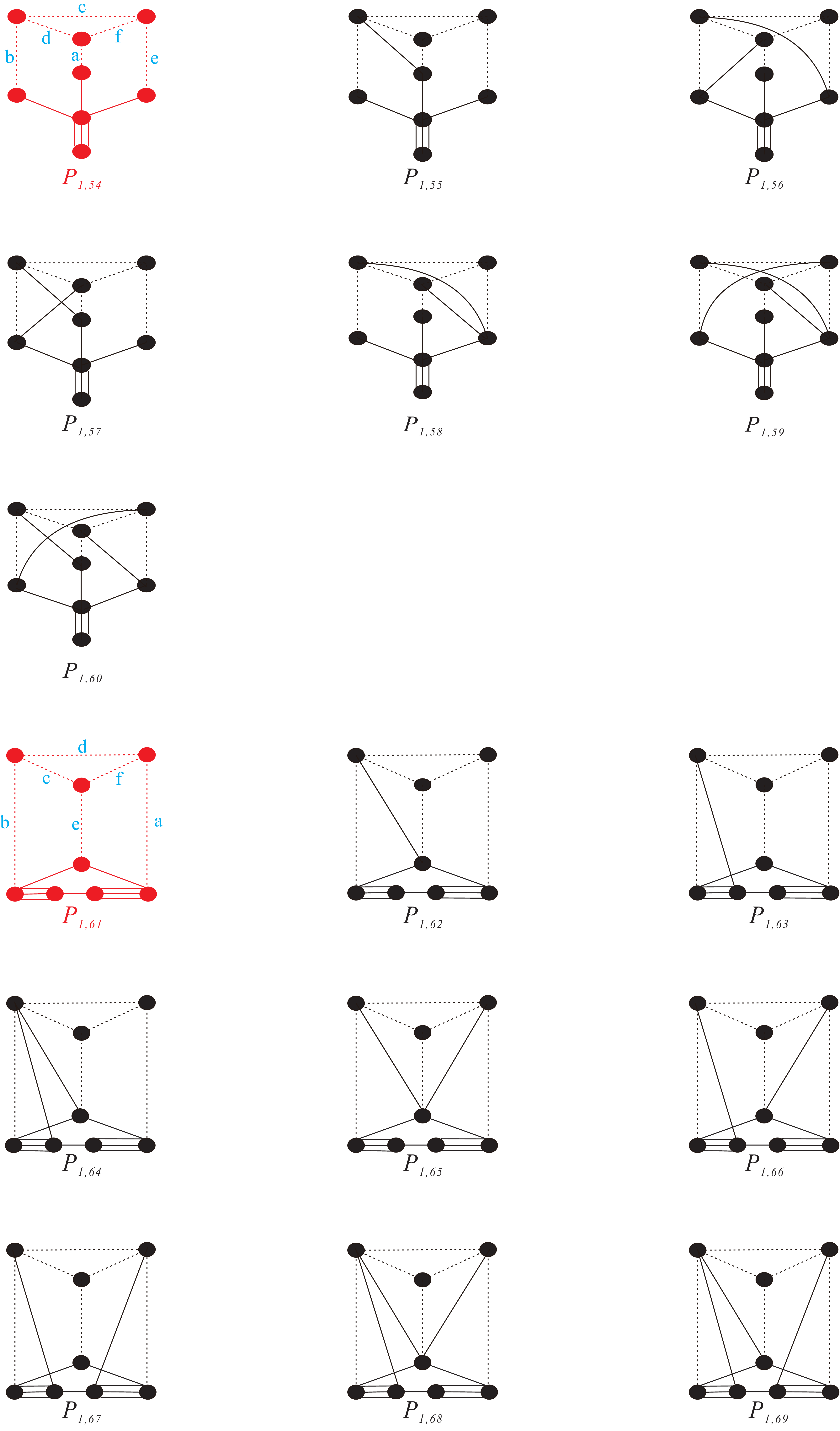}}
	\caption{$P_1$(4/8)} \label{figure:p813}
\end{figure}

\begin{figure}[H]
	\scalebox{0.3}[0.3]{\includegraphics {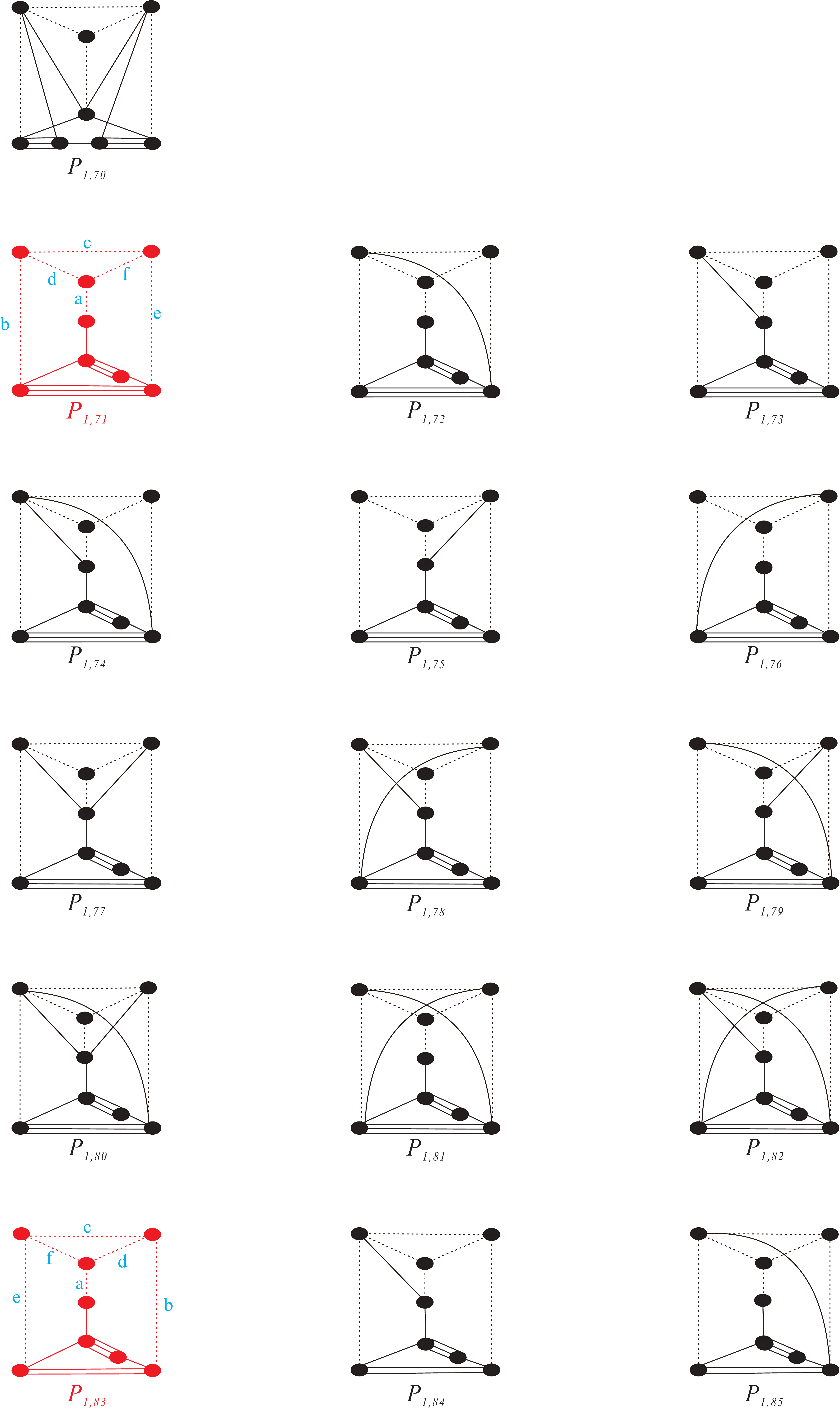}}
	\caption{$P_1$(5/8)} \label{figure:p814}
\end{figure}

\begin{figure}[H]
	\scalebox{0.3}[0.3]{\includegraphics {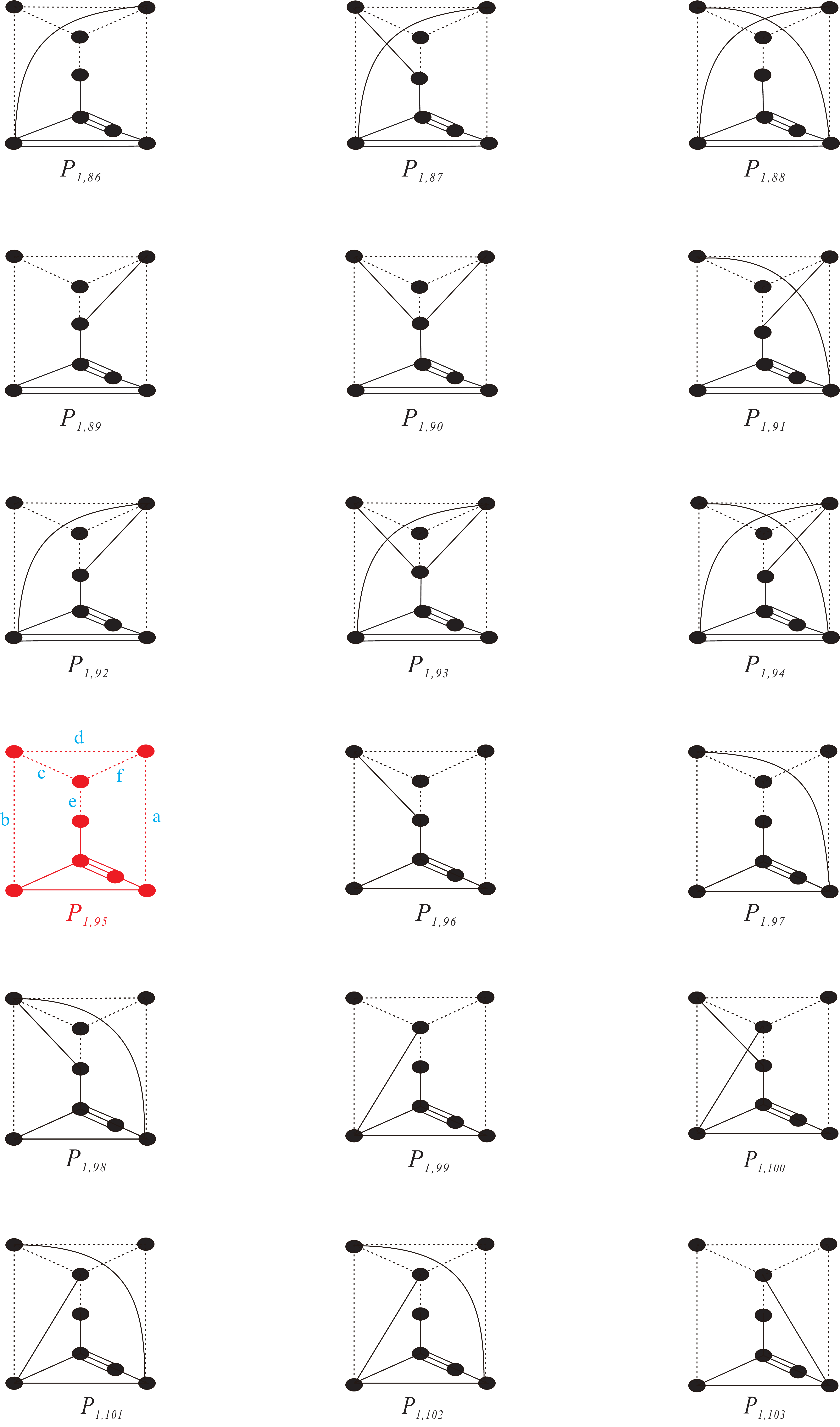}}
	\caption{$P_1$(6/8)} \label{figure:p815}
\end{figure}

\begin{figure}[H]
	\scalebox{0.3}[0.3]{\includegraphics {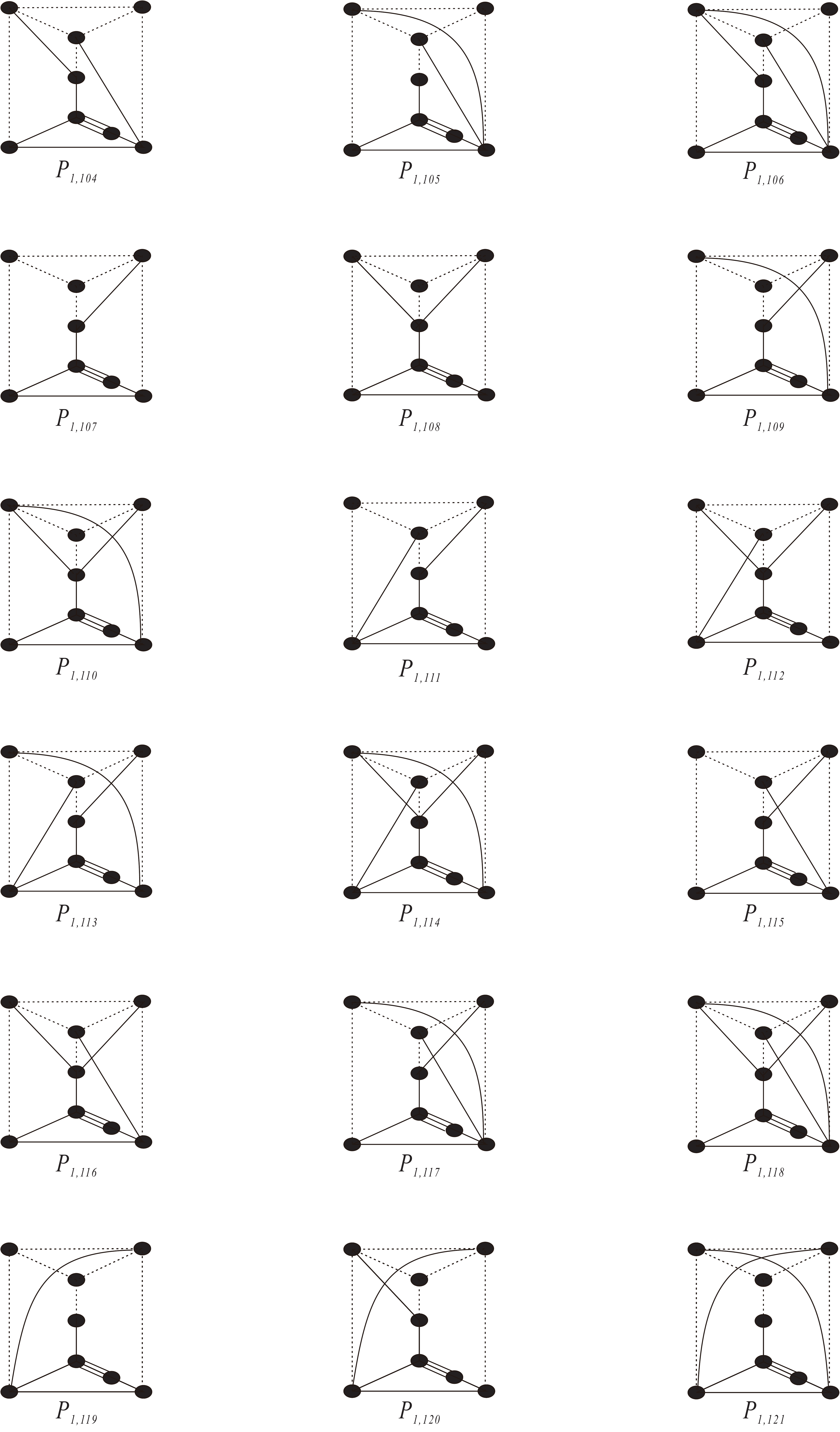}}
	\caption{$P_1$(7/8)} \label{figure:p816}
\end{figure}

\begin{figure}[H]
	\scalebox{0.3}[0.3]{\includegraphics {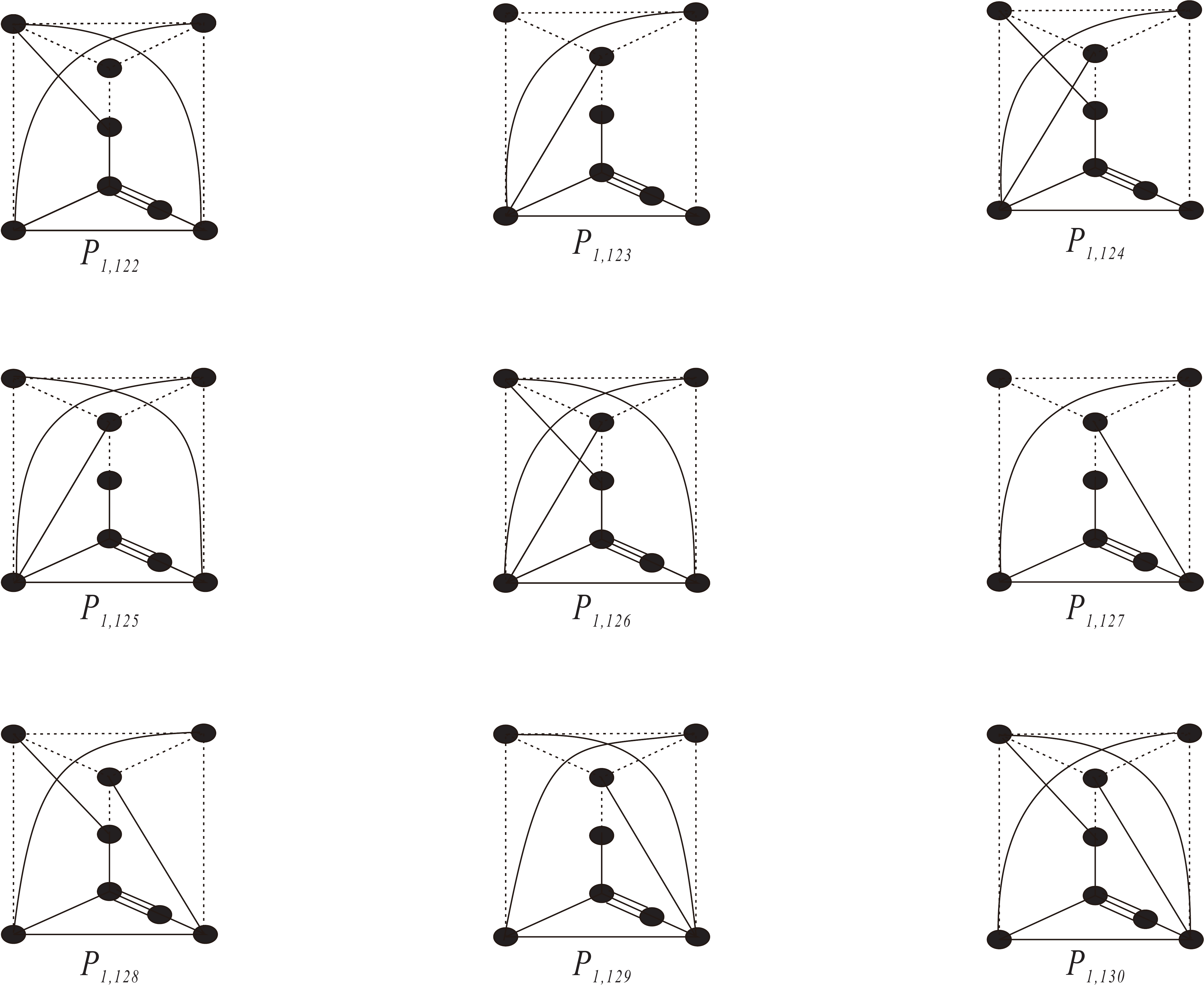}}
	\caption{$P_1$(8/8)} \label{figure:p817}
\end{figure}

3. Flikson and Turmarkin constructed 8 compact hyperbolic Coxeter polytopes with 8 facets in \cite{FT:14} as follows, which are exactly the eight bases of the polytope $P_1$ shown in red in Figure \ref{figure:p810}--\ref{figure:p817}.
\begin{figure}[H]
	\scalebox{0.35}[0.35]{\includegraphics {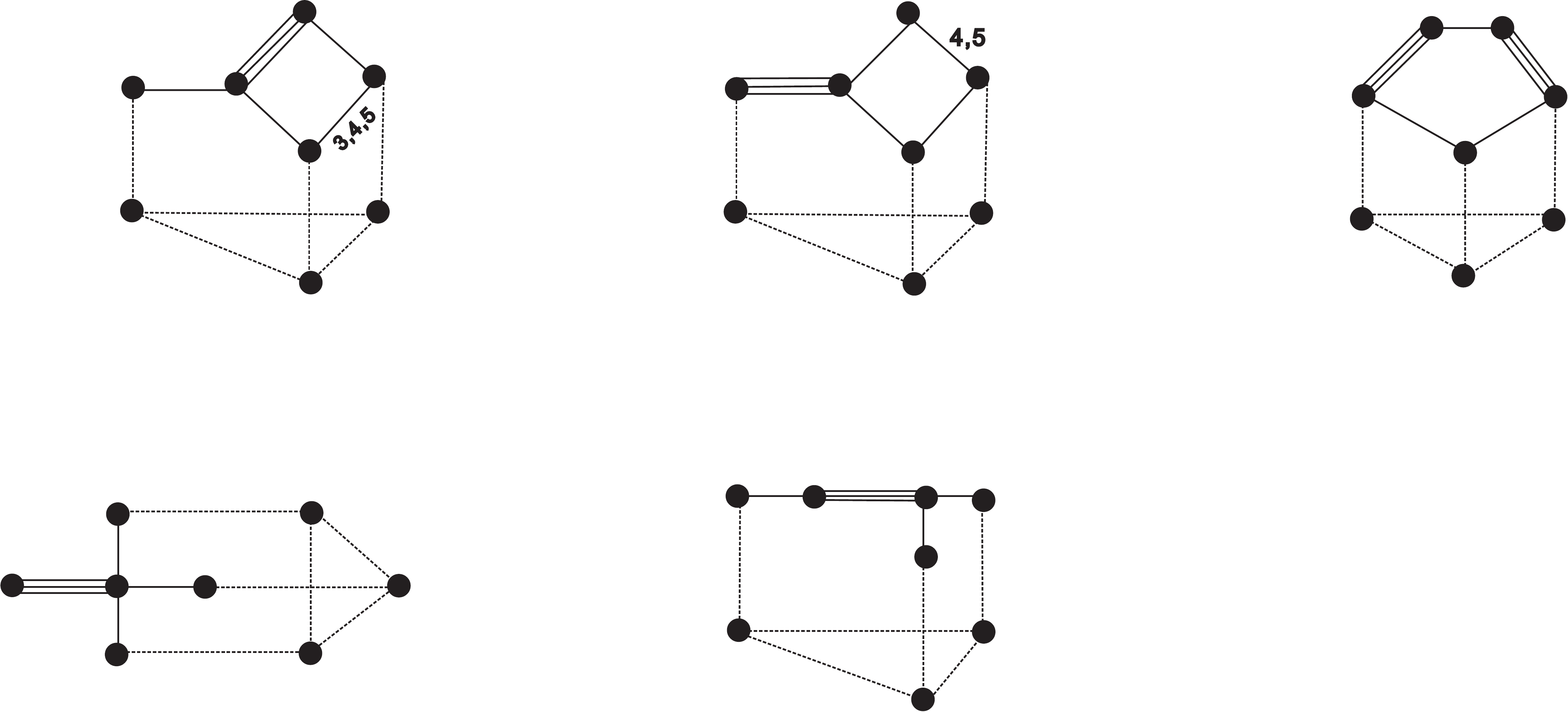}}
	\caption{Known cases from Flikson and Turmarkin} \label{figure:knowcase}
\end{figure}

4. A. Burcroff obtained independently and confirmed the same result in \cite{Amanda:2022} after our mutual check. The correspondence of notions of the polytopes, which admit a hyperbolic structure, between our list and Burcroff's are presented in Table \ref{comparison}.

\begin{table}[h]
	{\footnotesize
			\begin{tabular}{c|cccccccccccccc}
				\Xcline{1-15}{1.2pt}
				
				MZ&	1	&2&	3&	4&	6&	7&8&13&16&17&18&21&26&34\\
				\hline
				A. Burcroff&	$G_1$&	$G_3$	&$G_2$&	$G_7$&	$G_8$&	$G_9$&$G_6$& $G_5$ & $G_{13}$ & $G_{11}$& $G_{12}$& $G_{10}$ & $G_{14}$ &$G_{4}$\\
				\Xcline{1-15}{1.2pt}

\end{tabular}
}

\hspace*{0.5cm}
\caption{Notion correspondence between our list and the list in \cite{Amanda:2022}. }
\label{comparison}
\end{table}

\newpage
\newgeometry{left=0.5cm,right=0.5cm,top=1.5cm,bottom=2cm}

1. Coxeter diagrams for $P_1$

\vspace{1cm}

\restoregeometry

\newpage
\newgeometry{left=0.5cm,right=0.5cm,top=2cm}

\bgroup
\everymath{\displaystyle}
\begin{table}[H]
	\resizebox{17cm}{!}{
		\renewcommand{\arraystretch}{2.5}
		\begin{tabular}{|l|l|l|l|}
			\hline
			\multirow{2}{*}{}                         & a & b & c \\ \cline{2-4} 
			& d & e & f \\ \hline
			\multicolumn{1}{|c|}{\multirow{2}{*}{$P_{1,1}$}} & $  \frac{1}{2}\sqrt{5+3\sqrt{5}+4\sqrt{3+\sqrt{5}}}$ & $   \sqrt{\frac{1}{14}(11+6\sqrt{2}+\sqrt{5(9-4\sqrt{2})})}$  & $    \sqrt{\frac{1}{7}(10+4\sqrt{10}+\sqrt{253+80\sqrt{10}})}$ \\ \cline{2-4} 
			\multicolumn{1}{|c|}{}                    &$    \sqrt{\frac{1}{7}(5+4\sqrt{2})(2+\sqrt{5})}$ & $     \frac{1}{2}\sqrt{\frac{1}{2}(5+3\sqrt{5}+4\sqrt{3+\sqrt{5}})}$ & $     \frac{1}{2}\sqrt{9+4\sqrt{2}+\sqrt{5}}$ \\ \hline
			\multirow{2}{*}{$P_{1,19}$}                       & $  \sqrt{\frac{61}{62}+\frac{5\sqrt{5}}{62}}$ &   $   \sqrt{\frac{73}{38}+\frac{25\sqrt{5}}{38}}$  & $    2\sqrt{\frac{1}{19}(8+3\sqrt{5})}$ \\ \cline{2-4} 
			& $    3\sqrt{\frac{1}{589}(63+26\sqrt{5})}$  & $     \frac{1}{2}\sqrt{7+2\sqrt{5}}$ & $     \sqrt{\frac{74}{31}+\frac{33\sqrt{5}}{31}}$ \\ \hline
			\multirow{2}{*}{$P_{1,28}$}                       & $  \sqrt{\frac{20}{11}+\frac{6\sqrt{5}}{11}}$ & $   \sqrt{\frac{20}{11}+\frac{6\sqrt{5}}{11}}$ & $    \sqrt{\frac{16}{11}+\frac{7\sqrt{5}}{11}}$ \\ \cline{2-4} 
			& $    \frac{3}{11}(3+2\sqrt{5})$ & $  \frac{1}{2}\sqrt{5+\sqrt{5}}$ & $     \sqrt{\frac{16}{11}+ \frac{7\sqrt{5}}{11}}$ \\ \hline
			\multirow{2}{*}{$P_{1,54}$}                       & $  \frac{1}{2}(1+\sqrt{5}) $ &   $    \frac{1}{2}(1+\sqrt{5}) $  & $     \frac{1}{2}(1+\sqrt{5})$  \\ \cline{2-4} 
			& $     \frac{1}{2}(1+\sqrt{5})$  & $      \frac{1}{2}(1+\sqrt{5}) $ & $      \frac{1}{2}(1+\sqrt{5})$  \\ \hline
			\multirow{2}{*}{$P_{1,61}$ }                       & $  \sqrt{\frac{23}{11}+\frac{8\sqrt{5}}{11}}$ &   $   \sqrt{\frac{23}{11}+\frac{8\sqrt{5}}{11}}$  &  $    \sqrt{\frac{23}{11}+\frac{8\sqrt{5}}{11}}$ \\ \cline{2-4} 
			& $    \frac{1}{11}(28+15\sqrt{5})$  & $      \frac{1}{2}(1+\sqrt{5}) $  & $     \sqrt{\frac{23}{11}+\frac{8\sqrt{5}}{11}}$ \\ \hline
			\multirow{2}{*}{$P_{1,71}$}                       & $  \sqrt{\frac{17}{22}+\frac{13}{22\sqrt{5}}}$ &   $   \sqrt{\frac{42}{11}+\frac{17\sqrt{5}}{11}}$  &  $    \sqrt{\frac{65}{22}+\frac{25\sqrt{5}}{22}}$ \\ \cline{2-4} 
			& $    \sqrt{\frac{13}{11}+\frac{2}{\sqrt{5}}}$  &  $     \sqrt{\frac{19}{8}+\frac{7\sqrt{5}}{8}}$ & $     \sqrt{\frac{2}{55}(25+9\sqrt{5})}$ \\ \hline
			\multirow{2}{*}{$P_{1,83}$ }                       & $  \sqrt{\frac{1}{3}(1+\sqrt{10}-\sqrt{7-2\sqrt{10}})}$ &  $   \sqrt{\frac{1}{2}(2+\sqrt{5}+\sqrt{7+3\sqrt{5}})} $ &  $    \sqrt{\frac{1}{11}(24+6\sqrt{2}+5\sqrt{5}+4\sqrt{10})}$ \\ \cline{2-4} 
			& $    (\frac{13}{9}+\frac{4\sqrt{10}}{9})^{1/4}$ & $     \frac{2}{\sqrt{7+6\sqrt{2}-\sqrt{5}-4\sqrt{10}}}$ & $     \sqrt{\frac{1}{33}(1+28\sqrt{2}+4\sqrt{5}(2+\sqrt{2}))}$  \\ \hline
			\multicolumn{1}{|c|}{\multirow{2}{*}{$P_{1,95}$}} & $  \frac{1}{4}(5+\sqrt{5})$ &   $   \sqrt{\frac{65}{22}+\frac{25\sqrt{5}}{22}}$ &  $    \sqrt{\frac{431}{341}+\frac{170\sqrt{5}}{341}}$ \\ \cline{2-4} 
			\multicolumn{1}{|c|}{}                    & $    \sqrt{\frac{26}{11}+\frac{10\sqrt{5}}{11}}$  & $     \sqrt{\frac{5}{31}(6+\sqrt{5})}$ & $     \sqrt{\frac{57}{62}+\frac{25\sqrt{5}}{62}}$  \\ \hline
		\end{tabular}
	}
	
\end{table}
\egroup

\newpage
\newgeometry{left=0.5cm,right=0.5cm,top=1.5cm,bottom=2cm}

1. Coxeter diagrams for $P_2$

\vspace{1cm}
\begin{figure}[H]
	\scalebox{0.3}[0.3]{\includegraphics {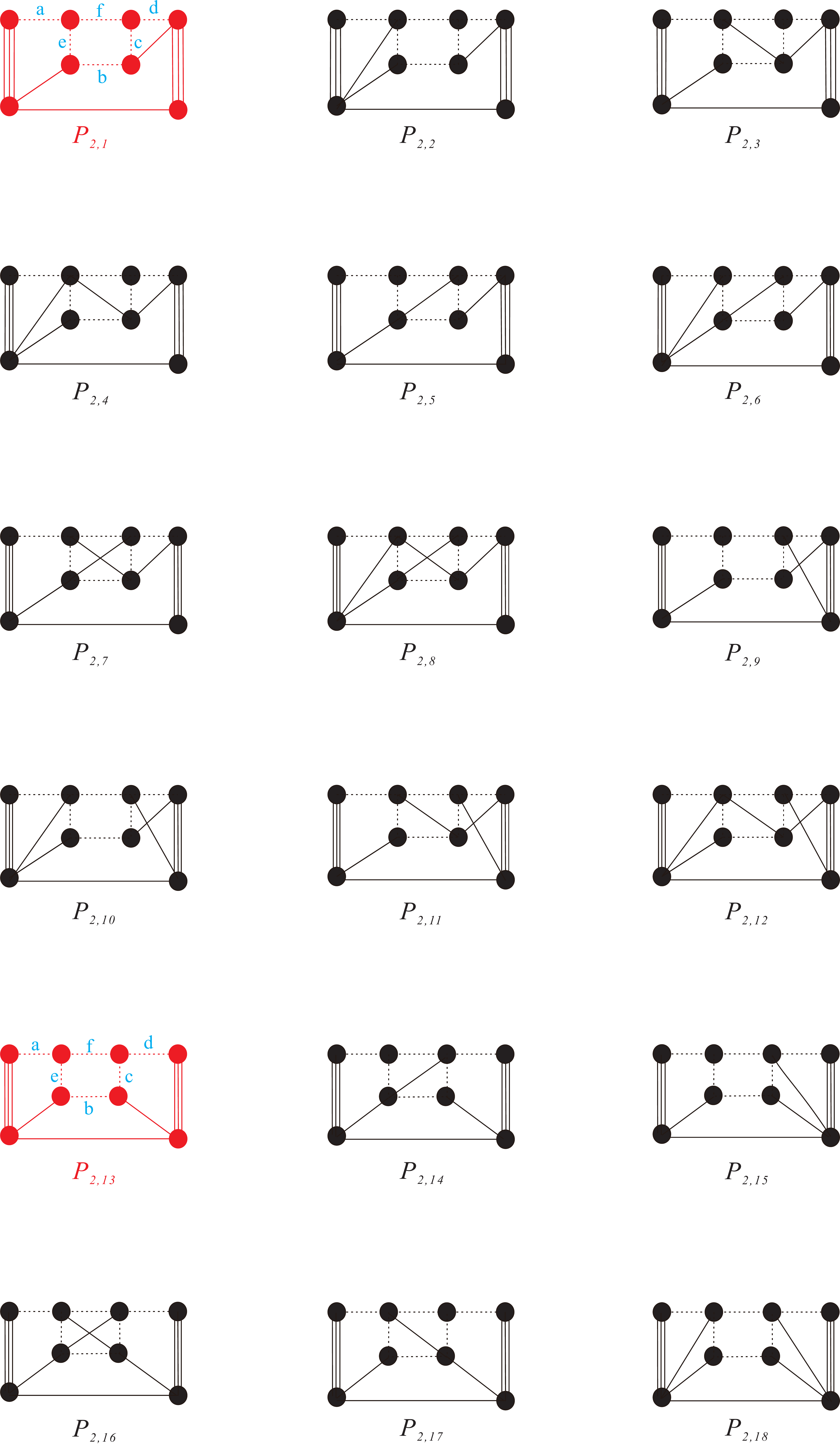}}
	\caption{$P_2$(1/4)} \label{figure:p821}
\end{figure}

\restoregeometry

\begin{figure}[H]
	\scalebox{0.3}[0.3]{\includegraphics {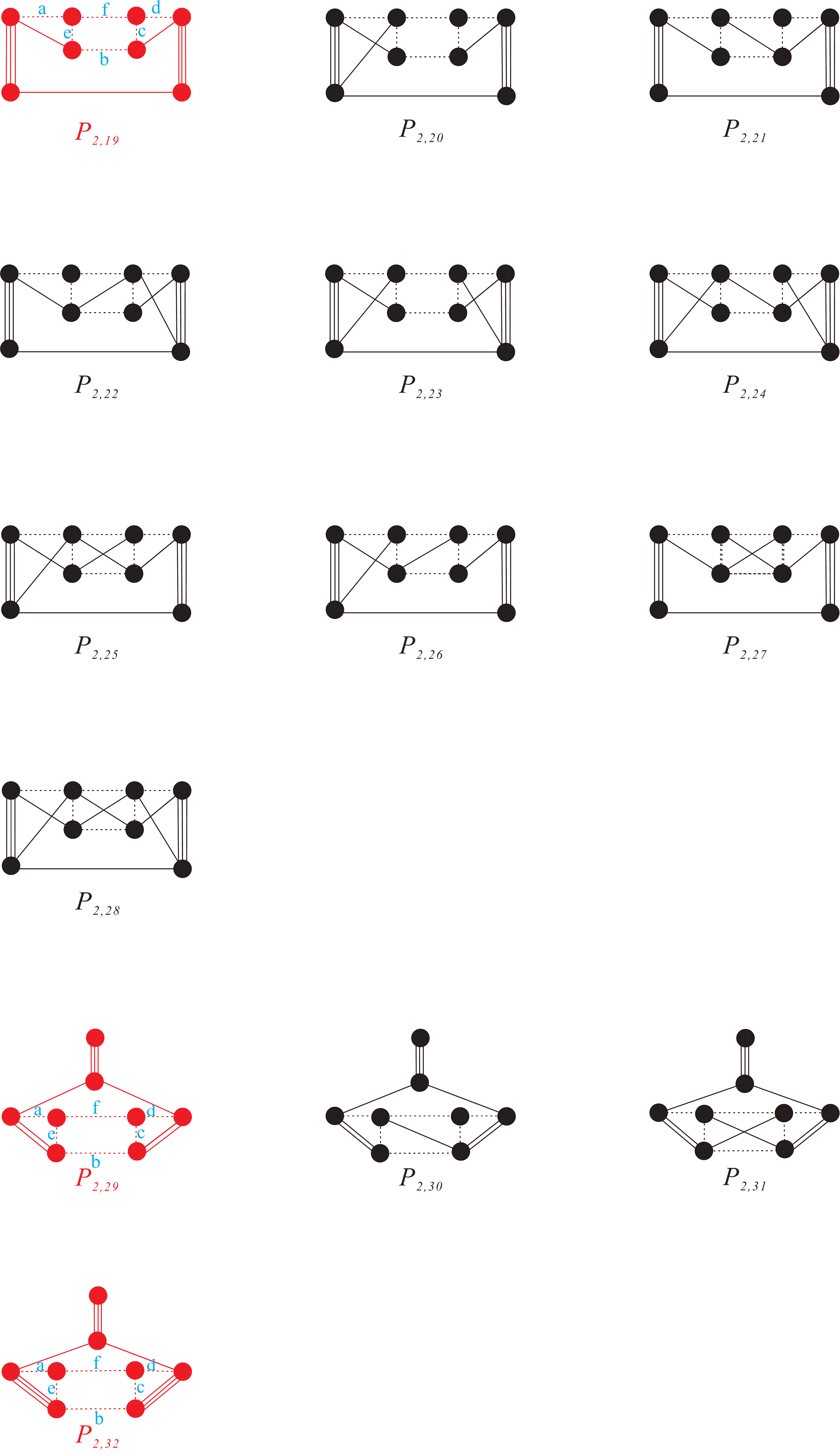}}
	\caption{$P_2$(2/4)} \label{figure:p822}
\end{figure}

\begin{figure}[H]
	\scalebox{0.3}[0.3]{\includegraphics {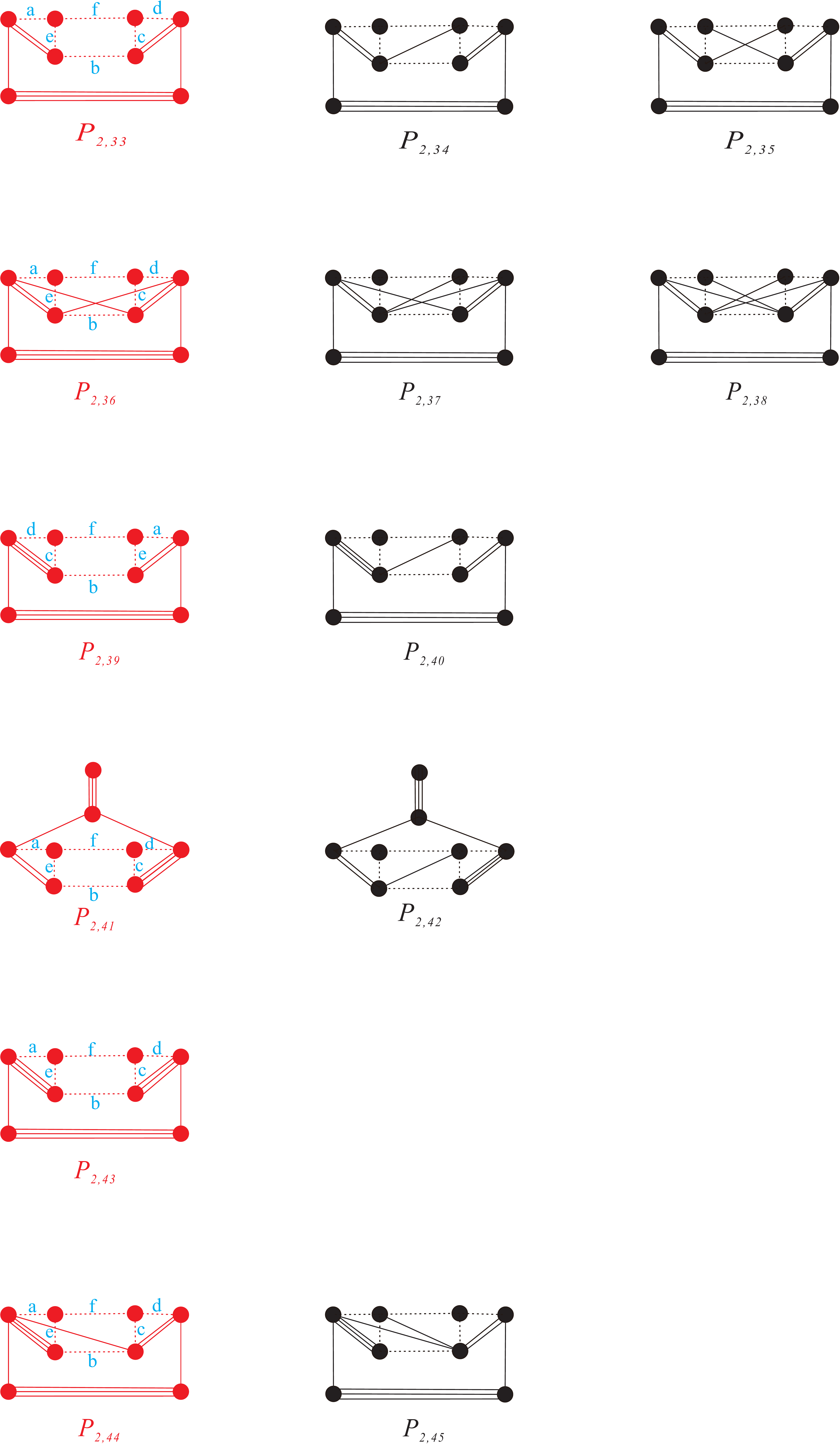}}
	\caption{$P_2$(3/4)} \label{figure:p823}
\end{figure}

\begin{figure}[H]
	\scalebox{0.33}[0.33]{\includegraphics {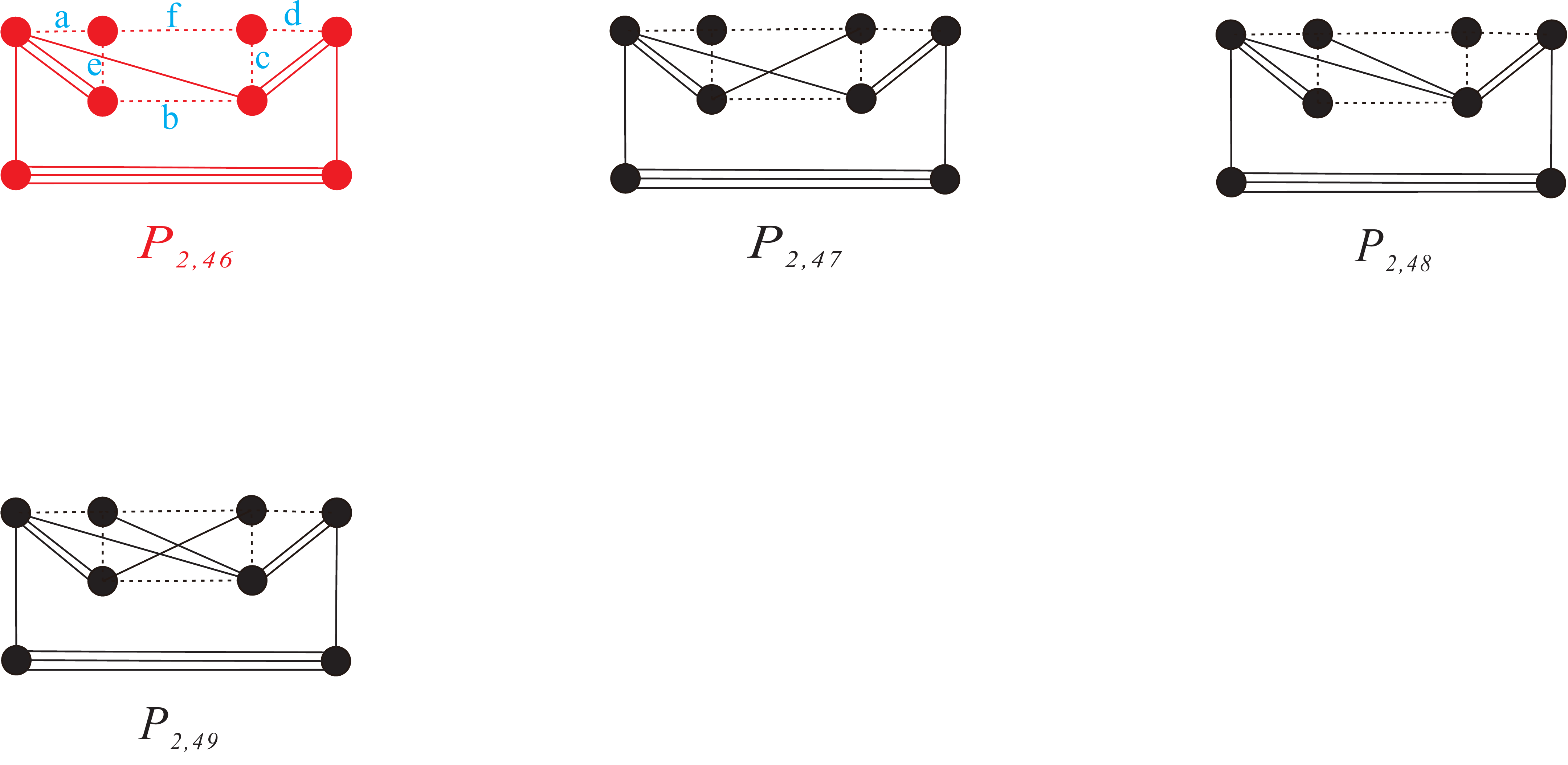}}
	\caption{$P_2$(4/4)} \label{figure:p824}
\end{figure}

\newpage
\newgeometry{left=0.5cm,right=0.5cm,top=0.5cm,bottom=0.5cm}
\bgroup
\everymath{\displaystyle}
\begin{table}[H]
	\resizebox*{19cm}{!}{
		\renewcommand{\arraystretch}{2.7}		
		\begin{tabular}{|c|c|c|c|}
			\hline
			\multirow{2}{*}{} &a  &b  &c  \\ \cline{2-4} 
			&d  &e  &f  \\ \hline
			\multirow{2}{*}{$P_{2,1}$} & $  \sqrt{\frac{2}{11}(7+\sqrt{5})}$ & $   \frac{1}{38}(16+6\sqrt{5}+19\sqrt{\frac{1728}{361}+\frac{496\sqrt{5}}{361}})$ & $    \frac{1}{19}(16+6\sqrt{5}+\sqrt{\frac{1}{2}(169-15\sqrt{5})})$  \\ \cline{2-4} 
			&$    \frac{1}{2}\sqrt{\frac{1}{2}(9+\sqrt{5})}$ & $     \frac{1}{19}\sqrt{\frac{1}{22}(22539+9889\sqrt{5}+4\sqrt{56071390+25075410\sqrt{5}})}$ & $     \frac{1}{19}\sqrt{\frac{1}{22}(76723+34293\sqrt{5}+4\sqrt{526827338+235604346\sqrt{5}})}$  \\ \hline
			\multirow{2}{*}{$P_{2,13}$} & $  \frac{1}{2}\sqrt{\frac{1}{2}(9+\sqrt{5})}$ & $   \frac{1}{2}(3+\sqrt{5})$ & $    \sqrt{\frac{39}{8}+\frac{17\sqrt{5}}{8}}$  \\ \cline{2-4} 
			& $    \frac{1}{2}\sqrt{\frac{1}{2}(9+\sqrt{5})}$ & $     \sqrt{\frac{39}{8}+\frac{17\sqrt{5}}{8}}$ & $     \frac{1}{4}(11+5\sqrt{5})$  \\ \hline
			\multirow{2}{*}{$P_{2,19}$ } & $  \sqrt{\frac{2}{11}(7+\sqrt{5})}$ &   $   \frac{1}{4}(5+\sqrt{5})$  &  $    \sqrt{\frac{26}{11}+\frac{10\sqrt{5}}{11}}$  \\ \cline{2-4} 
			& $     \sqrt{\frac{2}{11}(7+\sqrt{5})}$  &  $     \sqrt{\frac{26}{11}+\frac{10\sqrt{5}}{11}}$ & $     \frac{1}{11}(35+16\sqrt{5})$  \\ \hline
			\multirow{2}{*}{$P_{2,29}$}  &  $  \frac{1}{2}(1+\sqrt{5})$ &  $   1+\frac{\sqrt
				5}{2} $ &  $    \frac{1}{2}\sqrt{5(3+\sqrt{5})}$  \\ \cline{2-4} 
			& $    \frac{1}{2}(1+\sqrt{5})$ & $     \frac{1}{2}\sqrt{5(3+\sqrt{5})}$ & $     \frac{1}{2}(7+3\sqrt{5})$  \\ \hline
			\multirow{2}{*}{$P_{2,32}$} &  $  \sqrt{\frac{43}{38}+\frac{9\sqrt{5}}{38}}$ &   $   \frac{1}{8}(13+3\sqrt{5})$  &  $    \frac{1}{2}\sqrt{\frac{5}{38}(87+35\sqrt{5})}$ \\ \cline{2-4} 
			& $    \sqrt{\frac{43}{38}+\frac{9\sqrt{5}}{38}}$  & $      \frac{1}{2}\sqrt{\frac{5}{38}(87+35\sqrt{5})} $  & $     \frac{1}{38}(83+43\sqrt{5})$  \\ \hline
			
			\multirow{2}{*}{$P_{2,33}$} & $  \frac{1}{2}\sqrt{5+\sqrt{5}}$ & $   \frac{1}{2}(3+\sqrt{5})$ & $    \sqrt{5+2\sqrt{5}}$   \\ \cline{2-4} 
			&$    \frac{1}{2}\sqrt{5+\sqrt{5}}$ & $     \sqrt{5+2\sqrt{5}}$ & $     \frac{1}{2}(7+3\sqrt{5})$ \\ \hline
			\multirow{2}{*}{$P_{2,36}$} &  $  \sqrt{\frac{1}{2}(2+\sqrt{5}+\sqrt{7+3\sqrt{5}})} $ &   $    \frac{1}{4}(7+3\sqrt{5}+4\sqrt{3+\sqrt{5}}) $  & $ \frac{1}{2}\sqrt{61+43\sqrt{2}+\sqrt{5(1451+1026\sqrt{2})}} $  \\ \cline{2-4} 
			& $  \sqrt{\frac{1}{2}(2+\sqrt{5}+\sqrt{7+3\sqrt{5}})}$  & $  \frac{1}{2}\sqrt{61+43\sqrt{2}+27\sqrt{5}+19\sqrt{10}}$ & $ \sqrt{2}(2+\sqrt{5})+\frac{3}{2}(3+\sqrt{5})$   \\ \hline
			
			\multirow{2}{*}{$P_{2,39}$} & $ \sqrt{\frac{2}{19}(9+\sqrt{5})}$ & $   \frac{1}{22}\sqrt{1063+419\sqrt{5}+8\sqrt{29530+13204\sqrt{5}}}$ &  $ \frac{1}{22}(26+10\sqrt{5}+11\sqrt{\frac{340}{121}+\frac{124\sqrt{5}}{121}})$  \\ \cline{2-4} 
			& $   \frac{1}{2}\sqrt{5+\sqrt{5}}$ & $     \frac{1}{11}\sqrt{\frac{1}{19}(5426+1707\sqrt{5}+8\sqrt{356050+159164\sqrt{5}})}$  &$ \frac{11}{\sqrt{4628-2011\sqrt{5}-4\sqrt{52445-22999\sqrt{5}}}}$  \\ \hline
			\multirow{2}{*}{$P_{2,41}$}& $ \sqrt{\frac{43}{38}+\frac{9\sqrt{5}}{38}}$ &   $  \frac{1}{4}\sqrt{23+9\sqrt{5}+2\sqrt{202+90\sqrt{5}}}$  &  $\frac{1}{4}(1+\sqrt{5}+2\sqrt{8+3\sqrt{5}}) $  \\ \cline{2-4} 
			& $\frac{1}{2}(1+\sqrt{5})$  &  $  \frac{1}{2}\sqrt{\frac{1}{19}(143+37\sqrt{5}+2\sqrt{2442+1082\sqrt{5}})}$ & $   \sqrt{\frac{1}{38}(319+141\sqrt{5}+2\sqrt{43618+19506\sqrt{5}})}$  \\ \hline
			\multirow{2}{*}{$P_{2,43}$} &  $  \sqrt{\frac{2}{19}(9+\sqrt{5})} $ &   $    2+\frac{\sqrt{5}}{2} $  & $    \sqrt{\frac{78}{19}+\frac{34\sqrt{5}}{19}} $ \\ \cline{2-4} 
			& $      \sqrt{\frac{2}{19}(9+\sqrt{5})}$  & $       \sqrt{\frac{78}{19}+\frac{34\sqrt{5}}{19}}$ & $      \frac{1}{19}(47+20\sqrt{5})$   \\ \hline
			
			\multirow{2}{*}{$P_{2,44}$} &  $ 1.9923$ &   $4.04306 $  &  $ 4.7181$  \\ \cline{2-4} 
			& $ 1.08754$  & $      4.60491 $  & $  8.13984$  \\ \hline
			
			\multirow{2}{*}{$P_{2,46}$} &  $ 1.9923$ &  $ 3.71598 $ &  $   5.56213$  \\ \cline{2-4} 
			& $  1.345$ & $   4.25068$ & $  9.72857$  \\ \hline
		\end{tabular}
	}

\end{table}
\egroup

\restoregeometry

\newpage
\newgeometry{left=0.5cm,right=0.5cm,top=2cm}
\begin{remark}.
	\vspace{0.5cm}
	
	For $P_{2,44}$, the acute solutions are:
	
	$f=\displaystyle \frac{1}{11\sqrt{19}}
	\sqrt{22104 + 12082 \sqrt{2} + 9885 \sqrt{5} + 5400 \sqrt{10} + 
		4 \sqrt{86889872 + 59414261 \sqrt{2} + 38858338 \sqrt{5} + 
			26570859 \sqrt{10}}}$
	
	$e =\sqrt{1 - 31 f^2 + 14 \sqrt{5}  f^2}$

	$b=\displaystyle\frac{1}{248} (-42 - 49 \sqrt{2} - 38 \sqrt{5} - 3 \sqrt{10} + (5666 + 679 \sqrt{2} - 2538 \sqrt{5} - 295 \sqrt{10}) f^2)$
	
	$a=\displaystyle\frac{1}{1264} e(584 + 137 \sqrt{2} + 350 \sqrt{5} + 163 \sqrt{10} + (50886 + 13227 \sqrt{2} - 22792 \sqrt{5} - 5903 \sqrt{10}) f^2)$
	
	$d=\displaystyle\frac{1}{1896}e f (30344+ 21847 \sqrt{2} - 12830 \sqrt{5} -
	10679 \sqrt{10} + 19 (-203414 - 160301 \sqrt{2} + 90964 \sqrt{5} + 71693 \sqrt{10}) f^2) $
	
	$c=\displaystyle\frac{1}{117552} e f (-1719980 - 1190961 \sqrt{2} + 844320 \sqrt{5} + 
	495395 \sqrt{10} + (242831436 + 150419207 \sqrt{2} - 108605296 \sqrt{5}$ 
	\vspace{0.1cm}
	
	\hspace{0.6cm}$-67264185 \sqrt{10}) f^2)$

	\vspace{1cm}

	For $P_{2,46}$, the acute solutions are:
	
	$f=\displaystyle\frac{1}{11\sqrt{2}}\sqrt{3358+1857\sqrt{2}+1498\sqrt{5}+831\sqrt{10}+2\sqrt{2(3804500+2601475\sqrt{2}+1701426\sqrt{5}+1163413\sqrt{10}}}$
	
	$e=\displaystyle\sqrt{1-11f^2+5\sqrt{5}f^2}$
	
	$d=\displaystyle\frac{1}{712}ef(-1188-100\sqrt{5}+547\sqrt{2}+85\sqrt{10}-3f^2(-51628+23074\sqrt{5}-39527\sqrt{2}+17687\sqrt{10}))$
	
	$c=\displaystyle\frac{1}{22072}ef((-5147665-3402751\sqrt{2}+2303503\sqrt{5}+1520761\sqrt{10})f^2+(10517-15966\sqrt{2}+2905\sqrt{5}+935\sqrt{10}))$
	
	$b=\displaystyle\frac{1}{124}(-22-36\sqrt{2}-14\sqrt{5}-6\sqrt{10}-69f^2-313\sqrt{2}f^2+35\sqrt{5}f^2+139\sqrt{10}f^2)$
	
	$a=\displaystyle\frac{1}{356}e(171+313\sqrt{2}+9\sqrt{5}-21\sqrt{10}+(-10376-9016\sqrt{2}+4644\sqrt{5}+4027\sqrt{10})f^2)$
	
\end{remark}

\restoregeometry

\newgeometry{left=0.5cm,right=0.5cm,top=1.5cm,bottom=2cm}
3. Coxeter diagrams for $P_3$

\begin{figure}[H]
	\scalebox{0.3}[0.3]{\includegraphics {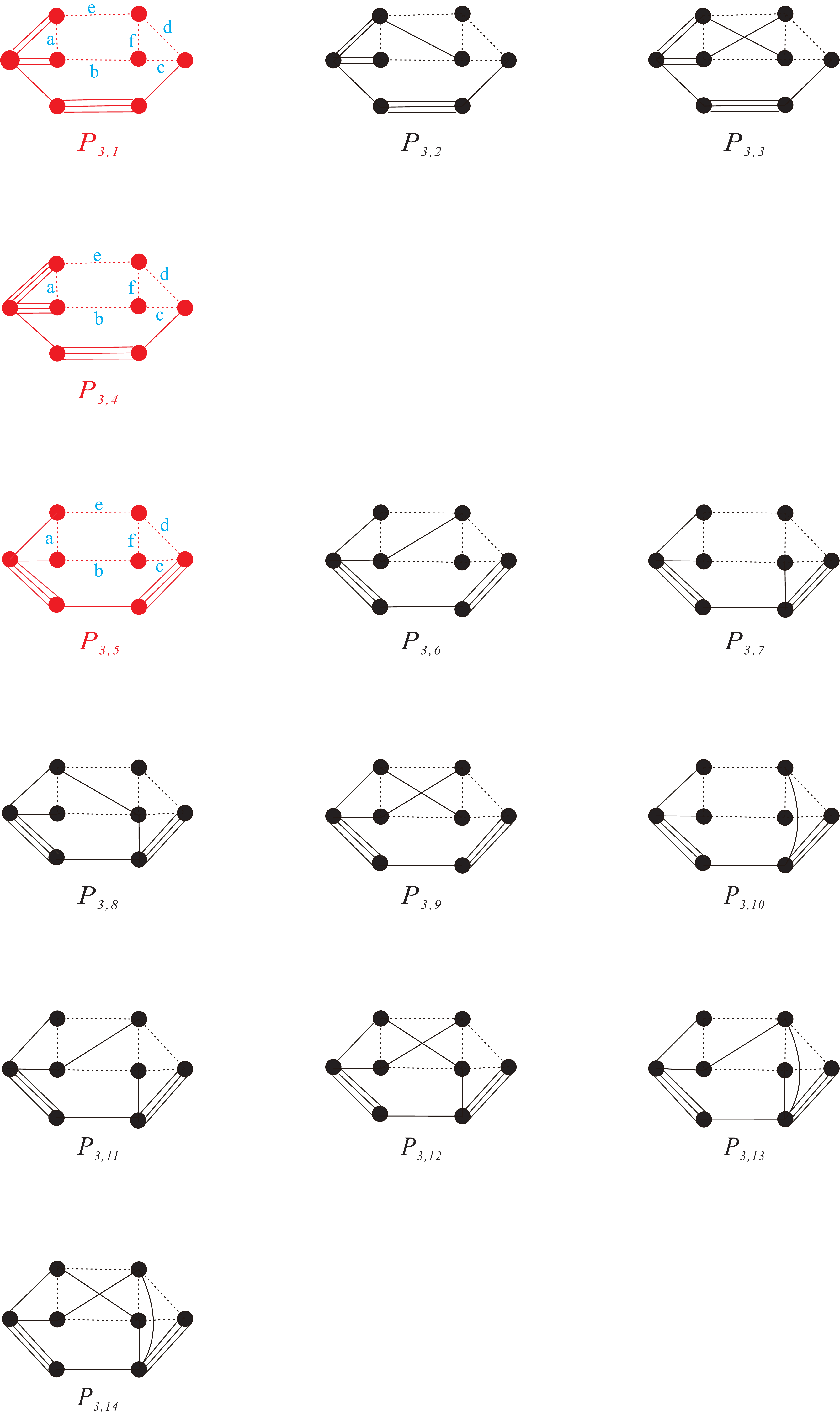}}
	\caption{$P_3$(1/8)} \label{figure:p831}
\end{figure}

\restoregeometry

\begin{figure}[H]
	\scalebox{0.3}[0.3]{\includegraphics {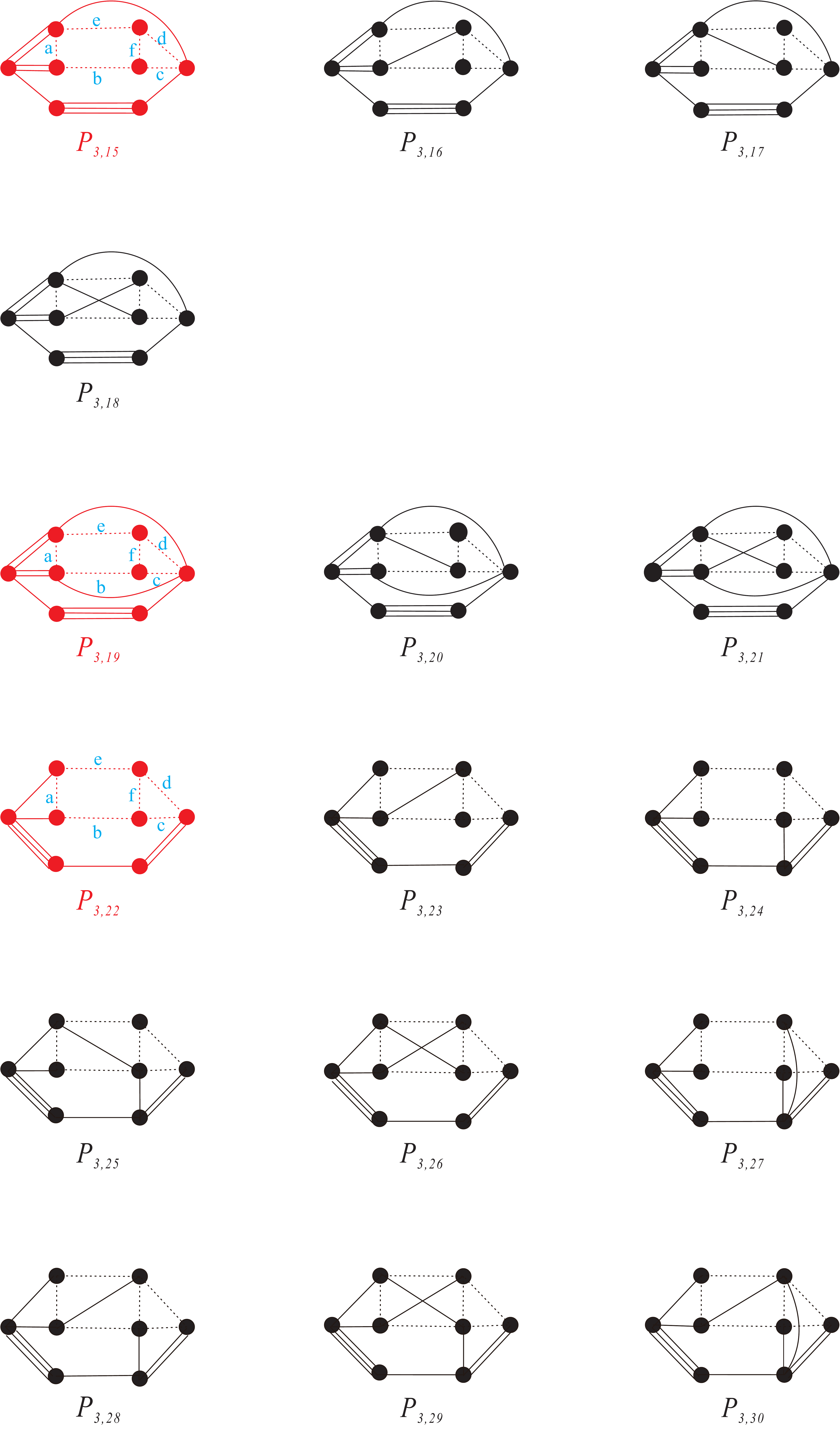}}
	\caption{$P_3$(2/8)} \label{figure:p832}
\end{figure}

\begin{figure}[H]
	\scalebox{0.3}[0.3]{\includegraphics {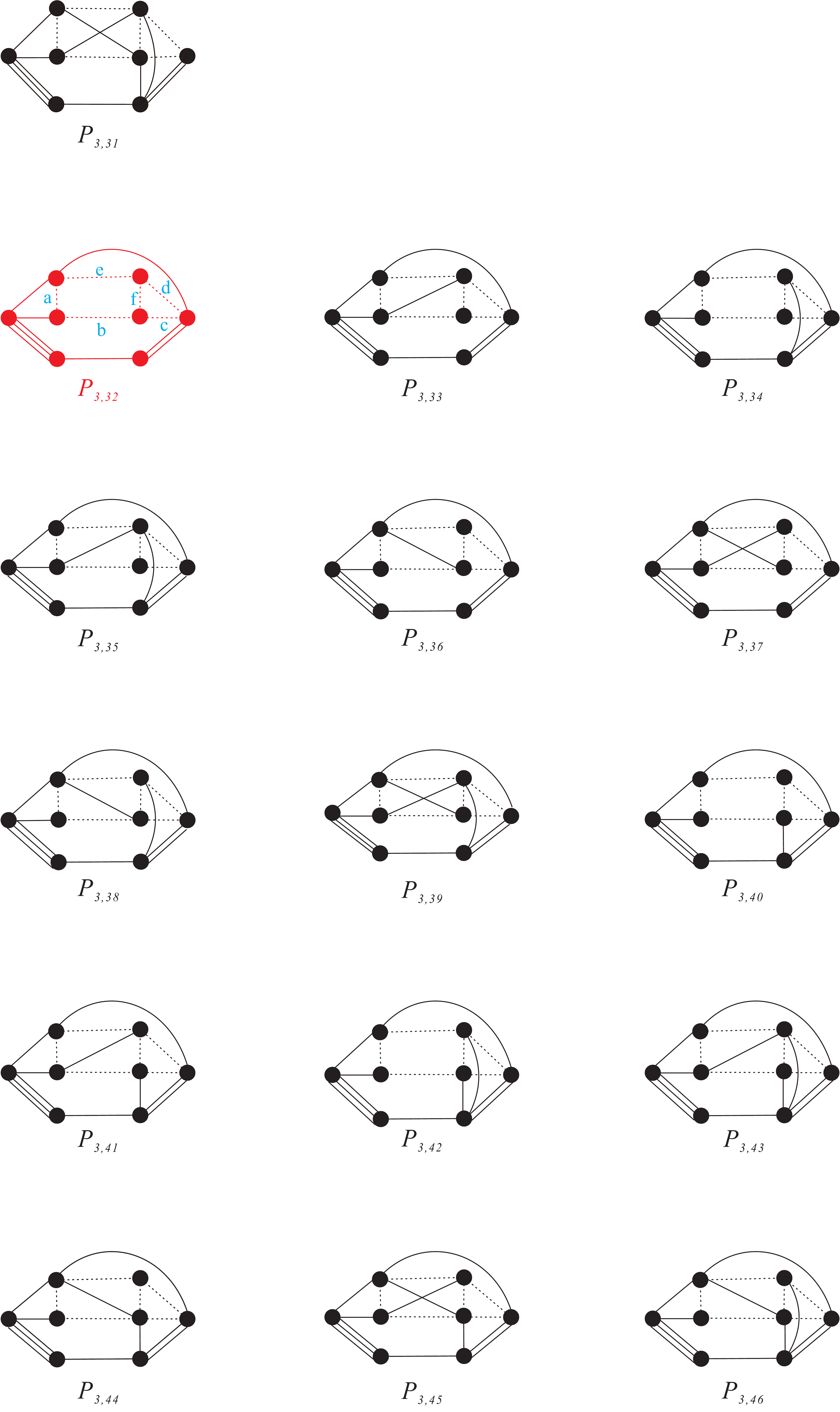}}
	\caption{$P_3$(3/8)} \label{figure:p833}
\end{figure}

\begin{figure}[H]
	\scalebox{0.3}[0.3]{\includegraphics {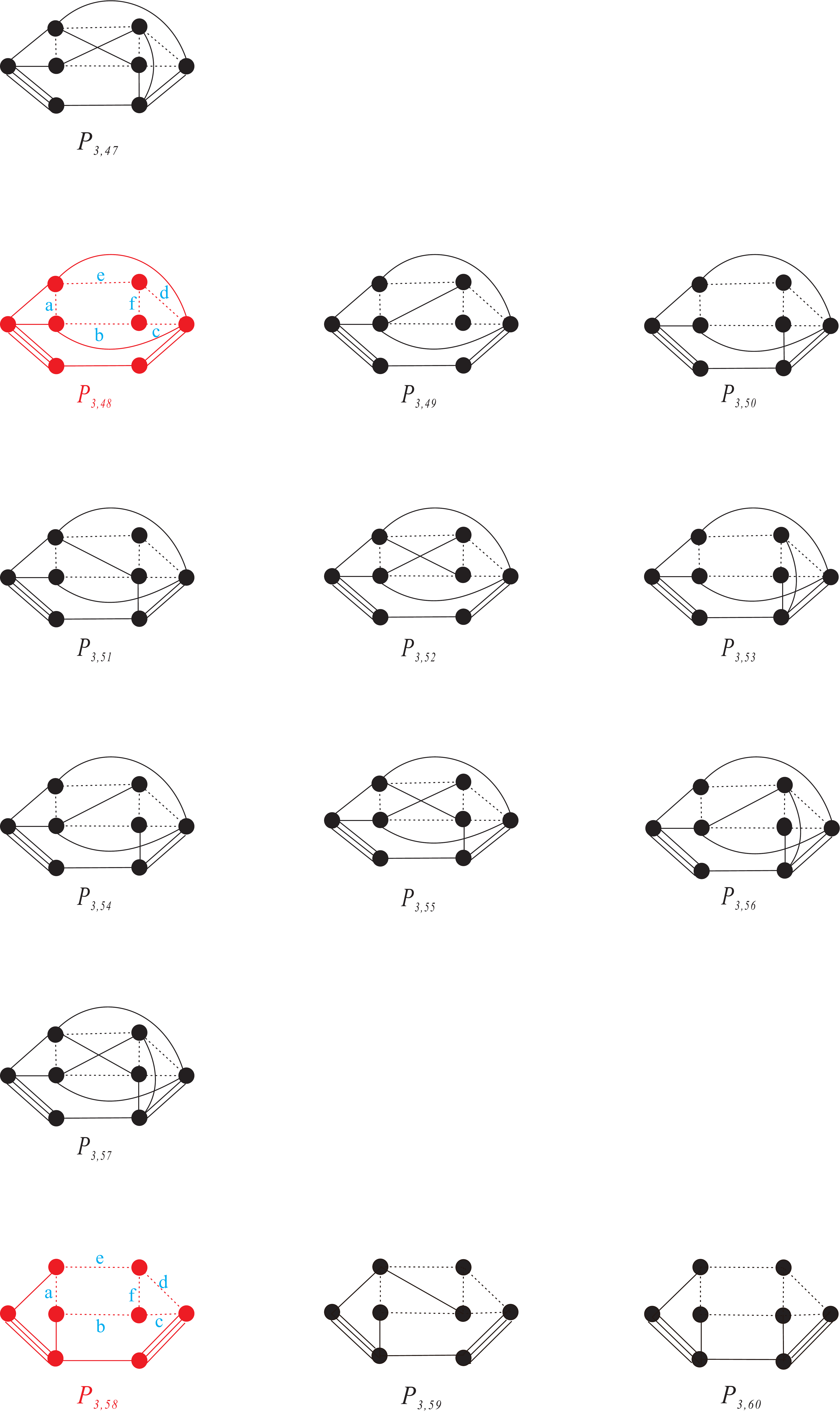}}
	\caption{$P_3$(4/8)} \label{figure:p834}
\end{figure}

\begin{figure}[H]
	\scalebox{0.3}[0.3]{\includegraphics {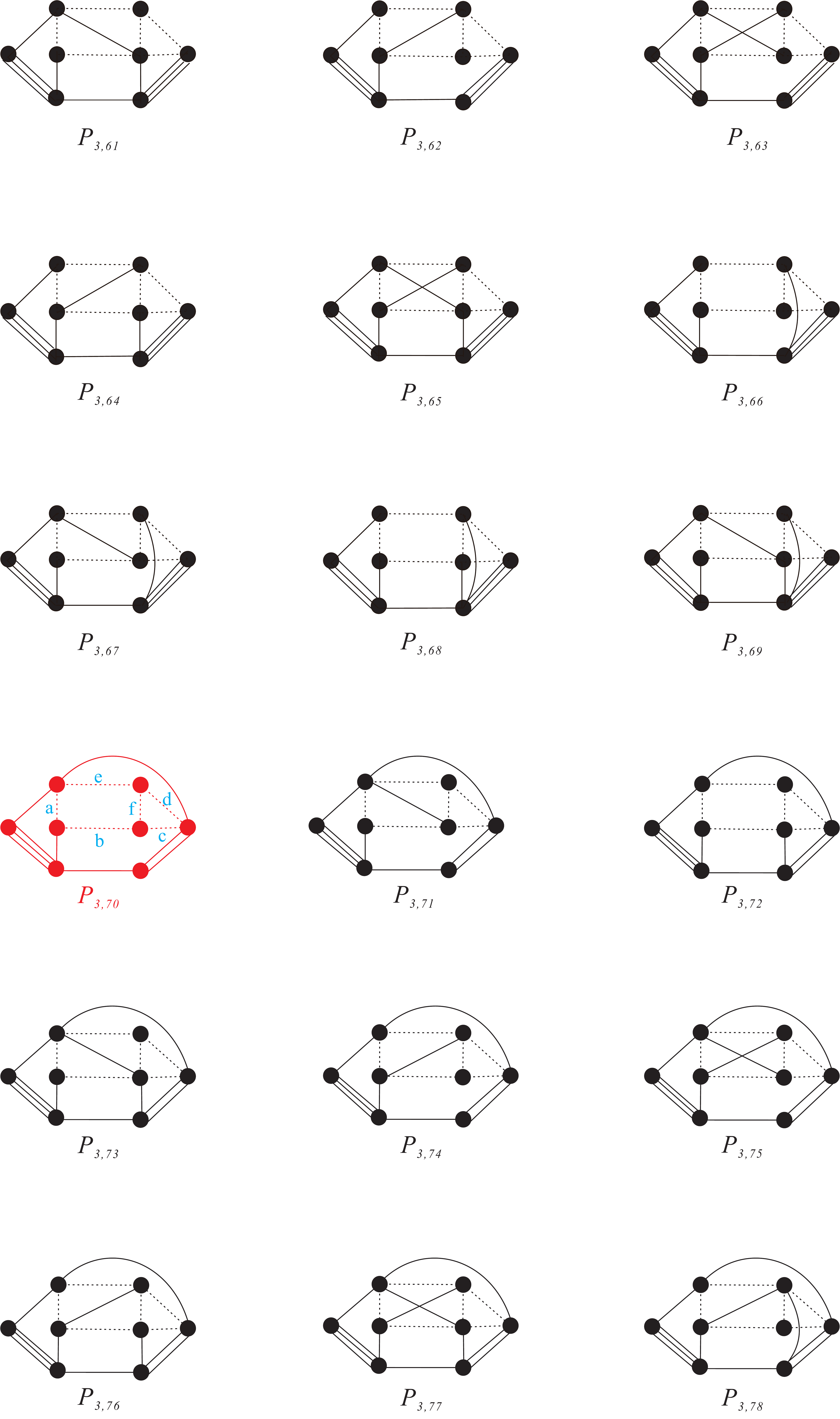}}
	\caption{$P_3$(5/8)} \label{figure:p835}
\end{figure}

\begin{figure}[H]
	\scalebox{0.3}[0.3]{\includegraphics {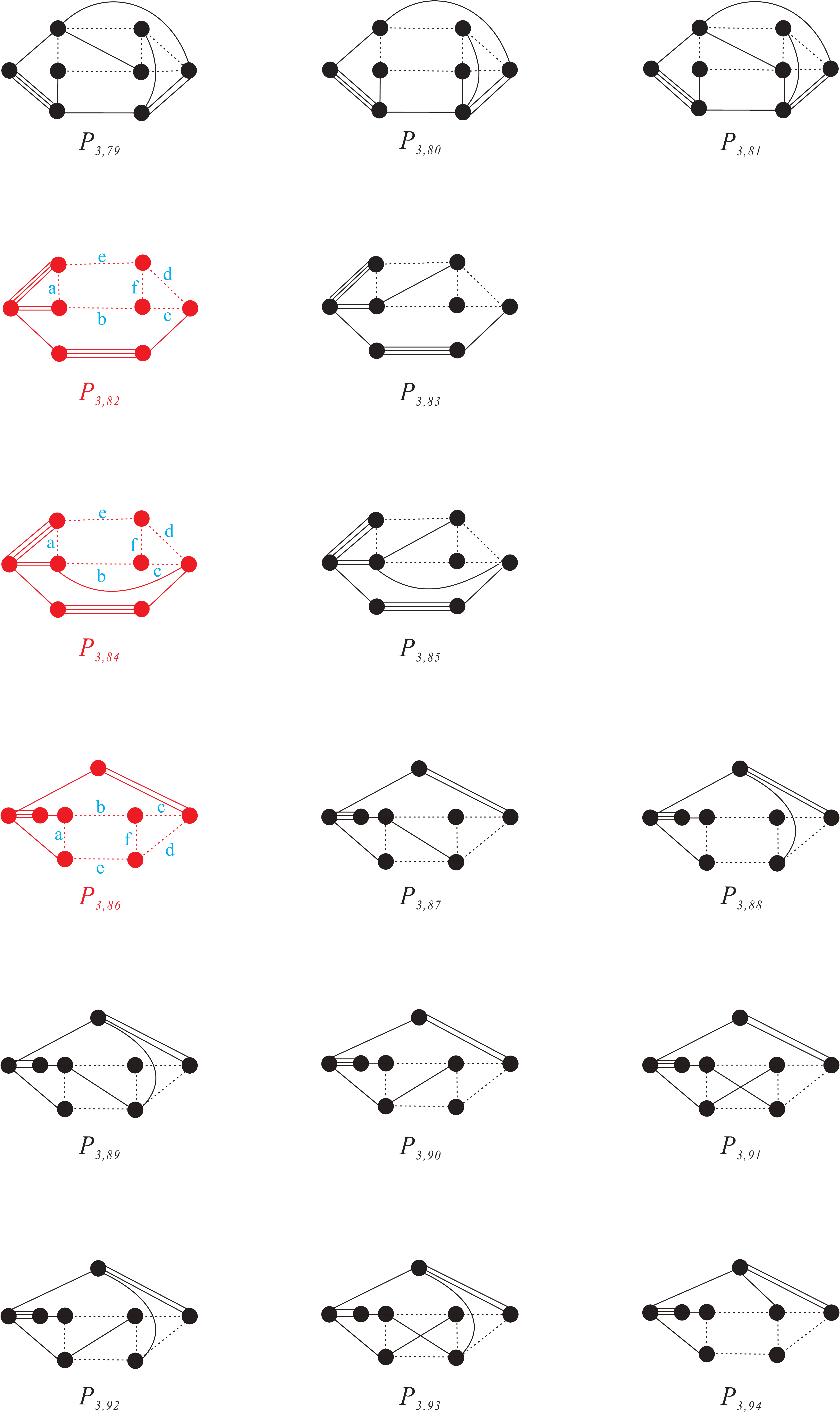}}
	\caption{$P_3$(6/8)} \label{figure:p836}
\end{figure}

\begin{figure}[H]
	\scalebox{0.3}[0.3]{\includegraphics {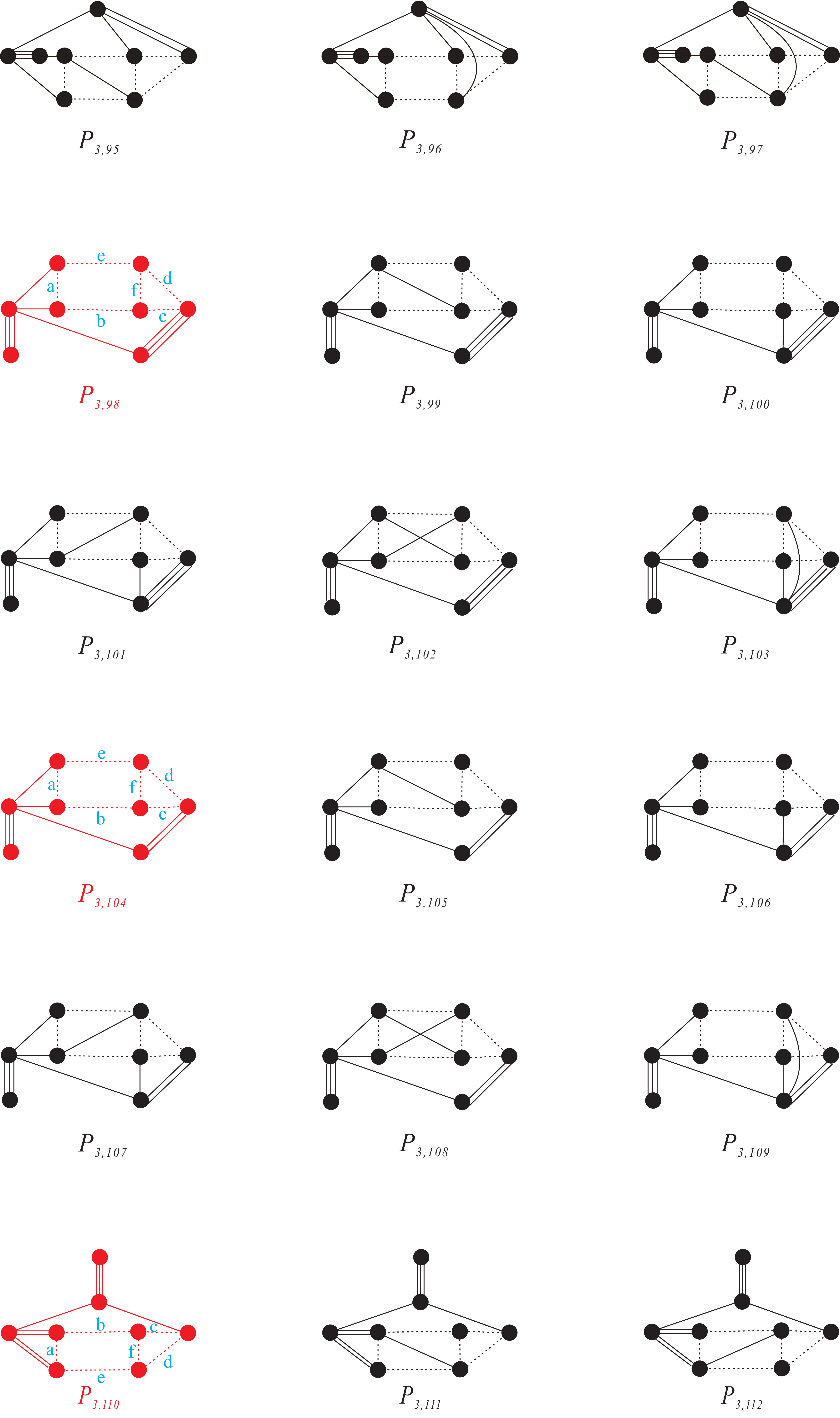}}
	\caption{$P_3$(7/8)} \label{figure:p837}
\end{figure}

\begin{figure}[H]
	\scalebox{0.3}[0.3]{\includegraphics {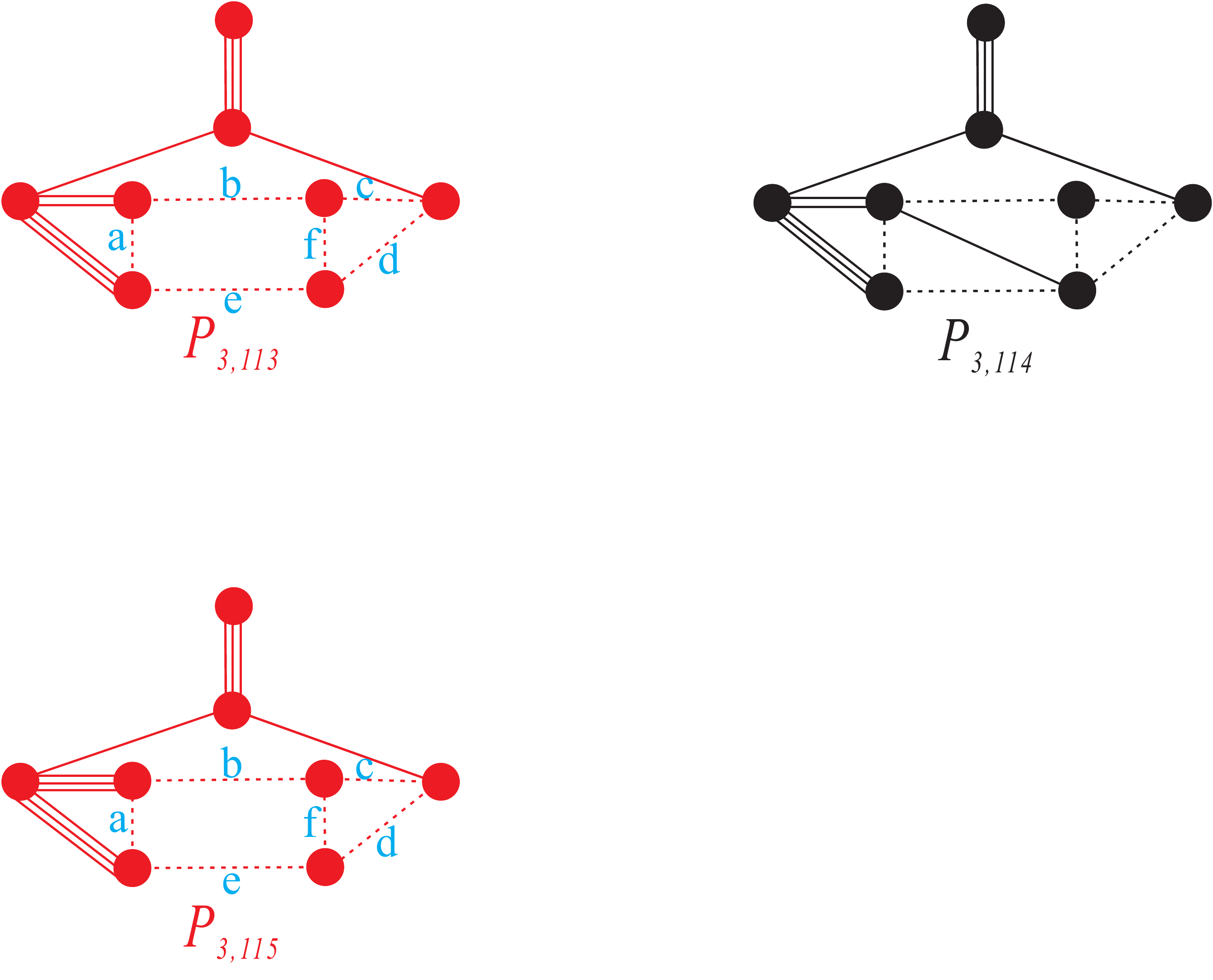}}
	\caption{$P_3$(8/8)} \label{figure:p838}
\end{figure}

\newgeometry{left=0.5cm,right=0.5cm,top=0.2cm,bottom=0.2cm}

\begin{table}[H]
	\resizebox*{\textwidth}{!}{
		\renewcommand{\arraystretch}{1.85}
		\begin{tabular}{|c|c|c|c|}
		\hline
		\multirow{2}{*}{} &$a$  &$b$  &$c$  \\ 
		\cline{2-4} 
		&$d$  &$e$  &$f$  \\ 
		\hline
		
		\multirow{2}{*}{$P_{3,1}$} & $\frac{3}{11}(3+2\sqrt{5})$  & $\frac{2}{11}\sqrt{265+118\sqrt{5}}$  &$\frac{1}{2}\sqrt{5+\sqrt{5}}$  \\ 
		\cline{2-4} 
		&$\frac{1}{2}\sqrt{5+\sqrt{5}}$  &$\frac{2}{11}\sqrt{265+118\sqrt{5}}$  & $\frac{7}{11}(3+2\sqrt{5})$ \\ \hline 
		
		\multirow{2}{*}{$P_{3,4}$} &$\frac{1}{11}(17+4\sqrt{5})$  &$\sqrt{\frac{17208}{2299}+\frac{7688\sqrt{5}}{2299}}  $&$\sqrt{\frac{2}{19}(9+\sqrt{5})}$  \\ 
		\cline{2-4} 
		&$\sqrt{\frac{2}{19}(9+\sqrt{5})}$  &$\sqrt{\frac{17208}{2299}+\frac{7688\sqrt{5}}{2299}}  $  &$\frac{1}{209}(257+208\sqrt{5})$  \\ 
		\hline
		
		\multirow{2}{*}{$P_{3,5}$} &$\frac{1}{19}(17+4\sqrt{5})$  &$\sqrt{\frac{20936}{3971}+\frac{9352\sqrt{5}}{3971}}$  &$\sqrt{\frac{2}{11}(7+\sqrt{5})}$  \\ 
		\cline{2-4} 
		&$\sqrt{\frac{2}{11}(7+\sqrt{5})}$  &$\sqrt{\frac{20936}{3971}+\frac{9352\sqrt{5}}{3971}}  $&$\frac{1}{209}(257+208\sqrt{5})$  \\ \hline
		
		\multirow{2}{*}{$P_{3,15}$} &$3.54605$  &$5.88536$  &$1.9923$  \\ 
		\cline{2-4} 
		&$1.345$  &$5.88536$  &$6.61923$  \\ 
		\hline
		
		\multirow{2}{*}{$P_{3,19}$} &$\frac{1}{11}(8+13\sqrt{2}+9\sqrt{5}+5\sqrt{10})$  &$\frac{1}{11}\sqrt{2074+1462\sqrt{2}+950\sqrt{5}+630\sqrt{10}}$  &$\sqrt{\frac{1}{2}(2+\sqrt{5}+\sqrt{7+3\sqrt{5}})}$  \\ 
		\cline{2-4} 
		&$\sqrt{\frac{1}{2}(2+\sqrt{5}+\sqrt{7+3\sqrt{5}})}$  &$\frac{1}{11}\sqrt{2074+1462\sqrt{2}+950\sqrt{5}+630\sqrt{10}}$  &$\frac{1}{11}(37+12\sqrt{2}+2\sqrt{5}(5+4\sqrt{2}))$  \\ 
		\hline
		
		\multirow{2}{*}{$P_{3,22}$} &$\frac{1}{2}(3+\sqrt{5})$  &$\sqrt{\frac{115}{11}+\frac{51\sqrt{5}}{11}}$  & $\sqrt{\frac{10}{11}+\frac{3\sqrt{5}}{11}}$ \\ 
		\cline{2-4} 
		&$\sqrt{\frac{10}{11}+\frac{3\sqrt{5}}{11}}$  &$\sqrt{\frac{115}{11}+\frac{51\sqrt{5}}{11}}$  &$\frac{7}{11}(3+2\sqrt{5})$  \\ 
		\hline
		
		\multirow{2}{*}{$P_{3,32}$} &$4.41427$  &$6.42282$  &$1.82559$  \\ 
		\cline{2-4} 
		&$1.23245$  &$6.42282$  &$6.61923$  \\ 
		\hline
		
		\multirow{2}{*}{$P_{3,48}$} &$1+\sqrt{5}+\sqrt{7+3\sqrt{5}}$  &$\sqrt{\frac{1}{11}(222+161\sqrt{2}+104\sqrt{5}+67\sqrt{10})}$  &$\sqrt{\frac{1}{22}(19+9\sqrt{5}+2\sqrt{147+65\sqrt{5}})}$  \\ 
		\cline{2-4} 
		&$\sqrt{\frac{1}{22}(19+9\sqrt{5}+2\sqrt{147+65\sqrt{5}})}$  &$\sqrt{\frac{1}{11}(222+161\sqrt{2}+104\sqrt{5}+67\sqrt{10})}$  &$\frac{1}{11}(37+12\sqrt{2}+2\sqrt{5}(5+4\sqrt{2}))$  \\ 
		\hline
		
		\multirow{2}{*}{$P_{3,58}$} & $\frac{1}{76}(25+7\sqrt{5}+8\sqrt{108+31\sqrt{5}})$  &$\frac{1}{19}\sqrt{\frac{1}{22}(30747+13689\sqrt{5}+8\sqrt{29359122+13129762\sqrt{5}})}$  &$\sqrt{\frac{2}{11}(7+\sqrt{5})}$  \\ 
		\cline{2-4} 
		&$\frac{1}{2}\sqrt{\frac{1}{2}(9+\sqrt{5})}$  &$\frac{1}{38}(32+12\sqrt{5}+19\sqrt{\frac{1592}{361}+\frac{711\sqrt{5}}{361}})$  &$\frac{1}{19}\sqrt{\frac{1}{22}(26791+8529\sqrt{5}+8\sqrt{11214278+5015082\sqrt{5}})}$  \\ 
		\hline
		
		\multirow{2}{*}{$P_{3,70}$} &$5.67675$  &$7.83841$  & $1.82559$ \\ 
		\cline{2-4} 
		&$1.14412$  &$6.04875$  &$6.24279$  \\ 
		\hline
		
		\multirow{2}{*}{$P_{3,82}$} &$\frac{1}{11}\sqrt{197+69\sqrt{5}+4\sqrt{2(845+358\sqrt{5})}}$  & $\frac{1}{11}\sqrt{\frac{1}{19}(7186+3203\sqrt{5}+4\sqrt{6387730+2856676\sqrt{5}}}$ &$\sqrt{\frac{2}{19}(9+\sqrt{5})}$  \\ 
		\cline{2-4} 
		&$\frac{1}{2}\sqrt{5+\sqrt{5}}$  &$\frac{1}{22}(26+10\sqrt{5}+11\sqrt{\frac{1385}{121}+\frac{619\sqrt{5}}{121})}$  & $\frac{1}{11}\sqrt{\frac{1}{19}(11128+3983\sqrt{5}+8\sqrt{2439905+1091149\sqrt{5}})}$ \\ 
		\hline
		
		\multirow{2}{*}{$P_{3,84}$} &$3.84864$  & $4.97946$ &$1.08754$  \\ 
		\cline{2-4} 
		&$1.9923$  &$6.47518$  &$5.62568$  \\ 
		\hline
		
		\multirow{2}{*}{$P_{3,86}$} &$\frac{3}{4}+\frac{\sqrt{5}}{4}+\sqrt{\frac{5}{2}+\sqrt{5}}$  & $1+\frac{\sqrt{5}}{2}+\sqrt{\frac{5}{2}+\sqrt{5}}$ &$\frac{1}{2}\sqrt{3+\sqrt{5}}$  \\ 
		\cline{2-4} 
		&$\frac{1}{\sqrt{2-\frac{3}{\sqrt{5}}}}$  &$\sqrt{\frac{1}{22}(173+75\sqrt{5}+4\sqrt{3625+1621\sqrt{5}})}$  & $\sqrt{\frac{1}{22}(123+49\sqrt{5}+4\sqrt{1385+619\sqrt{5}})}$ \\ 
		\hline
		
		\multirow{2}{*}{$P_{3,98}$} &$\frac{3}{19}(8+3\sqrt{5})$  &$\frac{1}{19}\sqrt{2442+1082\sqrt{5}}$  &$\frac{1}{2}\sqrt{\frac{1}{2}(9+\sqrt{5})}$  \\ 
		\cline{2-4} 
		& $\frac{1}{2}\sqrt{\frac{1}{2}(9+\sqrt{5})}$ & $\frac{1}{19}\sqrt{2442+1082\sqrt{5}}$  & $\frac{1}{19}(20+17\sqrt{5})$ \\ 
		\hline
		
		\multirow{2}{*}{$P_{3,104}$} &$2+\sqrt{5}$  &$3+\sqrt{5}$  &$\frac{1}{2}\sqrt{3+\sqrt{5}}$  \\ 
		\cline{2-4} 
		&$\frac{1}{2}\sqrt{3+\sqrt{5}}$  & $3+\sqrt{5}$  & $2+\sqrt{5}$ \\ 
		\hline  
		
		\multirow{2}{*}{$P_{3,110}$} &$\frac{1}{2}(1+\sqrt{5})$  &$\sqrt{7+3\sqrt{5}}$  &$\frac{1}{2}(1+\sqrt{5})$  \\ 
		\cline{2-4} 
		&$\frac{1}{2}(1+\sqrt{5})$  &$\sqrt{7+3\sqrt{5}}$  &$2+\sqrt{5}$  \\ 
		\hline
		
		\multirow{2}{*}{$P_{3,113}$} &$\frac{1}{2}\sqrt{\frac{3}{2}(3+\sqrt{5})+\sqrt{8+3\sqrt{5}}}$  & $\sqrt{\frac{1}{38}(94+40\sqrt{5}+\sqrt{2(8331+3725\sqrt{5})})}$ & $\sqrt{\frac{43}{38}+\frac{9\sqrt{5}}{38}}$ \\ 
		\cline{2-4} 
		& $\frac{1}{2}(1+\sqrt{5})$ &$\frac{1}{2}(2+\sqrt{5}+\sqrt{8+3\sqrt{5})}$  &$\sqrt{\frac{1}{19}(8+3\sqrt{5})(9+2\sqrt{8+3\sqrt{5})}}$  \\ 
		\hline
				
		\multirow{2}{*}{$P_{3,115}$} &$\frac{1}{4}(5+\sqrt{5})$  &$\sqrt{\frac{109}{19}+\frac{48\sqrt{5}}{19}}$  &$\sqrt{\frac{43}{38}+\frac{9\sqrt{5}}{38}}$  \\ 
		\cline{2-4} 
		& $\sqrt{\frac{43}{38}+\frac{9\sqrt{5}}{38}}$ &$\sqrt{\frac{109}{19}+\frac{48\sqrt{5}}{19}}$  &$\frac{1}{19}(20+17\sqrt{5})$   \\ 
		\hline
		\end{tabular}
		}
\end{table}

\restoregeometry

\newgeometry{left=0.5cm,right=0.5cm,top=2cm}
\begin{remark}.
	\vspace{0.5cm}
	
	For $P_{3,15}$, the acute solutions are:
	
	$f=\frac{1}{11}\sqrt{799+406\sqrt{2}+364\sqrt{5}+168\sqrt{10}+2\sqrt{2(211980+145015\sqrt{2}+94834\sqrt{5}+64817\sqrt{10}}}$
	
	$e=\frac{1}{2}\sqrt{-1-\sqrt{5}+f^2+\sqrt{5}f^2}$
	
	$d={\displaystyle\frac{1}{496}}e(1248+161\sqrt{2}+240\sqrt{5}+67\sqrt{10}+(140+81\sqrt{2}-3\sqrt{5}(36+7\sqrt{2}))f^2)$
	
	$c={\displaystyle\frac{1}{1488}}ef(4396+1661\sqrt{2}+436\sqrt{5}+511\sqrt{10}+(-896-1327\sqrt{2}+304\sqrt{5}+595\sqrt{10})f^2)$
	
	$b={\displaystyle\frac{1}{186}}ef(248+211\sqrt{2}-31\sqrt{5}-94\sqrt{10}+(-1000-443\sqrt{2}+443\sqrt{5}+200\sqrt{10})f^2)$
	
	$a={\displaystyle\frac{1}{248}}(-259+13\sqrt{2}-7\sqrt{5}-3\sqrt{10}+(37+69\sqrt{2}+\sqrt{5}-35\sqrt{10})f^2)$
	
	\vspace{1cm}

	For $P_{3,32}$, the acute solutions are:
	
	$f=\frac{1}{11}\sqrt{799+406\sqrt{2}+364\sqrt{5}+168\sqrt{10}+2\sqrt{2(211980+145015\sqrt{2}+94834\sqrt{5}+64817\sqrt{10}}}$
	
	$e={\displaystyle\frac{1}{2\sqrt{2}}}\sqrt{-1-3\sqrt{5}+f^2+3\sqrt{5}f^2}$
	
	$d={\displaystyle\frac{1}{2728}}e(4968+731\sqrt{2}+1464\sqrt{5}+315\sqrt{10}+11(20-28\sqrt{5}-3\sqrt{2}(-7+\sqrt{5}))f^2)$
	
	$c={\displaystyle\frac{1}{8184}}ef(16476+7091\sqrt{2}+3724\sqrt{5}+2619\sqrt{10}+11(-216-287\sqrt{2}+56\sqrt{5}+129\sqrt{10})f^2)$
	
	$b={\displaystyle\frac{1}{186}}ef(248+211\sqrt{2}-31\sqrt{5}-94\sqrt{10}+(-1000-443\sqrt{2}+443\sqrt{5}+200\sqrt{10})f^2)$
	
	$a={\displaystyle\frac{1}{496}}(-517+53\sqrt{2}-19\sqrt{5}-17\sqrt{10}+(127+329\sqrt{2}-15\sqrt{5}-157\sqrt{10})f^2)$
	
	\vspace{1cm}
	
	For $P_{3,70}$, the acute solutions are:
	
	$f={\displaystyle\frac{1}{\sqrt{22}}}\sqrt{135+64\sqrt{2}+57\sqrt{5}+28\sqrt{10}+4\sqrt{2743+1884\sqrt{2}+1231\sqrt{5}+838\sqrt{10}}}$
	
	$e={\displaystyle\frac{1}{2\sqrt{2}}}\sqrt{-1-3\sqrt{5}+f^2+3\sqrt{5}f^2}$
	
	$d={\displaystyle\frac{1}{23684}}e(36906+2881\sqrt{2}+11886\sqrt{5}+847\sqrt{10}+(4504+3808\sqrt{2}-4044\sqrt{5}-874\sqrt{10})f^2$
	
	$c={\displaystyle\frac{1}{2913132}}((2545015+877512\sqrt{2}+2315748\sqrt{5}+1194982\sqrt{10})ef+(2442218+1364747\sqrt{2}-1270198\sqrt{5}-$\vspace{0.1cm}
	
	\hspace{0.6cm}$586721\sqrt{10})ef^3)$
	
	$b={\displaystyle\frac{1}{3813}}((6454+4099\sqrt{2}-743\sqrt{5}-1864\sqrt{10})ef+\frac{1}{2}(-8619-7529\sqrt{2}+3601\sqrt{5}+3483\sqrt{10})ef^3)$
	
	$a={\displaystyle\frac{1}{248}}(-97-15\sqrt{2}-73\sqrt{5}+13\sqrt{10}+(21-53\sqrt{2}+19\sqrt{5}+17\sqrt{10})f^2)$

	\newpage
	
	For $P_{3,84}$, the acute solutions are:
	
	$f={\displaystyle\frac{1}{11\sqrt{19}}}\sqrt{11104+5702\sqrt{2}+5133\sqrt{5}+2276\sqrt{10}+8\sqrt{4842088+3311169\sqrt{2}+2165554\sqrt{5}+1480687\sqrt{10}}}$
	
	$e={\displaystyle\frac{1}{2}}\sqrt{-1-2\sqrt{5}+f^2+2\sqrt{5}f^2}$
	
	$	d={\displaystyle\frac{1}{24216}}e(-17104+5893\sqrt{2}-3602\sqrt{5}-20597\sqrt{10}+(-115066-24517\sqrt{2}+50384\sqrt{5}+12469\sqrt{10})f^2)$
	
	c=${\displaystyle\frac{1}{1937062056}}ef(9466659574+914636627\sqrt{2}-3948733168\sqrt{5}-1618127971\sqrt{10}+(-23059580156-23947802279\sqrt{2}+$\vspace{0.1cm}
	
	\hspace{0.6cm}$10158748906\sqrt{5}+10850958719\sqrt{10})f^2)$
	
	$b={\displaystyle\frac{1}{9918884}}ef(-31360330-1617047\sqrt{2}+13795674\sqrt{5}+3094731\sqrt{10}+(61953682+76754869\sqrt{2}-5\sqrt{5}(5508882+$\vspace{0.1cm}
	
	\hspace{0.6cm} $6899345\sqrt{2}))f^2)$
	
	$a={\displaystyle\frac{1}{248}}(-57+\frac{177}{\sqrt{2}}-265\sqrt{\frac{5}{2}}+37\sqrt{5}+(-679-\frac{117}{\sqrt{2}}+89\sqrt{\frac{5}{2}}+295\sqrt{5})f^2)$
\end{remark}
\restoregeometry

\newpage
\begin{figure}[H]
	\scalebox{0.3}[0.3]{\includegraphics {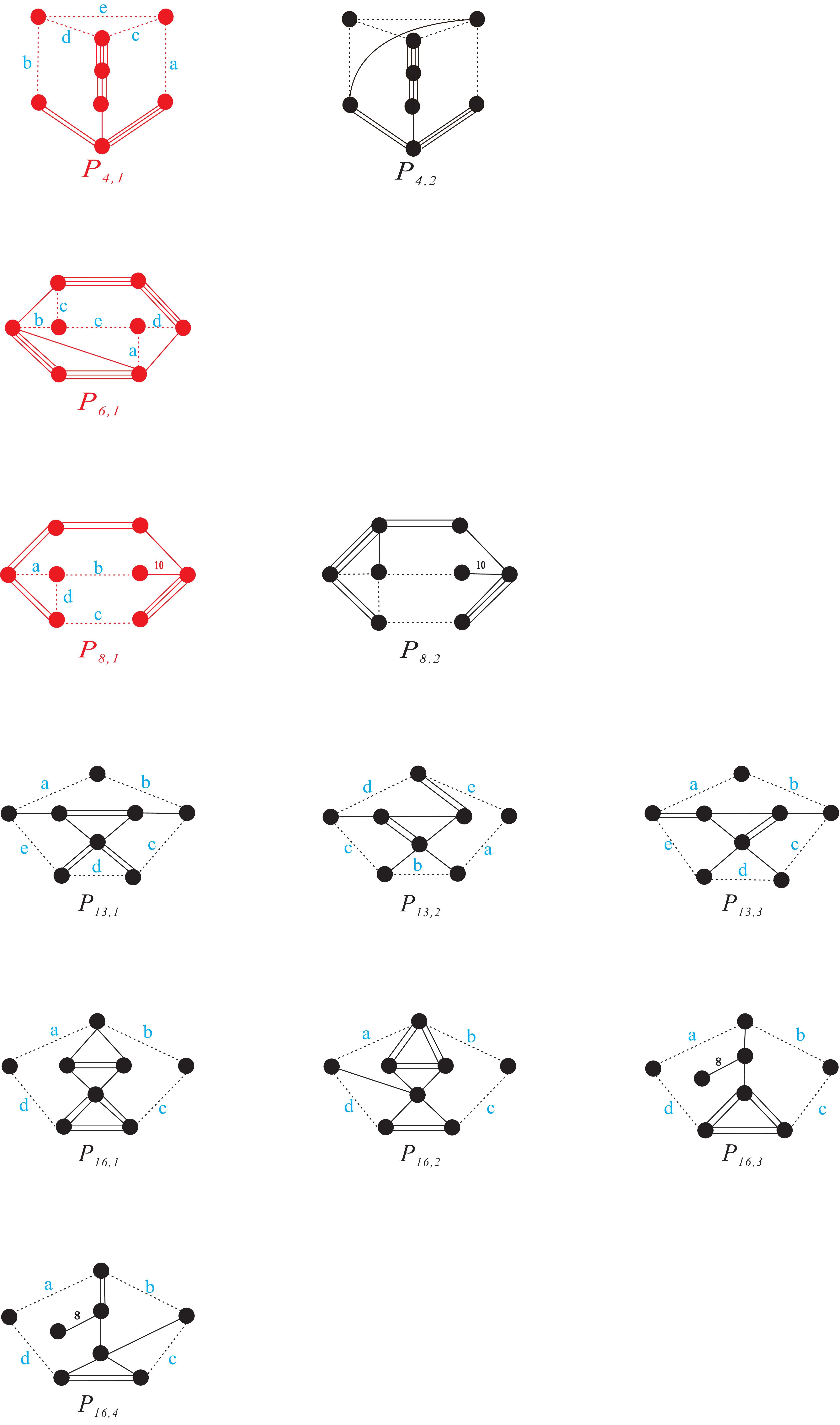}}
	\caption{$P_4,P_6,P_8,P_{13},P_{16}$} \label{figure:p4681316}
\end{figure}

\newgeometry{left=0.5cm,right=0.5cm,top=2cm,bottom=2cm}
\bgroup
\everymath{\displaystyle}

\begin{table}[H]
	\resizebox*{19cm}{!}{
		\renewcommand{\arraystretch}{2.7}		
		\begin{tabular}{|c |c | c | c | c | c |} 
				
					\hline
					&   a &   b &  c &  d &  e  \\ \hline
					$P_{4,1}$ & $\frac{ 1}{ 2} \sqrt{3(3+\sqrt{5})}$ & $\sqrt{\frac{ 3}{ 19} (8+3\sqrt{5})}$ & $\frac{ 1}{ 2} \sqrt{7 + 3\sqrt{5}}$ & $\sqrt{\frac{ 2}{19} (8 +3\sqrt{5})}$ & $\sqrt{\frac{ 5}{ 19} (8 +3\sqrt{5})}$\\[7pt] \hline
					$P_{4,2}$ &  $\frac{ 1}{ 4} \sqrt{57 + 23\sqrt{5} + 4\sqrt{30(9 + 4\sqrt{5})}}$ & $\sqrt{\frac{ 24}{ 19} + \frac{ 9\sqrt{5}}{19}}$ & $\frac{ 1}{ 4} (2(2+\sqrt{5})+\sqrt{21+9\sqrt{5}})$ & $\sqrt{\frac{ 16}{ 19} + \frac{ 6\sqrt{5}}{19}}$ &  $\frac{ 1}{ 2} \sqrt{\frac{ 1}{ 19} (249+91\sqrt{5} +60\sqrt{6} +32\sqrt{30})}$ \\[7pt] 
					\hline
					\specialrule{0.5em}{0.5pt}{0.5pt} \multicolumn{6}{l}{}\\
					\hline
					$P_{6,1}$ & $  \sqrt{\frac{ 5}{ 19} (8 +3\sqrt{5})} $ &   $   \sqrt{\frac{ 5}{ 19} (8 +3\sqrt{5})} $  &  $    \sqrt{\frac{ 5}{ 38} (9 +\sqrt{5})} $ & $    \sqrt{\frac{ 5}{ 38} (9 +\sqrt{5})} $  & $     \frac{ 25}{ 19} (9 +\sqrt{5})-9 $ \\[7pt]       
					\hline
				    \specialrule{0.5em}{0.5pt}{0.5pt} \multicolumn{6}{l}{}\\
					\hline
					$P_{8,1}$ & $\frac{ 1}{ 2} \sqrt{5+\sqrt{5}}$ & $\frac{ 1}{ 2} (3+\sqrt{5})$ & $\frac{ 1}{ 2} (1+\sqrt{5})$ & $\sqrt{5+2\sqrt{5}}$& \\[7pt]                    
					\hline
					$P_{8,2}$ & $\frac{ 1}{2} (1+\sqrt{5+2\sqrt{5}})$ & \footnotesize $\sqrt{7+3\sqrt{5} + \sqrt{85+38\sqrt{5}}}$ & $\frac{ 1}{2} (1+\sqrt{5})$ & \footnotesize $\frac{ 1}{2} \sqrt{43+19\sqrt{5} + 2\sqrt{890 + 398\sqrt{5}}}$& \\[7pt]                   
					\hline
					\specialrule{0.5em}{0.5pt}{0.5pt} \multicolumn{6}{l}{}\\
					\hline
					$P_{13,1}$ &  $\frac{ 1}{ 2} \sqrt{3+\sqrt{2}} $ &   $\frac{ 1}{ 2} \sqrt{3+\sqrt{2}} $  & $\frac{ 1}{ 2} (2+\sqrt{2}) $ & $1+\sqrt{2}$  & $\frac{ 1}{ 2} (2+\sqrt{2}) $  \\[7pt]                    
					\hline
					$P_{13,2}$ &  $\sqrt{\frac{ 17}{ 23} + \frac{ 8\sqrt{2}}{ 23}}$ &   $\frac{ 1}{ 2} +\frac{ 1}{ \sqrt{2}}$  &  $\frac{ 3}{ 2}+ \frac{ 1}{ \sqrt{2}}$ & $1 + \frac{ 1}{ \sqrt{2}}$  & $\sqrt{\frac{ 17}{ 23} + \frac{ 8\sqrt{2}}{ 23}}$ \\[7pt]                    
					\hline
					$P_{13,3}$ & $\sqrt{\frac{ 17}{ 23} + \frac{ 8\sqrt{2}}{ 23}}$ &   $\sqrt{\frac{ 17}{ 23} + \frac{ 8\sqrt{2}}{ 23}}$ &  $1+\frac{ 1}{ 2\sqrt{2}}$ & $1+\frac{3}{ 2\sqrt{2}}$  & $1+\frac{ 1}{\sqrt{2}}$  \\[7pt]                    
					\hline
					\specialrule{0.5em}{0.5pt}{0.5pt} \multicolumn{6}{l}{}\\
					\hline
					
					$P_{16,1}$ & $\sqrt{1+\frac{ 1}{  \sqrt{2}}}$ & $\sqrt{1+ \frac{1}{ \sqrt{2}}}$ & $\sqrt{1+\frac{ 1}{ \sqrt{2}}}$ & $\sqrt{1+\frac{ 1}{ \sqrt{2}}}$& \\[7pt]                    
					\hline
					
					$P_{16,2}$ &  $2+\sqrt{2}$ &   $\sqrt{\frac{ 13}{ 7} + \frac{ 9 \sqrt{2}}{ 7}}$  &  $\sqrt{\frac{  2}{ 7}(3+\sqrt{2})}$ & $1+\frac{ 1}{ \sqrt{2}}$& \\[7pt]  
					                  
					\hline
					$P_{16,3}$ & $\sqrt{1+\frac{ 1}{ \sqrt{2}}}$ & $\sqrt{1+\frac{ 1}{ \sqrt{2}}}$ & $\sqrt{1+\frac{ 1}{ \sqrt{2}}}$ & $\sqrt{1+\frac{ 1}{ \sqrt{2}}}$& \\[7pt] 
					     
					\hline
					$P_{16,4}$ & $\sqrt{\frac{ 13}{ 7} + \frac{ 9 \sqrt{2}}{ 7}}$ & $2+\sqrt{2}$ & $1+\frac{ 1}{ \sqrt{2}}$  &  $\sqrt{\frac{2}{ 7}(3+\sqrt{2})}$& \\[7pt]                    
					\hline
                	\specialrule{0.5em}{0.5pt}{0.5pt} \multicolumn{6}{l}{}\\
					
				\end{tabular}
			}
		\end{table}
\egroup
\restoregeometry

\newpage
\begin{figure}[H]
	\scalebox{0.3}[0.3]{\includegraphics {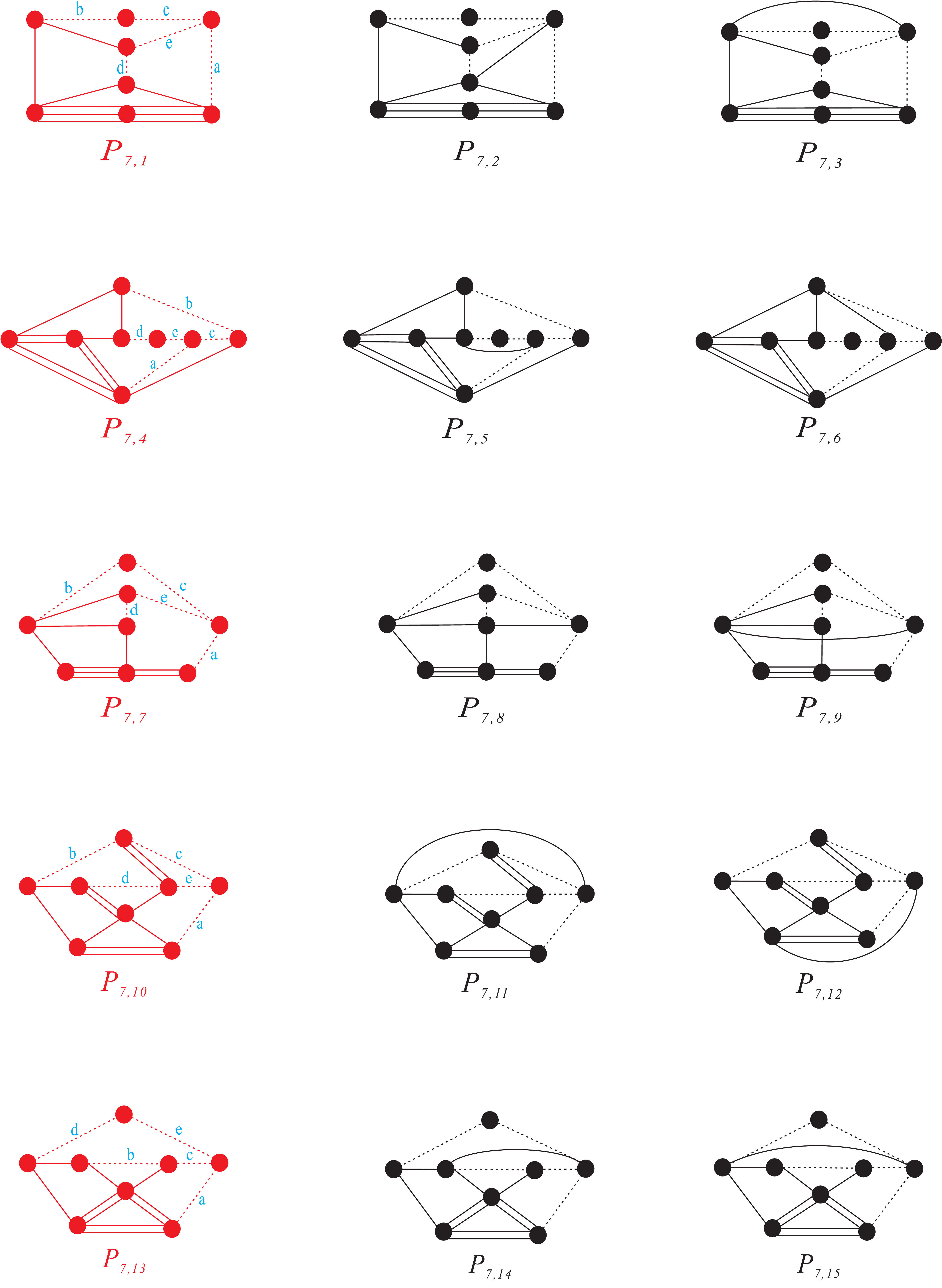}}
	\caption{$P_{7}$} \label{figure:p7}
\end{figure}

\newgeometry{left=0.5cm,right=0.5cm,top=2cm,bottom=2cm}
\bgroup
\everymath{\displaystyle}

\begin{table}[H]
	\renewcommand{\arraystretch}{2}
	\resizebox{17cm}{!}{ %
		\begin{tabular}{|l|l|l|l|}
			\hline
			\multirow{2}{*}{} &a  &b  &c  \\ \cline{2-4} 
			&d  &e  &  \\ \hline
			\multirow{2}{*}{$P_{7,1}$} &$\sqrt{\frac{23}{8}+\frac{9\sqrt{5}}{8}}$ & $\frac{1}{2}\sqrt{\frac{1}{2}(6+\sqrt{5})}$ & $\frac{1}{4}\sqrt{5(3+\sqrt{5})}$  \\ \cline{2-4} 
			&$\frac{1}{2}(1+\sqrt{5})$ & $\sqrt{\frac{23}{8}+\frac{9\sqrt{5}}{8}}$  &  \\ \hline
			
			\multirow{2}{*}{$P_{7,2}$} & $\frac{1}{8}(3(1+\sqrt{5})+\sqrt{206+86\sqrt{5}})$ &   $\frac{1}{2}\sqrt{\frac{1}{2}(6+\sqrt{5})}$  &  $\frac{1}{4}\sqrt{22+\frac{17\sqrt{5}}{2}+\sqrt{655+290\sqrt{5}}}$  \\ \cline{2-4} 
			&$\frac{1}{2}(1+\sqrt{5})$  &  $\frac{1}{8}(5+\sqrt{5}+\sqrt{206+86\sqrt{5}})$  &  \\ \hline
			
			\multirow{2}{*}{$P_{7,3}$} &  $\frac{1}{8}(5+3\sqrt{5}+\sqrt{206+86\sqrt{5}})$ & $\frac{1}{2}\sqrt{\frac{1}{2}(6+\sqrt{5})}$ & $\frac{1}{4}\sqrt{\frac{1}{2}(39+11\sqrt{5}+\sqrt{470+130\sqrt{5}})}$  \\ \cline{2-4} 
			& $\frac{1}{2}(1+\sqrt{5})$ &  $\frac{1}{8}(5+3\sqrt{5}+\sqrt{206+86\sqrt{5}})$  &  \\ \hline
			
			\multirow{2}{*}{$P_{7,4}$} &  $1+\sqrt{2} $ &   $1+\frac{1}{2\sqrt{2}} $  & $1+\frac{1}{\sqrt{2}} $  \\ \cline{2-4} 
			& $\sqrt{\frac{9}{14}+\frac{2\sqrt{2}}{7}}$  & $\sqrt{\frac{1}{7}(9+4\sqrt{2})} $  &  \\ \hline
						
			\multirow{2}{*}{$P_{7,5}$} &  $1+\frac{1}{\sqrt{2}}+\sqrt{\frac{26}{7}+\frac{18\sqrt{2}}{7}}$ &  $1+\frac{1}{2\sqrt{2}} $ &  $\frac{1}{14}(7+7\sqrt{2}+2\sqrt{91+63\sqrt{2}})$  \\ \cline{2-4} 
			& $\sqrt{\frac{9}{14}+\frac{2\sqrt{2}}{7}}$ & $\sqrt{\frac{2}{7}(6+3\sqrt{2}+2\sqrt{5+3\sqrt{2}})}$  &  \\ \hline
						
			\multirow{2}{*}{$P_{7,6}$} &  $1+\frac{1}{\sqrt{2}}+\sqrt{\frac{26}{7}+\frac{18\sqrt{2}}{7}}$ &   $1+\frac{1}{2\sqrt{2}} $  &  $\frac{1}{2}(\sqrt{2}+\sqrt{\frac{52}{7}+\frac{36\sqrt{2}}{7}})$  \\ \cline{2-4} 
			& $\sqrt{\frac{9}{14}+\frac{2\sqrt{2}}{7}}$  & $\sqrt{\frac{1}{7}(13+8\sqrt{2}+2\sqrt{54+38\sqrt{2}})}$ &  \\ \hline

			\multirow{2}{*}{$P_{7,7}$} & $\sqrt{\frac{1}{31}(29+10\sqrt{5})}$ &   $\frac{1}{2}\sqrt{4+\sqrt{5}}$ &  $\sqrt{\frac{1}{31}(23+9\sqrt{5})}$  \\ \cline{2-4} 
			& $\frac{1}{2}(1+\sqrt{5})$  & $\sqrt{\frac{2}{31}(29+10\sqrt{5})}$  &  \\ \hline
			
			\multirow{2}{*}{$P_{7,8}$} & $\frac{1}{62}\sqrt{4347+1577\sqrt{5}+8\sqrt{11(21727+9289\sqrt{5})}}$ &   $\frac{1}{2}\sqrt{4+\sqrt{5}}$ &  $\frac{1}{62}\sqrt{\frac{1}{2}(7899+2851\sqrt{5}+48\sqrt{29048+12674\sqrt{5}})}$  \\ \cline{2-4} 
			& $\frac{1}{2}(1+\sqrt{5})$  & $\frac{1}{124}(49+3\sqrt{5}+8\sqrt{448+178\sqrt{5}}$  &  \\ \hline
			
			\multirow{2}{*}{$P_{7,9}$} & $\frac{1}{31}\sqrt{1079+433\sqrt{5}+4\sqrt{75257+33535\sqrt{5}}}$ &   $\frac{1}{2}\sqrt{4+\sqrt{5}}$ &  $\frac{1}{31}\sqrt{968+213\sqrt{5}+4\sqrt{6613+787\sqrt{5}}}$  \\ \cline{2-4} 
			& $\frac{1}{2}(1+\sqrt{5})$  & $\frac{1}{62}(22+14\sqrt{5}+31\sqrt{\frac{7168}{961}+\frac{2848\sqrt{5}}{961}})$  &  \\ \hline

			\multirow{2}{*}{$P_{7,10}$} & $\sqrt{\frac{13}{7}+\frac{9\sqrt{2}}{7}}$ &   $\frac{1}{2}+\frac{1}{\sqrt{2}}$ &  $2\sqrt{\frac{1}{7}(3+\sqrt{2})}$  \\ \cline{2-4} 
			& $1+\frac{1}{\sqrt{2}}$  & $\sqrt{\frac{45}{7}+\frac{29\sqrt{2}}{7}}$  &  \\ \hline
			
			\multirow{2}{*}{$P_{7,11}$} & $2+\sqrt{2}$ &   $\frac{1}{2}+\frac{1}{\sqrt{2}}$ &  $1+\sqrt{2}$  \\ \cline{2-4} 
			& $1+\frac{1}{\sqrt{2}}$  & $3+\frac{5}{\sqrt{2}}$ &  \\ \hline
			
			\multirow{2}{*}{$P_{7,12}$} & $\frac{1}{14}(25+13\sqrt{2})$ &  $\frac{1}{2}+\frac{1}{\sqrt{2}}$  &  $\frac{8}{7}+\frac{17}{7\sqrt{2}}$  \\ \cline{2-4} 
			& $1+\frac{1}{\sqrt{2}}$  & $\frac{1}{14}(41+37\sqrt{2})$   &  \\ \hline

			\multirow{2}{*}{$P_{7,13}$} & $1+\sqrt{2}$ &  $\frac{1}{2}+\frac{1}{\sqrt{2}}$ &  $1+\sqrt{2}$  \\ \cline{2-4} 
			& $\sqrt{\frac{9}{14}+\frac{2\sqrt{2}}{7}}$  & $\sqrt{\frac{9}{7}+\frac{4\sqrt{2}}{7}}$  &  \\ \hline
			\multirow{2}{*}{$P_{7,14}$} & $1+\frac{1}{\sqrt{2}}+\sqrt{\frac{26}{7}+\frac{18\sqrt{2}}{7}}$ &   $\frac{1}{2}+\frac{1}{\sqrt{2}}$ &  $\frac{1}{2}+\frac{1}{\sqrt{2}}+\sqrt{\frac{26}{7}+\frac{18\sqrt{2}}{7}}$   \\ \cline{2-4} 
			& $\sqrt{\frac{9}{14}+\frac{2\sqrt{2}}{7}}$  & $\sqrt{\frac{1}{7}(13+8\sqrt{2}+2\sqrt{54+38\sqrt{2}})}$ &  \\ \hline
			\multirow{2}{*}{$P_{7,15}$} & $1+\frac{1}{\sqrt{2}}+\sqrt{\frac{26}{7}+\frac{18\sqrt{2}}{7}}$ &   $\frac{1}{2}+\frac{1}{\sqrt{2}}$  &  $1+\frac{1}{\sqrt{2}}+\sqrt{\frac{26}{7}+\frac{18\sqrt{2}}{7}}$   \\ \cline{2-4} 
			& $\sqrt{\frac{9}{14}+\frac{2\sqrt{2}}{7}}$  & $\sqrt{\frac{2}{7}(6+3\sqrt{2}+2\sqrt{5+3\sqrt{2}})}$   &  \\ \hline
		\end{tabular}
	}
\end{table}
\egroup
\restoregeometry

\begin{figure}[H]
	\scalebox{0.3}[0.3]{\includegraphics {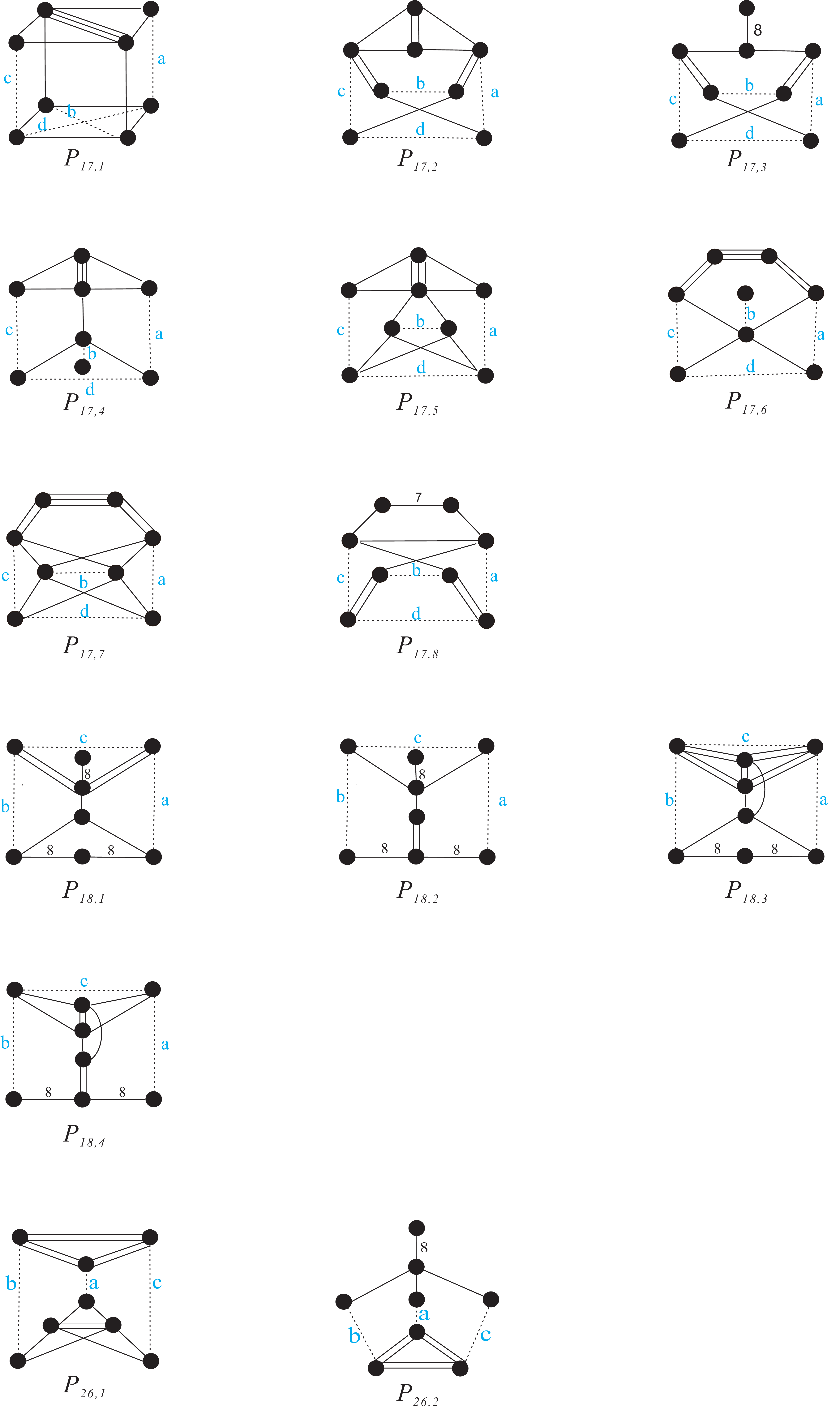}}
	\caption{$P_{17},P_{18},P_{26}$} \label{figure:p17182621}
\end{figure}

\newpage
\newgeometry{left=0.5cm,right=0.5cm,top=2cm,bottom=2cm}

\begin{table}[H]
	\resizebox*{17cm}{!}{
		\renewcommand{\arraystretch}{2.3}		
		\begin{tabular}{|c |c | c | c | c |} 
		\hline
       	& $a$ & $b$ & $c$ & $d$ \\                                      
       	\hline
       	$P_{17,1}$ & $\displaystyle\frac{1}{2} (1+\sqrt{5})$ & $\displaystyle\frac{1}{4} (3+\sqrt{5})$& $\displaystyle\frac{1}{2} (1+\sqrt{5})$ & $\displaystyle\frac{1}{2} (1+\sqrt{5})$\\[7pt]                    
       	\hline
       	
       	$P_{17,2}$ & $\displaystyle 1+\frac{1}{ \sqrt{2}}$ & $\displaystyle\frac{1}{2}+\frac{1}{\sqrt{2}}$   & $\displaystyle 1+\frac{1}{ \sqrt{2}}$  & $\displaystyle\frac{3}{2}+\frac{1}{ \sqrt{2}}$ \\[7pt]                   
       	\hline
       	\renewcommand\arraystretch{2} 
       	$P_{17,3}$ & $\displaystyle 1+\frac{1}{\sqrt{2}}$ & $\displaystyle\frac{1}{2}+\frac{1}{ \sqrt{2}}$   & $\displaystyle 1+\frac{1}{ \sqrt{2}}$   & $\displaystyle\frac{3}{2}+\frac{1}{\sqrt{2}}$  \\[7pt]                    
       	\hline
       	$P_{17,4}$ & $\displaystyle\frac{1}{2} (1+\sqrt{5})$ &   $\displaystyle\frac{1}{2}\sqrt{3+\frac{3}{ \sqrt{5}}}$ & $\displaystyle\frac{1}{2} (1+\sqrt{5})$ & $\displaystyle\frac{1}{2} (1+\sqrt{5})$ \\[7pt]                    
       	\hline
       	$P_{17,5}$ & $\displaystyle\frac{1}{2} (1+\sqrt{5})$ & $\displaystyle\frac{1}{10}(5+3\sqrt{5}) $& $\displaystyle\frac{1}{2} (1+\sqrt{5})$ & $\displaystyle\frac{1}{2} (1+\sqrt{5})$ \\[7pt]                  
       	\hline
       	$P_{17,6}$ & $\displaystyle\frac{1}{2} (1+\sqrt{5})$ & $\displaystyle \frac{1}{2} \sqrt{3+\sqrt{5}} $& $\displaystyle\frac{1}{2} (1+\sqrt{5})$ & $\displaystyle\frac{1}{2} (3+\sqrt{5})$\\[7pt]                   
       	\hline
       	$P_{17,7}$ & $\displaystyle\frac{1}{2} (1+\sqrt{5})$ & $\displaystyle\frac{1}{2} (1+\sqrt{5}) $& $\displaystyle\frac{1}{2} (1+\sqrt{5})$ & $\displaystyle\frac{1}{2} (3+\sqrt{5})$\\[7pt]                   
       	\hline
       	$P_{17,8}$ & $\displaystyle\frac{1}{\sqrt{2}}(2 \cos\frac{\pi}{7} (1 + 2 \cos\frac{\pi}{7})-1)$
       	& 
       	$\displaystyle 2 \cos^2\frac{\pi}{7}-\frac{1}{2}$ & 
       	$\displaystyle \frac{1}{\sqrt{2}}(2 \cos\frac{\pi}{7} (1 + 2 \cos\frac{\pi}{7})-1)$& 
       	$\displaystyle 2 \cos\frac{\pi}{7} (1 + 2 \cos\frac{\pi}{7})$
       	\\[7pt]                   
       	\hline
		\specialrule{0.5em}{0.5pt}{0.5pt} \multicolumn{5}{l}{}\\
		\hline	
		$P_{18,1}$ & $2+\sqrt{2}$ & $2+\sqrt{2}$ & $5+4\sqrt{2}$& \\[7pt]                    
		\hline
		
		$P_{18,2}$ & $\sqrt{2+\sqrt{2}}$ & $\sqrt{2+\sqrt{2}}$ & $1+\sqrt{2}$&  \\[7pt]                   
		\hline
		
		$P_{18,3}$ & $2+\sqrt{2}$ & $2+\sqrt{2}$ & $5+4\sqrt{2}$&  \\[7pt]       
		
		\hline
		$P_{18,4}$ & $\sqrt{2+\sqrt{2}}$ & $\sqrt{2+\sqrt{2}}$ & $1+\sqrt{2}$&
		\\[7pt]                   
		\hline
		\specialrule{0.5em}{0.5pt}{0.5pt} \multicolumn{5}{l}{}\\
		\hline
		$P_{26,1}$ &  $\sqrt{1+\frac{ 1}{ \sqrt{2}}}$ &    $\sqrt{1+\frac{ 1}{ \sqrt{2}}}$  &   $\sqrt{1+\frac{ 1}{ \sqrt{2}}}$ & \\[7pt]                    
		\hline
		$P_{26,2}$ & $\sqrt{1+\frac{ 1}{ \sqrt{2}}}$ &    $\sqrt{1+\frac{ 1}{ \sqrt{2}}}$  &   $\sqrt{1+\frac{ 1}{ \sqrt{2}}}$&
		 \\[7pt]                    
       	\hline
       	\specialrule{0.5em}{0.5pt}{0.5pt} \multicolumn{5}{l}{}\\
		\end{tabular}
		}
\end{table}

\restoregeometry

\begin{figure}[H]
	\scalebox{0.3}[0.3]{\includegraphics {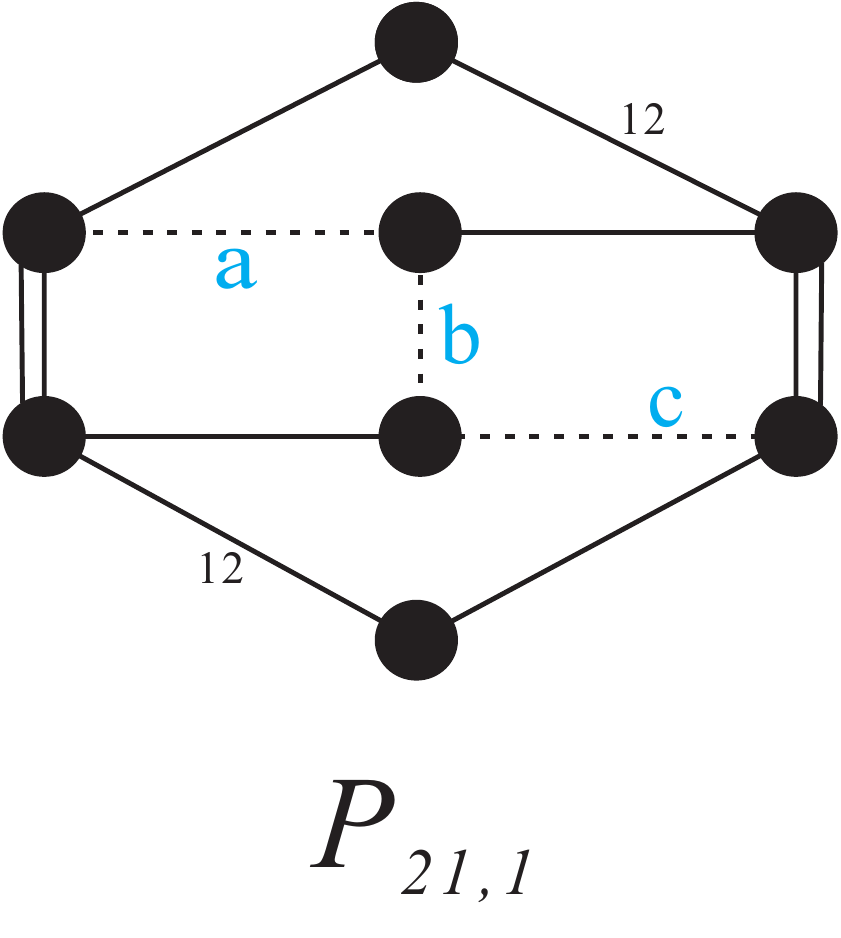}}
	\caption{$P_{21}$} \label{figure:p21}
\end{figure}
\restoregeometry

\newpage
\newgeometry{left=0.5cm,right=0.5cm,top=2cm,bottom=2cm}
\begin{table}[H]
	\resizebox*{16cm}{!}{
		\renewcommand{\arraystretch}{2.7}		
		\begin{tabular}{|c |c | c | c | c |} 
			\hline
			& $a$ & $b$ & $c$ & $d$ \\                                      
			\hline
			$P_{34,1}$ & $  \frac{1}{4}(1+\sqrt{13})$ &   $   \frac{1}{4}(1+\sqrt{13})$  &  $\sqrt{\frac{1}{6}{(5+\sqrt{13})}}$ & $ \frac{1}{4}(1+\sqrt{13})$ \\[7pt]      
			\hline					
			
			$P_{34,2}$ & $\frac{1}{2}\sqrt{3+\sqrt{5}}$ & $\frac{1}{2}\sqrt{3+\sqrt{5}}$ & 
			$\frac{1}{2}\sqrt{3+\sqrt{5}}$ & $\frac{1}{2}\sqrt{3+\sqrt{5}}$ \\[7pt]          
			\hline
			
			$P_{34,3}$ & $  \frac{\sqrt{5}}{2}$ & $   \sqrt{\frac{1}{10}(15-\sqrt{5})}$ & $    \frac{\sqrt{5}}{2}$ &$    \frac{\sqrt{5}}{2}$  \\[7pt]           					
			\hline
			
			$P_{34,4}$& $  \frac{1}{2}+\frac{1}{\sqrt{2}}$ & $   \frac{1}{2}+\frac{1}{\sqrt{2}}$ & $    \frac{1}{2}+\frac{1}{\sqrt{2}}$ & $    \frac{1}{2}+\frac{1}{\sqrt{2}}$ \\[7pt]            					
			\hline
			
			$P_{34,5}$ &  $  \frac{1}{2}(1+\sqrt{5})$ &  $   \frac{1}{2}\sqrt{3+\sqrt{5}}$ &  $    \frac{1}{2}\sqrt{3+\sqrt{5}}$  & $    \frac{1}{2}\sqrt{3+\sqrt{5}} $   \\[7pt]     					
			\hline
			
			$P_{34,6}$ &  $  \frac{1}{2}\sqrt{3+\sqrt{5}} $ &   $     \frac{1}{2}(1+\sqrt{5}) $  & $     \frac{1}{2}(1+\sqrt{5}) $ & $     \frac{1}{2}\sqrt{3+\sqrt{5}}$ \\[7pt]     					
			\hline
			
			$P_{34,7}$ &  $  \frac{1}{4}(1+\sqrt{13})$ &  $   \frac{1}{4}(1+\sqrt{13}) $ &  $    \frac{1}{4}(5+\sqrt{13})$  & $    \frac{1}{4}(1+\sqrt{13})$ \\[7pt]         				
			\hline
			
			$P_{34,8}$ &  $  \frac{1}{2}\sqrt{3+\sqrt{5}}$ &   $   \frac{1}{2}\sqrt{3+\sqrt{5}}$  &  $    \frac{1}{4}(3+\sqrt{5}+\sqrt{3+\sqrt{5}}$ & $   \frac{1}{2}\sqrt{3+\sqrt{5}}$   \\[7pt]                    
			\hline
			
			$P_{34,9}$&  $  \frac{\sqrt{5}}{2} $ &   $  \frac{\sqrt{5}}{2}$  &  $    \frac{1}{4}(5+\sqrt{5})$ & $    \frac{\sqrt{5}}{2}$ 
			\\[7pt]       				
			\hline
			
			$P_{34,10}$ &  $  \frac{1}{4}(1+\sqrt{13}) $ &   $    \frac{1}{4}(1+\sqrt{13}) $  & $    \frac{1}{3}(2+\sqrt{13})$ & $      \frac{1}{4}(1+\sqrt{13})$\\[7pt]         
			\hline
			
			$P_{34,11}$ & $  \frac{1}{2}\sqrt{3+\sqrt{5}}$ &   $  \frac{1}{2}(1+\sqrt{5})$  &  $    \frac{1}{2}\sqrt{3+\sqrt{5}}$ & $     \frac{1}{2}\sqrt{3+\sqrt{5}}$  \\[7pt]    			
			\hline
			
			$P_{34,12}$& $  \frac{\sqrt{5}}{2}$ & $   2-\frac{1}{\sqrt{5}}$ & $    \frac{\sqrt{5}}{2}$ & $    \frac{\sqrt{5}}{2}$  \\[7pt]        					
			\hline                   
			\specialrule{0.5em}{0.5pt}{0.5pt} \multicolumn{5}{l}{}\\
			\hline	
			$P_{21,1}$ &  $-\frac{ 1}{2} (2 \sqrt{3+\sqrt{3}}-(2+\sqrt{3})^{3/2})$ & $\sqrt{2(2+\sqrt{3})}$ &  $-\frac{ 1}{2} (2 \sqrt{3+\sqrt{3}}-(2+\sqrt{3})^{3/2})$& \\[7pt]\hline
			\specialrule{0.5em}{0.5pt}{0.5pt} \multicolumn{5}{l}{}\\
		\end{tabular}
	}
\end{table}
\restoregeometry


\begin{thebibliography}{10}

\bibitem[All06]{Allcock: 2006}
D. Allcock.
\newblock  Infinitely many hyperbolic Coxeter groups through dimension 19
\newblock  Geometry $\&$ Topology 10 (2006), 737--758.



\bibitem[Bur22]{Amanda:2022}
A. Burcroff.
\newblock Near classification of compact hyperbolic Coxeter d-polytopes with $d + 4$ facets and related dimension bounds.
\newblock arXiv:2201.03437.

\bibitem[And$70^{(1)}$]{Andreev1: 1970}
E. M. Andreev.
\newblock Convex polyhedra in Lobachevskii spaces (Russian).
\newblock Math. USSR Sbornik 10 (1970), 413--440.

\bibitem[And$70^{(2)}$]{Andreev2: 1970}
E. M. Andreev.
\newblock Convex polyhedra of finite volume in Lobachevskii space (Russian).
\newblock Math. USSR Sbornik 12 (1970), 255--259.




\bibitem[Bor98]{Borcherds: 1998}
R. Borcherds.
\newblock Coxeter groups, Lorentzian lattices, and K3 surfaces.
\newblock IMRN, 19 (1998), 1011--1031.


\bibitem[Bou68]{Bourbaki: 1968}
N. Bourbaki.
\newblock Groupes et alg`ebres de Lie.
\newblock Ch. IV--VI, Hermann, Paris, 1968.



\bibitem[Bug84]{Bugaenko:1984}
V. O. Bugaenko.
\newblock Groups of automorphisms of unimodular hyperbolic quadratic forms over the ring $\mathbb{Z}[\frac{\sqrt{5}+1}{2}]$
\newblock Moscow Univ. Math. Bull., 39 (1984), 6--14.

\bibitem[Bug92]{Bugaenko:1992}
V. O. Bugaenko.
\newblock Arithmetic crystallographic groups generated by reflections, and reflective hyperbolic lattices.
\newblock Adv. Sov. Math., 8 (1992), 33--55.

\bibitem[Cox34]{Coxeter:1934}	
H. S. M. Coxeter.
\newblock   Discrete groups generated by reflections.
\newblock   Ann. Math., 35 (1934), 588--621.

\bibitem[Ess94]{Esselmann}
F. Esselmann.
\newblock {\"{U}ber} kompakte hyperbolische Coxeter-Polytope mit wenigen Facetten.
\newblock Universit at Bielefeld, SFB 343, (1994) No. 94-087.


\bibitem[Ess96]{Esselmann: 1996}
F. Esselmann.
\newblock The classification of compact hyperbolic Coxeter $d$-polytopes with $d+2$ facets.
\newblock Comment. Math. Helvetici 71 (1996), 229--242.

\bibitem[F]{Annahomepage}
A. Flikson
\newblock http://www.maths.dur.ac.uk/users/anna.felikson/Polytopes/polytopes.html



\bibitem[FT$08^{(1)}$]{FT:08}
A. Felikson, P. Tumarkin.
\newblock On hyperbolic Coxeter polytopes with mutually intersecting facets.
\newblock J. Combin. Theory A 115 (2008), 121-146.


\bibitem[FT$08^{(2)}$]{FT:08s}
A. Felikson, P. Tumarkin.
\newblock On compact hyperbolic Coxeter $d$-polytopes with $d+4$ facets.
\newblock Trans. Moscow Math. Soc. 69 (2008), 105-151.

\bibitem[FT09]{FT:09}
A. Felikson, P. Tumarkin.
\newblock Coxeter polytopes with a unique pair of non-intersecting facets.
\newblock J. Combin. Theory A 116 (2009), 875-902.

\bibitem[FT14]{FT:14}
A. Felikson, P. Tumarkin.
\newblock Essential hyperbolic Coxeter polytopes.
\newblock Israel Journal of Mathematics (2014), 199, 1, 113–161.


\bibitem[Gan59]{G:1959}
F. R. Gantmacher.
\newblock  The theory of matrix.
\newblock  Chelsca Publishing Company (1959), New York.

\bibitem[Gr\"{u}67]{G:1967}
B. Gr\"{u}nbaum.
\newblock  Convex Polytopes.
\newblock Wiley, New York, (1967).2nd edn. Springer, Berlin (2003)

\bibitem[Gr\"{u}72]{G:1972}
B. Gr\"{u}nbaum.
\newblock Arrangements and Spreads. 
\newblock American Mathematical Society, New York (1972). 


\bibitem[GS67]{GS:1967}
B. Gr\"{u}nbaum, V.P. Sreedharan.
\newblock An enumeration of simplicial 4-polytopes with 8 vertices.
\newblock J. Combin. Theory 2 (1967), 437–465.




\bibitem[ImH90]{ImH: 1990}
H.-C. Im Hof.
\newblock Napier cycles and hyperbolic Coxeter groups, 
\newblock  Bull. Soc. Math. de Belg. S´erie A, XLII (1990), 523--545.


\bibitem [Jac17]{Jacquemet2017}
M. Jacquemet
\newblock On hyperbolic Coxeter n-cubes.
\newblock Europ. J. Combin. 59 (2017), 192-203.

\bibitem [JT18]{JT:2018}
 M. Jacquemet and S. Tschantz.
 \newblock All hyperbolic Coxeter n-cubes
 \newblock J. Combin. Theory A, 158 (2018), 387-406.



\bibitem[Kap74]{Kaplinskaya: 1974}
I. M. Kaplinskaya. 
\newblock Discrete groups generated by reflections in the faces
of simplicial prisms in Lobachevskian spaces.
\newblock  Math. Notes 15 (1974), 88--91.

\bibitem[KM13]{KM: 2013}
A. Kolpakov and B. Martelli.
\newblock Hyperbolic four-manifolds with one cusp.
\newblock Geom $\&$ Funct. Anal., 23, 2013, 1903--1933.

\bibitem[Kos67]{Koszul:1967}
J. L. Koszul.
\newblock Lectures on hyperbolic Coxeter group.
\newblock University of Notre Dame, 1967.

\bibitem[Lan50]{Lanner: 1950}
F. Lanner.
\newblock On complexes with transitive groups of automorphisms.
\newblock Comm. Sem. Math. Univ. Lund, 11 (1950), 1--71.


\bibitem[MZ18]{MZ:2018}
J. Ma and  F. Zheng.
\newblock  Orientable hyperbolic 4-manifolds ove the 120-cell.
\newblock ArXiv:math:GT/1801.08814.

\bibitem[Mak65]{Makarov: 1965}
V. S. Makarov.
\newblock On one class of partitions of Lobachevskian space.
\newblock Dokl. Akad. Nauk SSSR, 161 no.2, (1965), 277--278. English transl.: Sov. Math. Dokl. 6, 400-401 (1965), Zbl.135,209.

\bibitem[Mak66]{Makarov: 1966}
V. S. Makarov.
\newblock On one class of discrete groups of Lobachevskian space having an infinite fundamental region of finite measure.
\newblock Dokl. Akad. Nauk SSSR, 167 no.2, (1966), 30--33. English transl.: Sov. Math. Dokl. 7, 328-331 (1966), Zbl.146,165.

\bibitem[Mak68]{Makarov: 1968}
V. S. Makarov.
\newblock On Fedorov groups of four- and five-dimensional Lobachevsky spaces.
\newblock Issled. po obshch. algebre. no. 1, Kishinev St. Univ., 1970, 120--129 (Russian).




\bibitem[Mn\"{e}88]{Mne:1988}
N. E Mnëv. 
\newblock The universality theorems on the classification problem of configuration
varieties and convex polytopes varieties.
\newblock volume 1346 of Lecture Notes in
Mathematics, pages 527--543. Springer, 1988.

\bibitem[P1882]{Poincare: 1882}
H. Poincar\'{e}.
\newblock Theori\'{e} des groupes fuchsiens (French).
\newblock Acta Math. 1, 1882, no. 1, 1--76.	

\bibitem[Pro87]{Prokhorov:1987}	
M. N. Prokhorov.
\newblock The absence of discrete reflection groups with non-compact fundamental polyhedron of finite volume in Lobachevsky space of large dimension.
 \newblock Math. USSR Izv., 28 (1987), 401--411.


	

\bibitem[Rob15]{Roberts:15} 
M. Roberts.
\newblock A Classification of Non-Compact Coxeter Polytopes with $n+3$ Facets and One Non-Simple Vertex
\newblock arXiv:1511.08451.

\bibitem[Rus89]{Rusmanov: 1989}
O. P. Rusmanov.
\newblock  Examples of non-arithmetic crystallographic Coxeter groups in $n$-dimensional Lobachevsky space for $6\leq n\leq 10$. 
\newblock Problems in Group Theory and Homology
Algebra, Yaroslavl, 1989, 138--142 (Russian).

\bibitem[Tum03]{Tumarkin: n3fvs}
P. Tumarkin.
\newblock Non-compact hyperbolic Coxeter n-polytopes with n+3 facets, short version (3 pages).
\newblock  Russian Math. Surveys, 58 (2003), 805-806.

\bibitem[Tuma$04^{(1)}$]{Tumarkin: n2}
P. Tumarkin.
\newblock Hyperbolic Coxeter $n$-polytopes with $n+2$ facets. 
\newblock Math. Notes 75 (2004), 848--854.

\bibitem[Tum$04^{(2)}$]{Tumarkin: n3fv}
P. Tumarkin.
\newblock Hyperbolic Coxeter n-polytopes with $n+3$ facets.
\newblock Trans. Moscow Math. Soc. (2004), 235–25.

\bibitem[Tum07]{Tumarkin: n3}
P. Tumarkin.
\newblock Compact hyperbolic Coxeter n-polytopes with $n + 3$ facets.
\newblock the electronic journal of combinatorics 14 (2007).


\bibitem[Vin67]{Vinberg: 1967}
E. B. Vinberg.
\newblock Discrete groups generated by reflections in Lobachevskii
spaces. 
\newblock Mat. USSR Sb. 1 (1967), 429--444.

\bibitem[Vin69]{Vinberg: 1969}
E. B. Vinberg.
\newblock Some examples of crystallographic groups in Lobachevskian spaces. 
\newblock 78(120), no. 4,(1969), 633--639.

\bibitem[Vin72]{Vinberg:1972}
E. B. Vinberg.
\newblock On groups of unit elements of certain quadratic forms.
\newblock Math. USSR Sb., 16 (1972), 17--35.

\bibitem[Vin$85^{(1)}$]{Vinberg:1985}	
E. B. Vinberg. 
\newblock   The absence of crystallographic groups of reflections in Lobachevsky spaces of large dimension.
\newblock    Trans. Moscow Math. Soc., 47 (1985), 75--112.

\bibitem[Vin$85^{(2)}$]{Vinberg:1985s}
E. B. Vinberg.
\newblock Hyperbolic reflection groups. 
\newblock Russian Math. Surveys 40 (1985), 31--75.

\bibitem[Vin93]{Vinberg:1993}
E. B. Vinberg, editor.
\newblock  Geometry II - Spaces of Constant Curvature.
\newblock Springer, 1993.


\bibitem[VK78]{VK:1978}
E. B. Vinberg and I. M. Kaplinskaya.
\newblock On the groups $O_{18,1}(\mathbb{Z})$ and $O_{19,1}(\mathbb{Z})$.
\newblock Soviet Math. Dokl. 19 (1978), 194--197.






\end{thebibliography}
\end{document}